\documentclass[a4paper,12pt]{amsart}
\usepackage{amssymb,amscd,amsfonts,amsmath,times}
\usepackage{mathrsfs}
\usepackage[all]{xy}
\numberwithin{equation}{subsection}
     \newtheorem{thm}{Theorem}[subsection]
     \newtheorem{lem}[thm]{Lemma}
     \newtheorem{keylem}[thm]{Key Lemma}
     \newtheorem{sublem}[thm]{Sublemma}
     \newtheorem{prop}[thm]{Proposition}
     \newtheorem{cor}[thm]{Corollary}
     
     \newtheorem{conj}[thm]{Conjecture}
     
     \newtheorem{defn}[thm]{Definition}

     \newtheorem{rem}[thm]{Remark}
     \newtheorem{mthm}{Theorem}[section]
     \newtheorem{mprop}[mthm]{Proposition}
     \newtheorem{mcor}[mthm]{Corollary}
     \newtheorem{mfact}[mthm]{Fact}

\newenvironment{pf}{\par\smallskip\noindent{\it Proof.}}{\hfill$\square$\par\medskip}
\newenvironment{pf*}[1]{\par\smallskip\noindent{\bf #1. }}{\hfill$\square$\par\medskip}
\begin{document}
\title[$p$-adic regulator and finiteness]{A $\boldsymbol{p}$-adic regulator map and finiteness results \endgraf
 for arithmetic schemes$^{\,1}$}
\author[S. Saito \and K. Sato]{Shuji Saito \and Kanetomo Sato}
\address[Shuji Saito]{\;Department of Mathematics, University of Tokyo \endgraf
 Komaba, Meguro-ku, Tokyo 153-8914, JAPAN}
\email{sshuji@ms.u-tokyo.ac.jp}
\address[Kanetomo Sato]{\;Graduate School of Mathematics, Nagoya University \endgraf
 Furo-cho, Chikusa-ku, Nagoya 464-8602, JAPAN}
\email{kanetomo@math.nagoya-u.ac.jp}
\date{January 21, 2009}
\thanks{
\endgraf
(1) The earlier version was entitled `Torsion cycle class maps in codimension two of arithmetic schemes'.}
\begin{abstract}A main theme of the paper is a conjecture of Bloch-Kato on the image of $p$-adic regulator maps for a proper smooth variety $X$ over an algebraic number field $k$. The conjecture for a regulator map of particular degree and weight is related to finiteness of two arithmetic objects:
One is the $p$-primary torsion part of the Chow group in codimension $2$ of $X$.
Another is an unramified cohomology group of $X$.
As an application, for a regular model ${\mathscr X}$ of $X$ over the integer ring of $k$,
we show an injectivity result on torsion of a cycle class map from the Chow group in codimension $2$ of ${\mathscr X}$ to a new $p$-adic cohomology of ${\mathscr X}$ introduced by the second author, which is a candidate of the conjectural \'etale motivic cohomology with finite coefficients of Beilinson-Lichtenbaum.
\end{abstract}
\maketitle

\def\Br{{\mathrm{Br}}}      
\def\cd{{\mathrm{cd}}}      
\def\codim{{\mathrm{codim}}}
\def\ch{{\mathrm{ch}}}      
\def\CH{{\mathrm{CH}}}      
\def\Coker{{\mathrm{Coker}}}
\def\Ker{{\mathrm{Ker}}}    
\def\cont{{\mathrm{cont}}}  
\def\Cor{{\mathrm{cor}}}    
\def\Div{{\mathrm{Div}}}    
\def\e{{\epsilon}}          %
\def\et{{\mathrm{\acute{e}t}}}    
\def\Ext{{\mathrm{Ext}}}    
\def\F{F}        
\def\FF{{\mathscr F}}           
\def\Frac{{\mathrm{Frac}}}  
\def\Gal{{\mathrm{Gal}}}    
\def\gr{{\mathrm{gr}}}      
\def\gys{{\mathrm{Gys}}}    
\def\H{H}        
\def\Hom{{\mathrm{Hom}}}    
\def\id{{\mathrm{id}}}     
\def\Image{{\mathrm{Im}}}   
\def\ind{{\mathrm{ind}}}    
\def\K{{\mathscr K}}            
\def\ker{{\mathrm{Ker}}}    
\def\L{{\mathbb L}}            %
\def\LL{{\mathscr L}}            %
\def\loc{{\mathrm{loc}}}    
\def\NF{N}       
\def\N{{\mathbb N}}            
\def\Nis{\mathrm{Nis}}      
\def\NS{{\mathrm{NS}}}      
\def\num{{\mathrm{num}}}    
\def\O{{\mathscr O}}            
\def\OO{{\mathfrak o}}          
\def\ord{{\mathrm{ord}}}    
\def\pDiv{p\text{-}\Div}    
\def\P{{\mathbb P}}            
\def\Pic{{\mathrm{Pic}}}    
\def\ptor{{p\text{-}\mathrm{tors}}}
\def\Q{{\mathbb Q}}            
\def\reg{{\mathrm{reg}}}    
\def\res{{\mathrm{res}}}    
\def\rf{{\mathbb F}}           
\def\Spec{{\mathrm{Spec}}}  
\def\corank{{\mathrm{corank}}}    
\def\T{{\mathfrak T}}
\def\tor{{\mathrm{tors}}}   
\def\Cotor{{\mathrm{Cotor}}}   
\def\ur{{\mathrm{ur}}}      
\def\val{{\mathrm{val}}}    
\def\WW{{\mathscr W}}
\def\cX{{\mathscr X}}
\def\Z{{\mathbb Z}}            
\def\zar{{\mathrm{Zar}}}    

\def\cO{{\mathscr O}}
\def\cR{{\mathscr R}}
\def\cS{{\mathscr S}}
\def\btil{\wt{\phantom{a}}}
\def\ctil{\hspace{-4pt}\wt{\phantom{a}}}
\def\dtil{\hspace{-10pt}\wt{\phantom{a}}\hspace{3pt}}

%
\def\Lam{\varLambda}
\def\lam{\lambda}
\def\Te{\varTheta}
\def\te{\theta}
\def\vG{\varGamma}
\def\k{\kappa}
%
%
%
%
\def\ra{\rightarrow}
\def\lra{\longrightarrow}
\def\Lra{\Longrightarrow}
\def\lla{\longleftarrow}
\def\hra{\hookrightarrow}
\def\lmt{\longmapsto}
\def\sm{\setminus}
\def\wt#1{\widetilde{#1}}
\def\wh#1{\widehat{#1}}
\def\bs#1{\boldsymbol{#1}}
\def\ol#1{\overline{#1}}
\def\ul#1{\underline{#1}}
\def\us#1#2{\underset{#1}{#2}}
\def\os#1#2{\overset{#1}{#2}}
%
%
\def\isom{\hspace{9pt}{}^{\sim}\hspace{-16.5pt}\lra}
\def\lisom{\hspace{9pt}{}^{\sim}\hspace{-17.5pt}\lla}
%
%
%
\def\witt#1#2#3{W_{\hspace{-2pt}#2}{\hspace{1pt}}\Omega_{#1}^{#3}}
\def\mlogwitt#1#2#3{W_{\hspace{-2pt}#2}{\hspace{1pt}}\omega_{{#1},{\log}}^{#3}}
\def\logwitt#1#2#3{W_{\hspace{-2pt}#2}{\hspace{1pt}}\Omega_{{#1},{\log}}^{#3}}
%
%
\def\Gm{{\mathbb G}_{\hspace{-1pt}\mathrm{m}}}
\def\bF{{\mathbb F}}
\def\bR{{\mathbb R}}
\def\bZ{{\mathbb Z}}
\def\bQ{{\mathbb Q}}
\def\QpZp{{\Q}_p/{\Z}_p}
\def\qz{{\bQ}/{\bZ}}
\def\qzp{{\bQ_p}/{\bZ_p}}
\def\zl{{\bZ}_\ell}
\def\ql{{\bQ}_\ell}
\def\zp{{\bZ}_p}
\def\qp{{\bQ}_p}
\def\Qp{{\Q}_p}
\def\pnuz{\bZ/p^r\bZ}
\def\prz{\bZ/p^r\bZ}

\def\XU{\cX_U}
\def\XUp{\cX_U[p^{-1}]}
\def\Xgf{X}
\def\Xggf{\ol X}
\def\Xff{K}
\def\Zggf{\ol Z}

\def\knr{k^{\ur}}
\def\Xnr{\cX^{\ur}}
\def\Xknr{X^{\ur}}
\def\Ynr{\ol Y}
\def\kv{k_v}
\def\Fv{\rf_v}
\def\jv{j_v}
\def\iv{i_v}
\def\Xv{\cX_v}
\def\Yv{Y_v}
\def\Ynrv{\ol {\Yv}}
\def\Xkv{X_v}
\def\mupr#1{\mu_{p^r}^{\otimes #1}}
\def\mupnu#1{\mu_{p^r}^{\otimes #1}}
\def\adiv{\omega}
\def\bdiv{\beta_{\Div}}
\def\cdiv{\nu}
\def\regXkLam{\reg_{\Lam}}
\def\regXkqzp{\reg_{\qzp}}
\def\regXULam{\reg_{\cX_U,\Lam}}

\def\pnurX{{\mathfrak T}_r(n)}
\def\psnX{{\mathfrak T}_s(n)}
\def\pnmurX{{\mathfrak T}_{r+s}(n)}
\def\pnuX{{\mathfrak T}_r(2)}
\def\qzpX{{\mathfrak T}_{\infty}(2)}
\def\zptX{{\mathfrak T}_{\zp}(2)}
\def\qptX{{\mathfrak T}_{\qp}(2)}
\def\qzpnX{{\mathfrak T}_{\infty}(n)}
\def\zpnX{{\mathfrak T}_{\zp}(n)}
\def\qpnX{{\mathfrak T}_{\qp}(n)}

\def\dval{\partial^{\val}}
\def\dloc{\delta^{\loc}}
\def\db{\ol{{\mathfrak d}}}
\def\deltad{{\mathfrak d}}
\def\deltadd{\delta}
\def\tdelta{{\mathfrak d}_2}
\def\xdelta{{\mathfrak d}_3}
\def\tS{S'}
\def\tk{k'}
\def\tXggf{\ol X{}'}
\def\tX{\cX'}
\def\tY{Y'}
\def\tZ{\widetilde{Z}}

\def\xt{\ol {x}}

\def\qaq{\quad\text{and}\quad}

\def\Fb{\ol {\rf}}
\def\Xb{\ol X}
\def\Yb{\ol Y}
\def\zb{\ul z}

\def\Faa{\F^1\hspace{-1.5pt}}
\def\Fbb{\F^2\hspace{-1.3pt}}
\def\Naa{\NF^1\hspace{-1.5pt}}
\def\Nbb{\NF^2\hspace{-1.3pt}}
\def\NFF{N\hspace{-1.2pt}F}
\def\NFaa{\NFF^1\hspace{-1.5pt}}

\def\Hg#1#2{\H^1_{f,#1}(k,#2)}
\def\Hind#1#2{\H^1_{\ind}(#1,#2)}

\def\grF#1{\gr_{\F}^{#1}}
\def\varThetaur{\wt{\varTheta}}

\def\HnrXkk#1#2{\H_{\ur}^{#1}(\Xgf,\qzp(#2))}
\def\HnrXk{\H_{\ur}^{3}(\Xgf,\qzp(2))}
\def\HnrXq{\H^{n+1}_{\ur}(\Xff,\qzp(n))}
\def\HnrXz{\H^1_{\ur}(\Xff,\qzp(0))}
\def\HnrXo{\H^2_{\ur}(\Xff,\qzp(1))}
\def\HnrXt{\H^3_{\ur}(\Xff,\qzp(2))}
\def\HnrXd{\H^{d+1}_{\ur}(\Xff,\qzp(d))}
\def\HBX{\H^3_{\ur}(\Xff,\Xgf;\qzp(2))}

\def\HurXb{\H_{\ur}^{2}(K,\qz(1))}
\def\HurXc{\H_{\ur}^{3}(K,\qz(2))}
\def\bhZ{\widehat{\bZ}}

\def\GN{9}
\def\LN{10}

%
\section{Introduction}\label{intro}
\medskip
Let $k$ be an algebraic number field
and let $G_k$ be the absolute Galois group $\Gal(\ol k/k)$,
 where $\ol k$ denotes a fixed algebraic closure of $k$.
Let $X$ be a projective smooth variety over $k$ and put $\ol X:=X \otimes_k \overline{k}$.
Fix a prime $p$ and integers $r,m\geq 1$.
A main theme of this paper is a conjecture of Bloch and Kato
 concerning the image of the $p$-adic regulator map
\[ \reg^{r,m} \;:\; \CH^r(X,m)\otimes\Qp \lra \H^1_\cont(k,\H^{2r-m-1}_\et(\ol X,\Qp(r))) \]
from Bloch's higher Chow group to continuous Galois cohomology of $G_k$
(\cite{BK2}, Conjecture 5.3).
See $\S\ref{sect2}$ below for the definition of  this map in the case $(r,m)=(2,1)$.
This conjecture affirms that its image agrees with the subspace
$\H^1_g(k,\H^{2r-m-1}_\et(\ol X,\Qp(r)))$ defined in loc.\ cit.,
 and plays a crucial role in the so-called Tamagawa number conjecture 
on special values of $L$-functions attached to $\Xgf$.
In terms of Galois representations, the conjecture means that a 
$1$-extension of continuous $p$-adic representations of $G_k$
\[ 0 \lra \H^{2r-m-1}_\et(\ol X,\Qp(r)) \lra E \lra \Qp \lra 0 \]
arises from a $1$-extension of motives over $k$
\[ 0\lra h^{2r-m-1}(X)(r) \lra M \lra h(\Spec(k)) \lra 0, \]
 if and only if $E$ is a de Rham representation of $G_k$. 
There has been only very few known results on the conjecture.
In this paper we consider the following condition,
 which is the Bloch-Kato conjecture in the special case $(r,m)=(2,1)$:
\par
\vspace{8pt}
\begin{enumerate}
\item[]
{\bf H1:}
{\it The image of the regulator map
\[ \reg:=\reg^{2,1} : \CH^2(\Xgf,1) \otimes \qp \lra \H^1_{\cont}(k,\H^2_{\et}(\Xggf,\qp(2))). \]
agrees with $\H^1_g(k,\H^2_{\et}(\Xggf,\qp(2)))$.}
\end{enumerate}
\par
\vspace{8pt}
\noindent
We also consider a variant:
\par
\vspace{8pt}
\noindent
\begin{enumerate}
\item[]
{\bf H1*:}
{\it The image of the regulator map with $\qzp$-coefficients
\[ \regXkqzp : \CH^2(\Xgf,1) \otimes \qzp \lra
   \H^1_{\Gal}(k,\H^2_{\et}(\Xggf,\qzp(2))) \]
agrees with $\H^1_g(k,\H^2_{\et}(\Xggf,\qzp(2)))_{\Div}$
    $($see $\S\ref{sect1.1}$ for the definition of this group$)$}.
\end{enumerate}
\par
\vspace{8pt}
\noindent
We will show that {\bf H1} always implies {\bf H1*}, which is not straight-forward.
On the other hand the converse holds as well
  under some assumptions.
See Remark \ref{lem2-01} below for details.
\begin{mfact}\label{fact0-0}
The condition {\bf H1} holds in the following cases$:$
\begin{itemize}
\item[(1)]
$\H^2(\Xgf,\O_{\Xgf})=0$ $($\cite{CTR1}, \cite{CTR2}, \cite{Sal}$)$.
\item[(2)]
$\Xgf$ is the self-product of an elliptic curve over $k=\bQ$
with square-free conductor and
without complex multiplication,
   and $p \geq 5$ $($\cite{Md}, \cite{Fl}, \cite{LS}, \cite{La0}$)$.
\item[(3)]
$\Xgf$ is the elliptic modular surface of level $4$ over $k=\bQ$
   and $p \geq 5$ $($\cite{La1}$)$.
\item[(4)]
$\Xgf$ is a Fermat quartic surface over $k=\bQ$ or $\bQ(\sqrt{-1})$
   and $p \geq 5$ $($\cite{O}$)$.
\end{itemize}
\end{mfact}
\noindent
A main result of this paper relates the condition {\bf H1*} to finiteness of two arithmetic objects.
One is the $p$-primary torsion part of the Chow group $\CH^2(\Xgf)$ of algebraic cycles of 
codimension two on $\Xgf$ modulo rational equivalence.
Another is an unramified cohomology of $\Xgf$, which we are going to introduce in what follows.
\par
\bigskip
Let $\OO_k$ be the integer ring of $k$, and put $S:=\Spec(\OO_k)$.
We assume that there exists a regular scheme $\cX$ which is proper flat of finite type over $S$
 and whose generic fiber is $\Xgf$. We also assume the following:
\begin{enumerate}
\item[$(*)$]
{\it $\cX$ has good or semistable reduction at each closed point of $S$ of characteristic $p$.}
\end{enumerate}
Let $\Xff=k(X)$ be the function field of $\Xgf$.
For an integer $q\geq 0$, let $\cX^q$ be the set of all points
     $x \in \cX$ of codimension $q$.
Fix an integer $n \ge 0$. Roughly speaking, the unramified cohomology group
$\HnrXq$ is defined as the subgroup of $\H^{n+1}_{\et}(\Spec(\Xff),\qzp(n))$ 
consisting of those elements that are ``unramified" along all $y\in \cX^1$.
For a precise definition, we need the $p$-adic \'etale Tate twist $\T_r(n)=\T_r(n)_{\cX}$ introduced in \cite{Sat2}.
This object is defined in $D^b(\cX_{\et},\pnuz)$,
the derived category of bounded complexes of
     \'etale sheaves of $\pnuz$-modules on $\cX$,
 and expected to coincide with
   $\varGamma(2)_{\et}^\cX\otimes^{\mathbb L} \pnuz$.
 Here $\varGamma(2)_{\et}^\cX$ denotes the conjectural
   \'etale motivic complex of Beilinson-Lichtenbaum \cite{Be}, \cite{Li}.
We note that the restriction of $\T_r(n)$ to $\cX[p^{-1}]:=\cX\otimes_{\bZ}\bZ[p^{-1}]$
 is isomorphic to $\mu_{p^r}^{\otimes n}$,
where $\mu_{p^r}$ denotes the \'etale sheaf of $p^r\mathrm{th}$ roots of unity.
Then $\HnrXq$ is defined as the kernel of the boundary map of \'etale cohomology groups
\[ \H^{n+1}_{\et}(\Spec(\Xff),\qzp(n)) \lra \bigoplus_{x \in \cX^1} \
   \H^{n+2}_x(\Spec(\O_{\hspace{-2pt}\cX\hspace{-2pt},x}),\T_{\infty}(n)), \]
where $\T_{\infty}(n)$ denotes $\varinjlim{}_{r \geq 1}\ \T_r(n)$.
There are natural isomorphisms
\[ \HnrXz \simeq \H^1_{\et}(\cX\hspace{-2pt},\qzp) \qaq \HnrXo \simeq \Br(\cX)_{\ptor}, \]
where $\Br(\cX)$ denotes the Grothendieck-Brauer group $\H^2_{\et}(\cX\hspace{-2pt},\Gm)$,
 and for an abelian group $M$,
   $M_{\ptor}$ denotes its $p$-primary torsion part.
An intriguing question is as to whether the group $\HnrXq$ is finite,
which is related to several significant theorems and
conjectures in arithmetic geometry (see Remark \ref{rem:ur} below).
In this paper we are concerned with the case $n=2$.
A crucial role will be played by the following subgroup
  of $\HnrXt$:
\begin{align*}
  & \HBX  \\
  & := \Image\Big(  \H^3_{\et}(\Xgf,\qzp(2)) \ra
      \H^3_{\et}(\Spec(\Xff),\qzp(2))\Big) \cap \HnrXt.
\end{align*}
Our finiteness result is the following:
\begin{mthm}\label{thm0-0}
Let $X$ and $\cX$ be as above, and
assume $p\geq 5$. Then$:$
\begin{enumerate}
\item[(1)]
 {\bf H1*} implies that $\CH^2(\Xgf)_{\ptor}$ and $\HBX$ are finite.
\item[(2)]
Assume that the reduced part of every closed fiber of $\cX/S$ has simple normal crossings on $\cX$,
  and that the Tate conjecture holds in codimension $1$
       for the irreducible components of those fibers.
Then the finiteness of the groups $\CH^2(\Xgf)_{\ptor}$ and
\linebreak $\HBX$ implies {\bf H1*}.
\end{enumerate}
\end{mthm}
\noindent
The assertion (2) is a converse of (1) under the assumption of the Tate conjecture.
We obtain the following result from Theorem \ref{thm0-0}\,(1)
 (see also the proof of Theorem \ref{cor0-1} in \S\ref{sect5-1} below):
\begin{mcor}\label{cor0-0}
$\HBX$ is finite in the four cases in Fact $\ref{fact0-0}$.
\end{mcor}
\noindent
We will also prove variants of Theorem \ref{thm0-0} over local integer rings
 (see Theorems \ref{thm2-1}, \ref{thm4-1} and \ref{thm6-1} below).
As for the finiteness of $\HnrXt$ over local integer rings,
 Spiess proved that $\HnrXt=0$,
 assuming that $\OO_k$ is an $\ell$-adic local integer ring with $\ell \not =p$ and
 that either $\H^2(\Xgf,\O_{\Xgf})=0$ or $\cX$ is a product of two smooth elliptic curves over $S$
   (\cite{Sp}, \S4).
In \cite{SaSa}, the authors extended his vanishing result to a more general situation
   that $\OO_k$ is $\ell$-adic local with $\ell \not =p$ and that $\cX$ has generalized semistable reduction.
Finally we have to remark that there exists
   a smooth projective surface $X$ with $p_g(X) \ne 0$ over a local field $k$ for
   which the condition {\bf H1*} does not hold and such that $\CH^2(X)_{\tor}$ is infinite \cite{AS}.
\par
\bigskip
We next explain an application of the above finiteness result to a cycle class map of arithmetic schemes.
Let us recall the following fact due to Colliot-Th\'el\`ene, Sansuc, Soul\'e and Gros:
\begin{mfact}[{{\bf \cite{CTSS}, \cite{Gr}}}]
\label{thm0}
Let $X$ be a proper smooth variety over a finite field of characteristic 
$\ell > 0$.
Let $p$ be a prime number, which may be the same as $\ell$.
Then the cycle class map restricted to the $p$-primary torsion part
\[ \CH^2(X)_{\ptor} \lra \H^4_{\et}(X,\Z/p^r\Z(2)) \]
is injective for a sufficiently large $r > 0$.
If $\ell \ne p$, then $\Z/p^r\Z(2)$ denotes $\mu_{p^r}^{\otimes 2}$.
If $\ell = p$, then $\Z/p^r\Z(2)$ denotes $\logwitt X r 2 [-2]$
    with $\logwitt X r 2$ the \'etale subsheaf of the logarithmic part of
      the Hodge-Witt sheaf $\witt X r 2$ $($\cite{Bl0}, \cite{Il}$)$.
\end{mfact}
\noindent
In this paper, we study an arithmetic variant of this fact.
We expect that a similar result holds for proper regular arithmetic schemes, i.e.,
regular schemes which are proper flat of finite type over the integer ring of a number field or a local field.
To be more precise, let $k$, $\OO_k$, $\Xgf$ and $\cX$ be as in Theorem \ref{thm0-0}.
The $p$-adic \'etale Tate twist $\T_r(2)=\T_r(2)_\cX$ mentioned before
 replaces $\Z/p^r\Z(2)$ in Fact \ref{thm0}, and there is a cycle class map
\[ \varrho^2_r: \CH^2(\cX)/p^r \lra \H^{4}_{\et}(\cX\hspace{-2pt},\T_r(2)). \]
We are concerned with the induced map
\[  \varrho^2_{\ptor,r} : \CH^2(\cX)_{\ptor} \lra
    \H^{4}_{\et}(\cX\hspace{-2pt},\T_r(2)). \]
It is shown in \cite{Sat2} that the group on the right hand side is finite.
So the injectivity of this map is closely related with
  the finiteness of $\CH^2(\cX)_{\ptor}$.
The second main result of this paper concerns the injectivity of this map:
\begin{mthm}[{{\bf \S\ref{sect4}}}]\label{cor0-1}
Assume that $\H^2(\Xgf,\O_{\Xgf})=0$.
Then $\CH^2(\cX)_{\ptor}$ is finite and $\varrho^2_{\ptor,r}$ is injective for a sufficiently large $r > 0$.
\end{mthm}
\noindent
The finiteness of $\CH^2(\cX)_{\ptor}$ in this theorem is
 originally due to Salberger \cite{Sal}, Colliot-Th\'el\`ene and Raskind
 \cite{CTR1}, \cite{CTR2}. 
Note that there exists a projective smooth surface $V$ over a number field with $\H^2(V,\O_V)=0$
 whose torsion cycle class map
\[  \CH^2(V)_{\ptor} \lra \H^{4}_{\et}(V,\mu_{p^r}^{\otimes 2}) \]
 is not injective for some bad prime $p$ and any $r \ge 1$ \cite{Sur} (cf.\ \cite{PS}).
Our result suggests that
  we are able to recover the injectivity of torsion cycle class maps
  by considering a proper regular model of $V$ over the ring of integers in $k$.
The fundamental ideas of Theorem \ref{cor0-1} are the following.
A crucial point of the proof of Fact \ref{thm0} in \cite{CTSS} and \cite{Gr} is Deligne's proof of the Weil conjecture \cite{De}.
In the arithmetic situation, the role of the Weil conjecture is replaced by the condition {\bf H1}, 
which implies the finiteness of $\CH^2(\Xgf)_{\ptor}$ and $\HBX$ by Theorem \ref{thm0-0}\,(1).
The injectivity result in Theorem \ref{cor0-1} is derived from the finiteness of those objects.
\par
\bigskip
This paper is organized as follows.
In \S\ref{sect1}, we will review some fundamental facts on Galois cohomology groups and Selmer groups which will be used frequently in this paper.
In \S\ref{sect2}, we will prove the finiteness of $\CH^2(\Xgf)_{\ptor}$ in Theorem \ref{thm0-0}\,(1).
In \S\ref{sect3}, we will review $p$-adic \'etale Tate twists briefly and then provide some fundamental lemmas on cycle class maps and unramified cohomology groups.
In \S\ref{sect4}, we will first reduce Theorem \ref{cor0-1} to Theorem \ref{thm0-0}\,(1), and then reduce the finiteness of $\HBX$ in Theorem \ref{thm0-0}\,(1) to Key Lemma \ref{keylem3-1}.
In \S\ref{sect5}, we will prove that key lemma, which will complete the proof of Theorem \ref{thm0-0}\,(1).
\S\ref{sect6} will be devoted to the proof of Theorem \ref{thm0-0}\,(2).
In the appendix A, we will include an observation that
the finiteness of $\H^3_{\ur}(K,\qzp(2))$ is deduced from
the Beilinson--Lichtenbaum conjectures on motivic complexes.
\par
\bigskip
\noindent
{\it Acknowledgements.}
\hspace{1pt}
The research for this article was partially supported by
     JSPS Postdoctoral Fellowship for Research Abroad and EPSRC grant.
The second author expresses his gratitude to
 University of Southern California and The University of Nottingham for their great hospitality.
The authors also express their gratitude to Professors Wayne Raskind, Thomas Geisser and Ivan Fesenko
    for valuable comments and discussions.
\medskip
\newpage
\section*{Notation}
\par
\medskip
\noindent
\ref{intro}.6.
For an abelian group $M$ and a positive integer $n$,
    let ${}_nM$ and $M/n$ be the kernel and the cokernel of the map
       $M \os{\times n}{\lra} M$, respectively.
See \S\ref{sect1.3} below for other notation for abelian groups.
For a field $k$, let $\ol k$ be a fixed separable closure, and
   let $G_k$ be the absolute Galois group $\Gal(\ol k/k)$.
For a discrete $G_k$-module $M$, let $\H^*(k,M)$
   be the Galois cohomology groups $\H^*_{\Gal}(G_k,M)$,
     which is the same as the \'etale cohomology groups of $\Spec(k)$
       with coefficients in the \'etale sheaf associated with $M$.
\par
\bigskip
\noindent
\ref{intro}.7.
Unless indicated otherwise, all cohomology groups of schemes are
   taken over the \'etale topology.
For a scheme $X$, an \'etale sheaf $\FF$ on $X$
    (or more generally an object
     in the derived category of sheaves on $X_{\et}$) and a point $x \in X$,
        we often write
      $\H^*_x(X,\FF)$ for $\H^*_x(\Spec(\O_{X,x}),\FF)$.
For a pure-dimensional scheme $X$ and a non-negative integer $q$,
   let $X^q$ be the set of all points on $X$ of codimension $q$.
For a point $x \in X$, let $\kappa(x)$ be its residue field.
For an integer $n \geq 0$ and a noetherian excellent scheme $X$,
  $\CH_n(X)$ denotes the Chow group of algebraic cycles on $X$
     of dimension $n$ modulo rational equivalence.
If $X$ is pure-dimensional and regular,
     we will often write $\CH^{\dim(X)-n}(X)$ for this group.
For an integral scheme $X$ of finite type over
     $\Spec(\Q)$, $\Spec(\Z)$ or $\Spec(\Z_{\ell})$,
  we define $\CH^2(X,1)$ as the cohomology group, at the middle,
      of the Gersten complex of Milnor $K$-groups
\[ K^M_2(L) \lra \bigoplus_{y\in X^1} \ \kappa(y)^{\times}
            \lra \bigoplus_{x\in X^2} \ \bZ, \]
where $L$ denotes the function field of $X$.
As is well-known, this group coincides with a higher Chow group (\cite{Bl2}, \cite{Le2})
 by the localization theory (\cite{Bl3}, \cite{Le1}) and the Nesterenko-Suslin theorem \cite{NS} (cf. \cite{To}).
\par
\bigskip
\noindent
\ref{intro}.8.
In \S\S\ref{sect3}--\ref{sect6}, we will work under the following setting.
Let $k$ be an algebraic number field or its completion at a finite place.
Let $\OO_k$ be the integer ring of $k$ and put $S:=\Spec(\OO_k)$.
Let $p$ be a prime number, and let $\cX$ be a regular scheme
 which is proper flat of finite type over $S$ and satisfies the following condition:
\begin{quote}
{\it If $p$ is not invertible on $\cX$, then
 $\cX$ has good or semistable reduction at each closed point of $S$ of characteristic $p$.}
\end{quote}
Let $K$ be the function field of $\cX$.
We define $\HnrXt$ and $\HBX$ in the same way
  as in the introduction:
\begin{align*}
 & \HnrXt :=\ker\left(\H^3(\Xff,\qzp(2)) \lra \bigoplus{}_{y\in \cX^1} \
   \H^4_y(\cX\hspace{-2pt},\T_{\infty}(2))\right), \\
 & \HBX \\ & := \Image\Big( \H^3(\Xgf,\qzp(2)) \ra
      \H^3(\Xff,\qzp(2))\Big) \cap \HnrXt,
\end{align*}
where $\T_{\infty}(n)$ denotes $\varinjlim{}_{r \geq 1} \, \T_r(n)$.
If $k$ is an algebraic number field, then this setting is the same as that
  in the introduction.
\par
\bigskip
\noindent
\ref{intro}.\GN.
Let $k$ be an algebraic number field,
   and let $\cX \ra S=\Spec(\OO_k)$ be as in \ref{intro}.8.
In this situation,
    we will often use the following notation.
For a closed point $v\in S$,
   let $\OO_v$ (resp.\ $\kv$) be the completion of $\OO_k$ (resp.\ $k$) at
       $v$, and let $\Fv$ be the residue field of $\kv$.
We put
\[  \Xv:= \cX\otimes_{\OO_k} \OO_v, \qquad
\Xkv := \cX\otimes_{\OO_k} \kv, \qquad \Yv := \cX \otimes_{\OO_k} \Fv \]
 and write $\jv: \Xkv \hra \Xv$ (resp.\ $\iv: \Yv\hra \Xv$)
   for the natural open (resp.\ closed) immersion.
We put $\Ynrv:=\Yv\times_{\Fv}{\ol \Fv}$, and write $\varSigma$ for the set of all closed point on $S$
    of characteristic $p$.
\par
\bigskip
\noindent
\ref{intro}.\LN.
Let $k$ be an $\ell$-adic local field with $\ell$ a prime number,
   and let $\cX \ra S=\Spec(\OO_k)$ be as in \ref{intro}.8.
In this situation,
    we will often use the following notation.
Let $\rf$ be the residue field of $k$ and put
\[ \Xgf := \cX\otimes_{\OO_k} k, \qquad Y := \cX \otimes_{\OO_k} \rf. \]
We write $j:\Xgf \hra \cX$ (resp.\ $i:Y \hra \cX$) for the natural open
   (resp.\ closed) immersion.
Let $\knr$ be the maximal unramified extension of $k$, and
  let $\OO^{\ur}$ be its integer ring.
We put
\[ \Xnr:=\cX \otimes_{\OO_k} \OO^{\ur}, \qquad
  \Xknr:=\cX\otimes_{\OO_k} \knr, \qquad \Ynr:=Y\times_{\rf} \Fb. \]
\newpage
\section{Preliminaries on Galois cohomology}\label{sect1}
\medskip
In this section, we provide some preliminary lemmas
     which will be frequently used in this paper.
Let $k$ be an algebraic number field (global field) or its completion at a finite place (local field).
Let $\OO_k$ be the integer ring of $k$, and put $S:=\Spec(\OO_k)$.
Let $p$ be a prime number. If $k$ is global,
    we often write $\varSigma$ for the set of the closed points on $S$
     of characteristic $p$.
\subsection{Selmer group}\label{sect1.1}
Let $\Xgf$ be a proper smooth variety over $\Spec(k)$, and put $\Xggf : = \Xgf \otimes_k \ol k$.
If $k$ is global,
     we fix a non-empty open subset $U_0 \subset S \sm \varSigma$
         for which there exists a proper smooth morphism
             $\cX_{U_0} \ra U_0$ with $\cX_{U_0} \times_{U_0} k \simeq \Xgf$.
For $v \in S^1$, let $\kv$ and $\Fv$ be as in the notation \ref{intro}.\GN.
In this section we are concerned with $G_k$-modules
\[ 
 V:=\H^i(\Xggf,\qp(n)) \quad \hbox{ and } \quad A:=\H^i(\Xggf,\qzp(n)). \]
For $M=V$ or $A$ and a non-empty open subset $U\subset U_0$,
 let $\H^*(U,M)$ denote the \'etale cohomology groups with
     coefficients in the smooth sheaf on $U_{\et}$ associated to $M$.
\begin{defn}\label{def1-1}
\begin{itemize}
\item[(1)]
Assume that $k$ is local.
Let $\H^1_f(k,V)$ and $\H^1_g(k,V)$ be as defined in \cite{BK2}, $(3.7)$.
For $* \in \{ f,g \}$, we define
\[ \H^1_*(k,A) := \Image \big( \H^1_*(k,V) \lra \H^1(k,A) \big). \]
\item[(2)]
Assume that $k$ is global.
For $M \in \{V,A \}$ and a non-empty open subset $U \subset S$, we define
   the subgroup $\H^1_{f,U}(k,M) \subset \H^1_{\cont}(k,M)$
      as the kernel of the natural map
\[ \H^1_{\cont}(k,M) \lra \prod_{v \in U^1}\ \H^1_{\cont}(k_v,M)/\H^1_f(k_v,M)
    \times
    \prod_{v \in S \sm U}\ \H^1_{\cont}(k_v,M)/\H^1_g(k_v,M). \]
If $U\subset U_0$, we have
\[ \H^1_{f,U}(k,M) =
   \Ker\Big( \H^1(U,M) \lra \prod{}_{v\in S \sm U} \ \H^1_{\cont}(k_v,M)/\H^1_g(k_v,M) \Big). \]
We define the group $\H^1_g(k, M)$ and $\Hind k M$ as
\[ \H^1_g(k, M) := \us{U \subset U_0}\varinjlim~ \H^1_{f,U}(k,M),
   \qquad
   \Hind k M := \us{U \subset U_0}\varinjlim~ \H^1(U,M), \]
where $U$ runs through all non-empty open subsets of $U_0$.
These groups are independent of the choice of $U_0$ and $\cX_{U_0}$
      $($cf.\ \cite{ega4}, $8.8.2.5)$.
\item[(3)]
If $k$ is local,
 we define $\Hind k M$ to be $\H^1_{\cont}(k, M)$ for $M \in \{V,A \}$.
\end{itemize}
\end{defn}
\noindent
Note that $\H^1_{\ind}(k,A)=\H^1(k,A)$.
\subsection{$\bs{p}$-adic point of motives}\label{sect1.2}
We provide a key lemma from $p$-adic Hodge theory
      which play crucial roles in this paper (see Theorem \ref{thm1-1} below).
Assume that $k$ is a $p$-adic local field,
   and that there exists a regular scheme $\cX$
 which is proper flat of finite type over $S=\Spec(\OO_k)$ with $\cX \otimes_{\OO_k} k \simeq \Xgf$
 and which has semistable reduction.
Let $i$ and $n$ be non-negative integers.
Put
$$
V^i:=\H^{i+1}(\Xggf,\qp), \qquad
V^i(n):=V^i\otimes_{\qp}\qp(n),
$$
and
\[ \H^{i+1}(\cX\hspace{-2pt},\tau_{\leq n}Rj_*\qp(n)) :=
   \left\{
   \varprojlim{}_{r \geq 1}~ \H^{i+1}(\cX\hspace{-2pt},\tau_{\leq n} Rj_*\mupnu n)
   \right\}
    \otimes_{\zp} \qp, \]
where $j$ denotes the natural open immersion $\Xgf \hra \cX$.
There is a natural pull-back map
\[ \alpha : \H^{i+1}(\cX\hspace{-2pt},\tau_{\leq n}Rj_*\qp(n)) \lra
      \H^{i+1}(\Xgf,\qp(n)). \]
Let $\H^{i+1}(\cX\hspace{-2pt},\tau_{\leq r}Rj_*\qp(n))^0$ be the kernel of
      the composite map
\[ \alpha' : \H^{i+1}(\cX\hspace{-2pt},\tau_{\leq n}Rj_*\qp(n))
       \os{\alpha}\lra \H^{i+1}(\Xgf,\qp(n))
       \lra \big( V^{i+1} (n) \big)^{G_k}. \]
For this group, there is a composite map
\[ \ol \alpha{} : \H^{i+1}(\cX\hspace{-2pt},\tau_{\leq n}Rj_*\qp(n))^0
      \lra \Faa\H^{i+1}(\Xgf,\qp(n))
      \lra \H^1_{\cont}(k,V^i(n)), \]
whose first arrow is induced by $\alpha$.
The second arrow is an edge homomorphism a Hochschild-Serre spectral sequence
\[ E_2^{u,v}:=\H^u_{\cont}(k,V^v(n)) \Lra \H^{u+v}_{\cont}(\Xgf,\qp(n))(\simeq\H^{u+v}(\Xgf,\qp(n))). \]
and $\F^{\bullet}$ denotes the filtration on $\H^{i+1}(\Xgf,\qp(n))$
 resulting from this spectral sequence.
Concerning the image of $\ol \alpha$, we show the following:
\begin{thm}\label{thm1-1}
Assume that $p \geq n+2$.
Then $\Image(\ol \alpha)=\H^1_g(k,V^i(n))$.
\end{thm}
\begin{pf}
We use the following comparison theorem of log syntomic complexes and $p$-adic vanishing cycles due to
 Kato, Kurihara and Tsuji (\cite{Ka1}, \cite{Ku}, \cite{Ka2}, \cite{Ts2}).
Let $Y$ be the closed fiber of $\cX \ra S$ and let $\iota : Y \hra \cX$ be the natural closed immersion.
\begin{thm}[{{\bf Kato\,/\,Kurihara\,/\,Tsuji}}]\label{thm1-2}
For integers $n,r$ with $0 \leq n \leq p-2$ and $r \geq 1$,
  there is a canonical isomorphism
\[ \eta : s_r^{\log}(n) \isom \iota_*\iota^*(\tau_{\leq n} Rj_*\mupnu n) \quad \hbox { in } \; D^b(\cX_{\et},\pnuz), \]
where $s_r^{\log}(n)=s_r^{\log}(n)_\cX$ is the log syntomic complex defined by Kato \cite{Ka1} $($cf.\ \cite{Ts1}$)$.
\end{thm}
Put
$$
   \H^{*}(\cX\hspace{-2pt}, s_{\qp}^{\log}(n)):=
      \left\{
      \varprojlim{}_{r \geq 1}  \ \H^{*}(\cX\hspace{-2pt},s_r^{\log}(n))
      \right\} \otimes_{\zp} \qp,
$$
and define $\H^{i+1}(\cX\hspace{-2pt},s_{\qp}^{\log}(n))^0$
   as the kernel of the composite map
\[  \H^{i+1}(\cX\hspace{-2pt},s_{\qp}^{\log}(n)) \us{\eta}{\isom}
      \H^{i+1}(\cX\hspace{-2pt},\tau_{\leq n}Rj_*\qp(n))
        \os{\alpha'}\lra \big(V^{i+1}(n) \big)^{G_k}, \]
where we have used the properness of $\cX$ over $S$.
There is an induced map
\[ \ol \eta : \H^{i+1}(\cX\hspace{-2pt},s_{\qp}^{\log}(n))^0 \us{\eta}{\isom}
  \H^{i+1}(\cX\hspace{-2pt},\tau_{\leq n}Rj_*\qp(n))^0 \os{\ol \alpha}\lra \H^1_{\cont}(k,V^i(n)). \]
On the other hand, we have the following fact (\cite{La2}, \cite{Ne2}, Theorem 3.1):
\begin{thm}[{{\bf Langer\,/\,Nekov\'a\v{r}}}]\label{thm1-3}
   $\Image(\ol\eta)$ agrees with $\H^1_g(k,V^i(n))$.
\end{thm}
\noindent
By these facts, we obtain Theorem \ref{thm1-1}.
 \end{pf}
\noindent
\begin{rem}
\begin{enumerate}
\item[(1)]
Theorem $\ref{thm1-3}$
    is an extension of
      the $p$-adic point conjecture raised by
        Schneider in the good reduction case \cite{Sch}.
This conjecture
was proved by Langer-Saito \cite{LS} in a special case and
      by Nekov\'a\v{r} \cite{Ne1} in the general case.
\item[(2)]
Theorem $\ref{thm1-3}$ holds unconditionally on $p$, if we define the space
 $\H^{i+1}(\cX\hspace{-2pt},s_{\qp}^{\log}(n))$ using Tsuji's version of log syntomic complexes $\cS_r\,\ctil(n)\,(r\geq 1)$
 in \cite{Ts1},\,\S$2$.
\end{enumerate}
\end{rem}
\subsection{Elementary facts on $\bs{\zp}$-modules}\label{sect1.3}
For an abelian group $M$, let $M_{\Div}$ be
     its maximal divisible subgroup.
For a torsion abelian group $M$,
     let $\Cotor(M)$ be the cotorsion part $M/M_{\Div}$.
We say that a $\zp$-module $M$ is {\it cofinitely generated
      over $\zp$} (or simply, {\it cofinitely generated}),
         if its Pontryagin dual $\Hom_{\zp}(M,\qzp)$
            is a finitely generated $\zp$-module.
\begin{lem}\label{lem1-3}
Let $0 \ra L \ra M \ra N \ra 0$
    be a short exact sequence of $\zp$-modules.
\begin{enumerate}
\item[(1)]
Assume that $L$, $M$ and $N$ are cofinitely generated.
Then there is a positive integer $r_0$
    such that for any $r \geq r_0$
    we have an exact sequence of finite abelian $p$-groups
\[ 0 \ra {}_{p^r} L \ra {}_{p^r} M \ra {}_{p^r} N
     \ra \Cotor(L) \ra \Cotor(M) \ra \Cotor(N) \ra 0. \]
Consequently, taking the projective limit
     of this exact sequence with respect to $r \geq  r_0$
   there is an exact sequence of
      finitely generated $\zp$-modules
\[ 0 \ra T_p(L) \ra T_p(M) \ra T_p(N)
     \ra \Cotor(L) \ra \Cotor(M) \ra \Cotor(N) \ra 0, \]
where for an abelian group $A$, $T_p(A)$ denotes its $p$-adic Tate module.
\item[(2)]
Assume that $L$ is cofinitely generated
    up to a group of
     finite exponent, i.e.,
       $L_{\Div}$ is cofinitely generated
         and $\Cotor(L)$ has a finite exponent.
Assume further that $M$ is divisible, and that
         $N$ is cofinitely generated and divisible.
Then $L$ and $M$ are cofinitely generated.
\item[(3)]
Assume that $L$ is cofinitely generated
    up to a group of
     finite exponent.
Then for a divisible subgroup $D \subset N$
        and its inverse image $D' \subset M$,
       the induced map $(D')_{\Div} \ra D$ is surjective.
In particular, the natural map $M_{\Div} \ra N_{\Div}$ is surjective.
\item[(4)]
If $L_{\Div}=N_{\Div}=0$, then we have $M_{\Div}=0$.
\end{enumerate}
\end{lem}
\begin{pf}
(1)
There is a commutative diagram with exact rows
\[ \xymatrix{
0 \ar[r] & L \ar[r] \ar[d]_{\times p^r} & M \ar[r] \ar[d]_{\times p^r}
 & N \ar[r] \ar[d]_{\times p^r} & 0 \\
0 \ar[r] & L \ar[r] & M \ar[r] & N \ar[r] & 0.
} \]
One obtains the assertion by applying the snake lemma to this diagram,
 noting $\Cotor(A) \simeq A/p^r$ for a cofinitely generated $\zp$-module $A$
          and a sufficiently large $r \geq 1$.
\par
(2)
Our task is to show that $\Cotor(L)$ is finite.
By a similar argument as for (1),
   there is an exact sequence
      for a sufficiently large $r \geq 1$
\[ 0 \lra {}_{p^r} L \lra {}_{p^r} M \lra {}_{p^r} N \lra \Cotor(L) \lra 0, \]
     where we have used the assumptions on $L$ and $M$.
Hence the finiteness of $\Cotor(L)$ follows from
     the assumption that $N$ is cofinitely generated.
\par
(3)
We have only to show the case $D=N_{\Div}$.
For a $\zp$-module $A$, we have
\[  A_{\Div} = \Image \left(\Hom_{\zp}(\qp,A) \ra A \right) \]
by \cite{Ja1}, Lemma (4.3.a).
Since $\Ext^1_{\zp}(\qp,L)=0$ by the assumption on $L$,
    the following natural map is surjective:
\[  \Hom_{\zp}(\qp,M) \lra \Hom_{\zp}(\qp,N). \]
By these facts, the natural map $M_{\Div} \ra N_{\Div}$ is
      surjective.
\par
(4)
For a $\zp$-module $A$, we have
\[  A_{\Div} =0 \Longleftrightarrow \Hom_{\zp}(\qp,A)=0 \]
by \cite{Ja1}, Remark (4.7).
The assertion follows from this fact and the exact sequence
\[  0 \lra \Hom_{\zp}(\qp,L) \lra \Hom_{\zp}(\qp,M) \lra \Hom_{\zp}(\qp,N) \]
This completes the proof of the lemma.
\end{pf}
\subsection{Divisible part of $\bs{\H^1(k,A)}$}
Let the notation be as in \S\ref{sect1.1}.
We prove here the following general lemma,
    which will be used frequently in \S\S\ref{sect2}--\ref{sect6}:
\begin{lem}\label{lem1-1}
Under the notation in Definition $\ref{def1-1}$ we have
\begin{align*}
\Image\big(\Hind k V \ra \H^1(k,A)\big) & = \H^1(k,A)_{\Div},\\
\Image\big(\H^1_g(k,V) \ra \H^1(k,A)\big) & = \H^1_g(k,A)_{\Div}.
\end{align*}
\end{lem}
\begin{pf}
The assertion is clear if $k$ is local.
Assume that $k$ is global.
Without loss of generality we may assume that $A$ is divisible.
We prove only the second equality and omit the first one
     (see Remark \ref{rem1-1}\,(2) below).
Let $U_0 \subset S$ be as in \S\ref{sect1.1}.
We have
\stepcounter{equation}
\begin{equation}\label{equal:div1}
   \Image\big(\Hg U V \to \H^1(U,A)\big)=  {\Hg U A}{}_{\Div}
\end{equation}
for non-empty open $U \subset U_0$.
This follows from a commutative diagram with exact rows
\[ \xymatrix{
0  \ar[r] & \H^1_{f,U}(k,V) \ar[r] \ar[d] & \H^1(U,V) \ar[r] \ar[d]_{\alpha} &
            \prod{}_{v\in S \sm U} \ \H^1_{\cont}(k_v,V)/\H^1_g(k_v,V) \ar[d]_{\beta} \\
0  \ar[r] &  \H^1_{f,U}(k,A) \ar[r] & \H^1(U,A) \ar[r] &
            \prod{}_{v\in S \sm U} \ \H^1(k_v,A)/\H^1_g(k_v,A)\\
} \]
and the facts that $\Coker(\alpha)$ is finite and that
$\Ker(\beta)$ is finitely generated over $\zp$.
By \eqref{equal:div1},
     the second equality of the lemma is reduced to the following assertion:
\begin{equation}\label{equal:div2}
   \varinjlim_{U\subset U_0}~ \big({\Hg U A}{}_{\Div} \big)=
    \left(\varinjlim_{U\subset U_0}~ \Hg U A)\right){}_{\Div}.
\end{equation}
To show this equality, we will prove the following sublemma:
\addtocounter{thm}{2}
\begin{sublem}\label{sublem2-1}
For an open subset $U \subset U_0$, put
\[ C_U:=\Coker\big(\Hg {U_0} A \to \Hg U A \big). \]
Then there exists a non-empty open subset $U_1 \subset U_0$ such that
   the quotient
    $C_U/C_{U_1}$ is divisible
     for any open subset $U \subset U_1$.
\end{sublem}
We first finish our proof of \eqref{equal:div2} admitting this sublemma.
Let $U_1 \subset U_0$ be a non-empty open subset as in Sublemma \ref{sublem2-1}.
Noting that $\Hg U A$ is cofinitely generated,
   there is an exact sequence of finite groups
\[ \Cotor \big(\Hg {U_1} A \big)
    \lra \Cotor\big(\Hg U A \big) \lra \Cotor(C_U/C_{U_1}) \lra 0 \]
for open $U \subset U_1$ by Lemma \ref{lem1-3}\,(1).
By this exact sequence and Sublemma \ref{sublem2-1},  the natural map
     $\Cotor(\Hg {U_1} A) \ra \Cotor(\Hg U A)$
       is surjective for any open $U \subset U_1$,
     which implies that the inductive limit
\[ \us{U\subset U_0}\varinjlim \Cotor({\Hg U A}) \]
is a finite group.
The equality \eqref{equal:div2} follows easily from this.
\par
\medskip
\noindent
{\it Proof of Sublemma \ref{sublem2-1}.}
We need the following general fact:
\begin{sublem}\label{sublem2-2}
Let $\ul N=\{N_\lambda\}_{\lambda\in \varLambda}$ be an inductive system of
cofinitely generated $\zp$-modules indexed by a filtered set $\varLambda$
such that $\Coker(N_{\lambda}\to N_{\lambda'})$ is divisible for
any two $\lambda,\lambda'\in \varLambda$ with $\lambda' \geq \lambda$.
Let $L$ be a cofinitely generated $\zp$-module and
$\{f_\lambda:N_\lambda\to L\}_{\lambda\in \varLambda}$
   be $\zp$-homomorphisms compatible with
   the transition maps of $\ul N$.
Then there exists $\lambda_0\in \varLambda$ such
that $\Coker\big(\Ker(f_{\lambda_0})\to \Ker(f_{\lambda})\big)$
is divisible for any $\lambda \geq \lambda_0$.
\end{sublem}
\begin{pf*}{\it Proof of Sublemma \ref{sublem2-2}}
Let $f_{\infty}:N_\infty\to L$ be the limit of $f_\lambda$.
The assumption on $\ul N$ implies that
for any two $\lambda, \lambda' \in \varLambda$ with $\lambda' \geq \lambda$,
   the quotient $\Image(f_{\lambda'})/\Image(f_\lambda)$ is divisible, so that
\addtocounter{equation}{2}
\begin{equation}\label{surj:div}
\Cotor(\Image(f_\lambda)) \to \Cotor(\Image(f_{\lambda'}))
    \hbox{ is surjective}.
\end{equation}
By the equality $\Image(f_\infty)
   = \varinjlim{}_{\lambda \in \varLambda}\ \Image(f_{\lambda})$,
  there is a short exact sequence
\[   0  \lra  \us{\lambda \in \varLambda}{\varinjlim}\
               \big(\Image(f_\lambda)_{\Div}\big)
      \lra  \Image(f_\infty)
      \lra  \us{\lambda \in \varLambda}{\varinjlim}\ \Cotor(\Image(f_\lambda))
      \lra  0, \]
and the last term is finite by the fact \eqref{surj:div}
   and the assumption that
$L$ is cofinitely generated.
Hence we get
\[ \us{\lambda \in \varLambda}{\varinjlim} \ \big(\Image(f_\lambda)_{\Div}\big)= \Image(f_\infty)_{\Div}. \]
Since $\Image(f_\infty)_{\Div}$ has finite corank,
there exists an element $\lambda_0\in \varLambda$ such that
$\Image(f_\lambda)_{\Div} = \Image(f_\infty)_{\Div}$
for any $\lambda \geq \lambda_0$.
This fact and \eqref{surj:div} imply
   the equality
\begin{equation}\label{surj2:div}
     \Image(f_\lambda) = \Image(f_{\lambda_0})\;\;
    \hbox{ for any } \lambda \geq \lambda_0.
\end{equation}
Now let $\lambda \in \varLambda$ satisfy $\lambda \geq \lambda_0$.
Applying the snake lemma to the commutative diagram
\[ \xymatrix{
  & N_{\lambda_0} \ar[r] \ar[d]_{f_{\lambda_0}}
 & N_\lambda \ar[r] \ar[d]_{f_{\lambda}}
 & N_{\lambda}/N_{\lambda_0} \ar[r] \ar[d] & 0\\
0 \ar[r] & L \ar@{=}[r] & L \ar[r] & \, 0,
} \]
we get an exact sequence
\[ \Ker(f_{\lambda_0})
    \lra \Ker(f_{\lambda_0}) \lra  N_{\lambda}/N_{\lambda_0} \os{0}{\lra}
       \Coker(f_{\lambda_0}) \isom \Coker(f_{\lambda}), \]
which proves Sublemma \ref{sublem2-2}, beucase $N_{\lambda}/N_{\lambda_0}$ is divisible by
assumption.
\end{pf*}
We now turn to the proof of Sublemma \ref{sublem2-1}.
For non-empty open $U\subset U_0$,
   there is a commutative diagram with exact rows
\[ {\small \begin{CD}
@. @. \H^1(U_0,A) @. ~ \lra ~ @. \H^1(U,A) @. \; \lra \; @.
    \underset{v\in U_0 \sm U}{\bigoplus} A(-1)^{G_{\rf_v}}
       @. \; \os{\beta_U}{\lra} \;\, @.  \H^2(U_0,A) \\
@. @. @V{r_{U_0}}VV @. @V{r_U}VV @. @V{\alpha_U}VV @. @.\\
0 @. \; \lra \; @.
\underset{v\in S \sm U_0}{\bigoplus} ~\H^1_{/g}(k_v,A) @. \;\, \lra \; @.
\underset{v\in S \sm U}{\bigoplus}~\H^1_{/g}(k_v,A) @. \;\, \lra \; @.
\underset{v\in U_0 \sm U}{\bigoplus}~\H^1_{/g}(k_v,A),
\end{CD} } \]
where we put
\[ \H^1_{/g}(k_v,A):= \H^1(k_v,A)/\H^1_{g}(k_v,A) \]
for simplicity.
The upper row is obtained from a localization exact sequence
      of \'etale cohomology groups
     and the isomorphism
\[ \H^2_v(U_0,A)\simeq \H^1(k_v,A)/\H^1(\rf_v,A)
    \simeq  A(-1)^{G_{\rf_v}}
     \quad\hbox{ for } v \in U_0 \sm U, \]
where we have used the fact that
   the action of $G_k$ on $A$ is unramified at $v\in U_0$.
The map $\alpha_U$ is obtained from the facts that
   $\H^1_g(k_v,A)=\H^1(k_v,A)_{\Div}$
      if $v \not\in \varSigma$ and
     that $\H^1(\rf_v,A)$ is divisible
       (recall that $A$ is assumed to be divisible).
It gives
\begin{equation}\label{equal:div3}
\Ker(\alpha_U)=\bigoplus_{v\in U_0 \sm U} \
    \left(A(-1)^{G_{\rf_v}}\right){}_{\Div}.
\end{equation}
Now let $\phi_U$ be the composite map
\[ \phi_U : \Ker(\alpha_U) \; \hra \;
 \bigoplus_{v\in U_0 \sm U}~ A(-1)^{G_{\rf_v}}
      \os{\beta_U}\lra \H^2(U_0,A), \]
and let
\[ \psi_U: \Ker(\phi_U) \lra \Coker(r_{U_0}) \]
be the map induced by the above diagram.
Note that
\[ C_U \simeq \Ker(\psi_U), \; \hbox{ since } \; \Hg U A =\Ker(r_U).  \]
By \eqref{equal:div3}, the inductive system
$\{\Ker(\alpha_U)\}_{U\subset U_0}$
    and the maps $\{ \phi_U \}_{U\subset U_0}$ satisfy
   the assumptions in Sublemma \ref{sublem2-2}.
Hence there exists a non-empty open subset
     $U' \subset U_0$ such that
$\Ker(\phi_U)/\Ker(\phi_{U'})$ is divisible for any open $U \subset U'$.
Then applying Sublemma \ref{sublem2-2} again to the inductive system
   $\{\Ker(\phi_U)\}_{U \subset U'}$ and the maps
     $\{ \psi_U \}_{U\subset U'}$, we conclude that
        there exists a non-empty open subset $U_1 \subset U'$ such that
    the quotient
\[ \Ker(\psi_U)/\Ker(\psi_{U_1})=C_U/C_{U_1} \] is divisible for any open subset $U \subset U_1$.
This completes the proof of Sublemma \ref{sublem2-1} and Lemma \ref{lem1-1}.
\end{pf}
\addtocounter{thm}{3}
\begin{rem}\label{rem1-1}
\begin{enumerate}
\item[(1)]
By the argument after Sublemma $\ref{sublem2-1}$, $\Cotor(\H^1_g(k,A))$ is finite.
\item[(2)]
One obtains the first equality in Lemma $\ref{lem1-1}$
     by replacing the local terms $\H^1_{/g}(\kv,A)$
         in the above diagram with $\Cotor(\H^1(\kv,A))$.
\end{enumerate}
\end{rem}
\subsection{Cotorsion part of $\bs{\H^1(k,A)}$}\label{sect1.5}
Assume that $k$ is global, and
   let the notation be as in \S\ref{sect1.1}.
We investigate here the boundary map
\[ \delta_{U_0} : \H^1(k,A)  \lra
        \bigoplus_{v\in (U_0)^1} \ A(-1)^{G_{\Fv}} \]
arising from the localization theory in \'etale topology and the purity for discrete valuation rings.
Concerning this map, we prove the following standard lemma,
    which will be used in our proof of
      Theorem \ref{thm0-0}:
\begin{lem}\label{lem1-2}
\begin{enumerate}
\item[(1)]
The map
\[ \delta_{U_0,\Div} : \H^1(k,A)_{\Div}  \lra
        \bigoplus_{v\in (U_0)^1} \ \big(A(-1)^{G_{\Fv}} \big){}_{\Div} \]
induced by $\delta_{U_0}$ has cofinitely generated cokernel.
\item[(2)]
The map
\[ \delta_{U_0,\Cotor} : \Cotor(\H^1(k,A))  \lra
        \bigoplus_{v\in (U_0)^1} \ \Cotor\big(A(-1)^{G_{\Fv}} \big) \]
induced by $\delta_{U_0}$ has finite kernel and
    cofinitely generated cokernel.
\end{enumerate}
\end{lem}
\noindent
We have nothing to say about the finiteness of the cokernel of these maps.
\begin{pf}
For a non-empty open $U \subset U_0$,
    there is a commutative diagram of
     cofinitely generated $\zp$-modules
\[ \xymatrix{
    & \H^1(U,A)_{\Div} \ar[r]^{\gamma_U\qquad\;\;\;\;} \ar@{_{(}->}[d]
    & \bigoplus{}_{v\in U_0 \sm U}\ \big((A(-1)^{G_{\Fv}}\big){}_{\Div} \ar@{_{(}->}[d] \\
\H^1(U_0,A) \ar[r]
  & \H^1(U,A)  \ar[r]^{\alpha_U\qquad\;\;} & \bigoplus{}_{v\in U_0 \sm U}\ A(-1)^{G_{\Fv}} \ar[r]^{\quad\;\;\beta_U}
  & \H^2(U_0,A), } \]
where the lower row is obtained from
     the localization theory in \'etale cohomology and the purity
       for discrete valuation rings,
       and $\gamma_U$ is induced by $\alpha_U$.
Let
\[ f_U : \Cotor(\H^1(U,A)) \lra
      \bigoplus_{v\in U_0 \sm U}~ \Cotor\big(A(-1)^{G_{\Fv}} \big) \]
be the map induced by $\alpha_U$.
By a diagram chase, we obtain an exact sequence
\[ \ker(f_U) \lra \Coker(\gamma_U) \lra \Coker(\alpha_U) \lra \Coker(f_U) \lra 0. \]
Taking the inductive limit with respect to all non-empty open
     subsets $U \subset U_0$,
       we obtain an exact sequence
\[ \ker(\delta_{U_0,\Cotor}) \lra \Coker(\delta_{U_0,\Div}) \lra
         \us{U \subset U_0}{\varinjlim}~\Coker(\alpha_U)
           \lra \Coker(\delta_{U_0,\Cotor}) \lra 0, \]
where we have used Lemma \ref{lem1-1} to obtain the equalities
    $\ker(\delta_{U_0,\Cotor}) = \varinjlim{}_{U \subset U_0}~\ker(f_U)$
     and $\Coker(\delta_{U_0,\Div})= \varinjlim{}_{U \subset U_0}~
      \Coker(\gamma_U)$.
Since $\varinjlim{}_{U \subset U_0}~\Coker(\alpha_U)$
    is a subgroup of $\H^2(U_0,A)$,
       it is cofinitely generated.
Hence the assertions in Lemma \ref{lem1-2} are
      reduced to showing that $\ker(\delta_{U_0,\Cotor})$ is finite.
We prove this finiteness assertion.
The lower row of the above diagram yields exact sequences
\stepcounter{equation}
\begin{align}
&\Cotor(\H^1(U_0,A)) \lra \Cotor(\H^1(U,A))  \lra
       \Cotor(\Image(\alpha_U)) \lra 0,
\label{UC.1.1} \\
&T_p(\Image(\beta_U)) \lra
       \Cotor(\Image(\alpha_U)) \lra
       \bigoplus_{v\in U_0 \sm U} \ \Cotor\big(A(-1)^{G_{\Fv}} \big),
\label{UC.1.2}
\end{align}
where the second exact sequence arises from the short exact sequence
\[ 0 \lra \Image(\alpha_U) \lra
       \bigoplus_{v\in U_0 \sm U}~A(-1)^{G_{\Fv}}
          \lra  \Image(\beta_U) \lra 0 \]
(cf.\ Lemma \ref{lem1-3}\,(1)).
Taking the inductive limit of \eqref{UC.1.1} with respect to
    all non-empty open $U \subset U_0$,
      we obtain the finiteness of the kernel of the map
\begin{equation}\notag
   \Cotor(\H^1(k,A))  \lra
       \varinjlim_{U \subset U_0} \ \Cotor(\Image(\alpha_U)).
\end{equation}
Taking the inductive limit of \eqref{UC.1.2} with respect to
    all non-empty open $U \subset U_0$,
    we see that the kernel of the map
\begin{equation}\notag
     \us{U \subset U_0}{\varinjlim}~ \Cotor(\Image(\alpha_U)) \lra
     \bigoplus_{v\in (U_0)^1} \ \Cotor\big(A(-1)^{G_{\Fv}} \big),
\end{equation}
is finite, because we have
\[ \us{U \subset U_0}{\varinjlim}~ T_p(\Image(\beta_U)) \subset T_p(\H^2(U_0,A)) \]
and the group on the right hand side is
   a finitely generated $\zp$-module.
Thus $\ker(\delta_{U_0,\Cotor})$ is finite and we obtain
     Lemma \ref{lem1-2}.
\end{pf}
\subsection{Local-global principle}\label{sect1.6}
Let the notation be as in \S\ref{sect1.1}.
If $k$ is local, then the Galois cohomological dimension
      $\cd(k)$ is $2$ (cf.\ \cite{Se}, II.4.3).
In the case that $k$ is global,
   we have $\cd(k)=2$ either if $p \geq 3$ or if $k$ is totally imaginary.
Otherwise, $\H^q(k,A)$ is finite $2$-torsion
        for $q \geq 3$ (cf.\ loc.\ cit., II.4.4, Proposition 13, II.6.3, Theorem B).
As for the second Galois cohomology groups,
 the following local-global principle due to Jannsen (\cite{Ja2}, \S4, Theorem 4)
 plays a fundamental role in this paper (see also loc.\ cit., \S7, Corollary 7):
\begin{thm}[{{\bf Jannsen}}]\label{thm1-4}
Assume that $k$ is global and that $i \not= 2(n-1)$.
Let $P$ be the set of all places of $k$.
Then the map
\[ \H^2(k,\H^i(\Xggf,\qzp(n))) \lra \bigoplus_{v \in P} \ \H^2(k_v,\H^i(\Xggf,\qzp(n))) \]
has finite kernel and cokernel,
   and the map
\[ \H^2(k,\H^i(\Xggf,\qzp(n))_{\Div}) \lra \bigoplus_{v \in P}\ \H^2(k_v,\H^i(\Xggf,\qzp(n))_{\Div}) \]
is bijective.
\end{thm}
\noindent
We apply these facts to
    the filtration $\F^{\bullet}$ on $\H^*(\Xgf,\qzp(n))$
       resulting from the Hochschild-Serre spectral sequence
\stepcounter{equation}
\begin{equation}\label{TOR.2.3}
   E_2^{u,v}=\H^u(k,\H^v(\Xggf,\qzp(n)))
     \Longrightarrow \H^{u+v}(\Xgf,\qzp(n)).
\end{equation}
\stepcounter{thm}
\begin{cor}\label{cor1-5}
Assume that $k$ is global and that $i \not= 2n$.
Then$:$
\begin{enumerate}
\item[(1)]
$\Fbb\H^i(\Xgf,\qzp(n))_{\Div}$ is cofinitely generated
     and $\Cotor(\Fbb\H^i(\Xgf,\qzp(n)))$ has a finite exponent.
\item[(2)]
For $v \in P$, put $\Xkv:=\Xgf \otimes_k \kv$. Then the natural maps
\begin{align*}
\Fbb\H^i(\Xgf,\qzp(n)) & \lra
    \bigoplus_{v \in P} \ \Fbb\H^i(\Xkv,\qzp(n)),\\
\Fbb\H^i(\Xgf,\qzp(n))_{\Div} & \lra
    \bigoplus_{v \in P} \ \Fbb\H^i(\Xkv,\qzp(n))_{\Div}
\end{align*}
have finite kernel and cokernel
$($and the second map is surjective$)$.
\end{enumerate}
\end{cor}
\begin{pf}
Let $\OO_k$ be the integer ring of $k$, and put $S:=\Spec(\OO_k)$.
Note that the set of all finite places of $k$ agrees with $S^1$.
\par
(1)
The group $\H^2(\kv,\H^{i-2}(\Xggf,\qzp(n))_{\Div})$
            is divisible and cofinitely generated for any $v \in S^1$,
   and it is zero if
       $p \hspace{-2pt} \not \hspace{-1pt}\vert \hspace{2pt} v$
             and $\Xgf$ has good reduction at $v$,
               by the local Poitou-Tate duality \cite{Se}, II.5.2,
                 Th\'eor\`eme 2 and
       Deligne's proof of the Weil conjecture \cite{De}
         (see \cite{sato2}, Lemma 2.4 for details).
The assertion follows from this fact and Theorem \ref{thm1-4}.
\par
(2) We prove the assertion only for the first map.
The assertion for the second map is similar and left to the reader.
For simplicity, we assume that
\begin{quote}
  $(\sharp)$ $p \geq 3$ or $k$ is totally imaginary.
\end{quote}
Otherwise one can check the assertion by repeating the same
       arguments as below in the category of abelian groups
          modulo finite abelian groups.
By $(\sharp)$, we have $\cd_p(k)=2$ and there is a commutative diagram
\[ \xymatrix{
\H^2(k,\H^{i-2}(\Xggf,\qzp(n))) \ar[r] \ar@{->>}[d] &
    \bigoplus{}_{v \in S^1}\ \H^2(k_v,\H^{i-2}(\Xggf,\qzp(n))) \ar@{->>}[d] \\
\Fbb\H^i(\Xgf,\qzp(n)) \ar[r] &
    \bigoplus{}_{v \in S^1}\ \Fbb\H^i(\Xkv,\qzp(n)), }\]
where the vertical arrows are edge homomorphisms
     of Hochschild-Serre spectral sequences and these arrows are surjective.
Since
\[ \H^2(k_v,\H^{i-2}(\Xggf,\qzp(n)))=0 \quad \hbox{ for archimedean places $v$} \]
by $(\sharp)$, the top horizontal arrow has finite kernel and cokernel by Theorem \ref{thm1-4}.
Hence it is enough to show that the right vertical arrow has finite kernel.
For any $v \in S^1$, the $v$-component of this map
      has finite kernel by Deligne's criterion \cite{De0}
         (see also \cite{sato2}, Remark 1.2).
If $v$ is prime to $p$ and $\Xgf$ has good reduction at $v$,
    then the $v$-component is injective.
Indeed, there is an exact sequence resulting from a Hochschild-Serre spectral sequence
 and the fact that $\cd(\kv)=2$:
\begin{align*} \H^{i-1}(\Xkv,\qzp(n)) & \os{d}{\lra} \H^{i-1}(\Xggf,\qzp(n)))^{G_{\kv}} \\
 & \lra \H^2(k_v,\H^{i-2}(\Xggf,\qzp(n))) \lra \Fbb\H^i(\Xkv,\qzp(n)). \end{align*}
The edge homomorphism $d$ is surjective by the commutative square
\[\xymatrix{
   \H^{i-1}(\Yv,\qzp(n)) \ar@{->>}[r] \ar[d]
 & \H^{i-1}(\Ynrv,\qzp(n)))^{G_{\Fv}} \ar[d]^{\hspace{-2pt}\wr} \\
   \H^{i-1}(\Xkv,\qzp(n)) \ar[r]^{d \quad} & \H^{i-1}(\Xggf,\qzp(n)))^{G_{\kv}}. }\]
Here $\Yv$ denotes the reduction of $\Xgf$ at $v$ and $\Ynrv$ denotes $\Yv \otimes_{\Fv} \ol {\Fv}$.
The left (resp.\ right) vertical arrow  arises from the proper base change theorem
         (resp.\ proper smooth base change theorem),
  and the top horizontal arrow is surjective by the fact that $\cd(\Fv)=1$.
Thus we obtain the assertion.
\end{pf}

\newpage
\section{Finiteness of torsion in a Chow group}\label{sect2}
\medskip
\subsection{Finiteness of $\bs{\CH^2(\Xgf)_{\ptor}}$}\label{sect2.1}
Let $k,S,p$ and $\varSigma$ be as in the beginning of \S\ref{sect1},
    and let $\Xgf$ be a proper smooth geometrically integral
      variety over $\Spec(k)$.
We introduce the following technical condition:
\par
\vspace{8pt}
\begin{enumerate}
\item[]
{\bf H0:}
{\it The group $\H^3_{\et}(\Xggf,\qp(2))^{G_k}$ is trivial}.
\end{enumerate}
\par
\vspace{8pt}
\noindent
If $k$ is global,
{\bf H0} always holds by Deligne's proof of
    the Weil conjecture \cite{De}.
When $k$ is local,
 {\bf H0} holds if $\dim(\Xgf)=2$ or if $X$ has good reduction (cf.\ \cite{CTR2}, \S6);
  it is in general a consequence of the monodromy-weight conjecture.
The purpose of this section is to show the following result,
     which is a generalization of a result of Langer \cite{La3},
       Proposition 3 and implies the finiteness assertion on 
$\CH^2(\Xgf)_{\ptor}$
  in Theorem \ref{thm0-0}\,(1):
\begin{thm}\label{thm2-1}
Assume {\bf H0}, {\bf H1*} and
      either $p \geq 5$ or the equality
\begin{equation}\tag{$*_g$}\label{intro.1'}
\H^1_g(k,\H^2(\Xggf,\qzp(2)))_{\Div} =
    \H^1(k,\H^2(\Xggf,\qzp(2)))_{\Div}.
\end{equation}
Then $\CH^2(\Xgf)_{\ptor}$ is finite.
\end{thm}
\begin{rem}\label{rem2-1}
\begin{enumerate}
\item[(1)]
\eqref{intro.1'} holds if $\H^2(\Xgf,\O_{\Xgf})=0$ or if $k$ is $\ell$-adic local with $\ell \not= p$.
\item[(2)]
Crucial facts to this theorem are Lemmas $\ref{lem2-0}$, $\ref{lem2-1}$ and $\ref{lem2-1a}$ below.
The short exact sequence in Lemma $\ref{lem2-0}$ is an important consequence of
 the Merkur'ev-Suslin theorem \cite{MS}.
\end{enumerate}
\end{rem}
\subsection{Regulator map}\label{sect2.2}
We recall here the definition of the regulator maps
\begin{equation}\label{TOR.1.1}
\begin{CD}
   \regXkLam: \CH^2(\Xgf,1)\otimes\Lam
     @>>> \Hind k {\H^2(\Xggf,\Lam(2))}
\end{CD}
\end{equation}
with $\Lam = \qp$ or $\qzp$, assuming {\bf H0}.
The general framework on \'etale Chern class maps and regulator maps
is due to Soul\'e \cite{So1}, \cite{So2}.
We include here a more elementary construction of $\regXkLam$,
   which will be useful in this paper.
Let $K:=k(\Xgf)$ be the function field of $\Xgf$.
Take an open subset $U_0 \subset S \sm \varSigma= S[p^{-1}]$
      and a smooth proper scheme $\cX_{U_0}$ over $U_0$
         satisfying
         $\cX_{U_0} \times _{U_0} \Spec(k) \simeq \Xgf$.
For an open subset $U \subset U_0$, put $\XU:=\cX_{U_0} \times_{U_0} U$
    and define
$$
\Naa\H^3(\XU,\mupnu 2) :=
    \Ker\big(\H^3(\XU,\mupnu 2) \to \H^3(K,\mupnu 2)\big).
$$
\stepcounter{thm}
\begin{lem}\label{lem2-0}
For an open subset $U \subset U_0$,
    there is an exact sequence
\[ 0 \lra \CH^2(\XU,1)/p^r \lra \Naa\H^3(\XU,\mupnu 2) \lra {}_{p^r} \CH^2(\XU) \lra 0 \]
See $\S\ref{intro}.7$ for the definition of $\CH^2(\XU,1)$.
\end{lem}
\begin{pf}
The following argument is due to Bloch [Bl], Lecture 5.
We recall it for the convenience of the reader.
There is a localization spectral sequence
\stepcounter{equation}
\begin{equation}\label{TOR.1.2}
E_1^{u,v}=\bigoplus_{x \in (\XU)^u} \
              \H^{u+v}_{x}(\XU,\mupnu 2)
                   \Lra \H^{u+v}(\XU,\mupnu 2).
\end{equation}
By the relative smooth purity, there is an isomorphism
\[ E_1^{u,v} \simeq
    \bigoplus_{x \in (\XU)^u} \, \H^{v-u}(x, \mupnu {2-u}),
\]
which implies that $\Naa\H^3(\XU,\mupnu 2)$ is isomorphic to
the cohomology of the Bloch-Ogus complex
\[ \H^2(K,\mupnu 2) \lra \bigoplus_{y\in (\XU)^1} \ \H^1(y,\mu_{p^r}) \lra \bigoplus_{x\in (\XU)^2} \ \pnuz. \]
By the Merkur'ev-Suslin theorem \cite{MS}, this complex is
     isomorphic to the Gersten complex
\[ K^M_2(K)/p^r \lra \bigoplus_{y\in (\XU)^1} \ k(y)^{\times}/p^r
          \lra \bigoplus_{x\in (\XU)^2} \ \pnuz. \]
On the other hand, there is an exact sequence obtained by a diagram chase
\[ 0 \lra \CH^2(\XU,1)\otimes\pnuz \lra \CH^2(\XU,1;\pnuz) \lra
      {}_{p^r}\CH^2(\XU)  \lra 0. \]
Here $\CH^2(\XU,1;\pnuz)$ denotes the cohomology of the above Gersten complex
  and it is isomorphic to $\Naa\H^3(\XU,\mupnu 2)$.
Thus we obtain the lemma.
\end{pf}
\medskip
Put 
\[ M^q:=\H^q(\Xggf,\Lam(2)) \quad \hbox{ with } \quad  \Lam \in \{ \qp,\qzp \}. \]
For an open subset $U \subset U_0$
  let $\H^*(U,M^q)$ be the \'etale cohomology with coefficients in the smooth sheaf associated with $M^q$.
There is a Leray spectral sequence
\[ E_2^{u,v}= \H^u(U,M^v) \Lra \H^{u+v}(\XU,\Lam(2)). \]
By Lemma \ref{lem2-0}, there is a natural map
\[ \CH^2(\XU,1)\otimes\Lam \lra \H^3(\XU,\Lam(2)). \]
Noting that $E_2^{0,3}$ is zero or finite by {\bf H0}, we define the map
\[ \regXULam: \CH^2(\XU,1) \otimes \Lam \lra \H^1(U,M^2) \]
as the composite of the above map with an edge homomorphism
  of the Leray spectral sequence.
Finally we define $\regXkLam$ in \eqref{TOR.1.1} by passing to the limit over
     all non-empty open $U\subset U_0$.
Our construction of $\regXkLam$ does not depend on the choice of $U_0$ or $\cX_{U_0}$.
\stepcounter{thm}
\begin{rem}\label{lem2-01}
By Lemma $\ref{lem1-1}$, {\bf H1} always implies {\bf H1*}.
If $k$ is local, {\bf H1*} conversely implies {\bf H1}.
As for the case that $k$ is global,
     one can check that {\bf H1*} implies {\bf H1},
        assuming that the group
   $\Ker(\CH^2(\cX_{U_0}) \to \CH^2(\Xgf))$
     is finitely generated up to torsion
   and that the Tate conjecture for divisors holds for almost all closed fibers
      of $\cX_{U_0}/U_0$.
\end{rem}
\subsection{Proof of Theorem \ref{thm2-1}}
We start the proof of Theorem \ref{thm2-1},
     which will be completed in \S\ref{sect2.5} below.
By Lemma \ref{lem2-0},
    there is an exact sequence
\begin{equation}\label{TOR.2.1}
 0 \lra \CH^2(\Xgf,1)\otimes\qzp \os{\phi}\lra
   \Naa\H^3(\Xgf,\qzp(2)) \lra \CH^2(\Xgf)_{\ptor} \lra 0,
\end{equation}
where we put
\[ \Naa\H^3(\Xgf,\qzp(2)):= \Ker(\H^3(\Xgf,\qzp(2))\to \H^3(K,\qzp(2))). \]
In view of \eqref{TOR.2.1},
     Theorem \ref{thm2-1} is reduced to the following two propositions:
\stepcounter{thm}
\begin{prop}\label{prop2-2}
\begin{enumerate}
\item[(1)]
If $k$ is local,
    then $\CH^2(\Xgf)_{\ptor}$ is cofinitely generated over $\zp$.
\item[(2)]
Assume that $k$ is global, and that $\Coker\big(\regXkqzp\big){}_{\Div}$
    is cofinitely generated over $\Z_p$.
Then $\CH^2(\Xgf)_{\ptor}$ is cofinitely generated over $\zp$.
\end{enumerate}
\end{prop}
\begin{prop}\label{prop2-1}
Assume {\bf H0}, {\bf H1*} and either $p\geq 5$ or \eqref{intro.1'}.
Then we have
\[ \Image(\phi) = \Naa\H^3(\Xgf,\qzp(2))_{\Div}. \]
\end{prop}
\noindent
We will prove Proposition \ref{prop2-2} in \S\ref{sect2.4} below,
     and Proposition \ref{prop2-1} in \S\ref{sect2.5} below.
\begin{rem}\label{rem2-3}
\begin{enumerate}
\item[(1)]
If $k$ is local, then
    $\H^3(\Xgf,\qzp(2))$ is cofinitely generated.
Hence Proposition $\ref{prop2-2}\,(1)$ immediately follows
     from the exact sequence \eqref{TOR.2.1}.
\item[(2)]
When $k$ is global,
   then $\H^1(k,A)_{\Div}/\H^1_g(k,A)_{\Div}$ with $A:=\H^2(\Xggf,\qzp(2))$ is
      cofinitely generated by Lemma $\ref{lem1-1}$.
Hence {\bf H1*} implies the second assumption
  of Proposition $\ref{prop2-2}\,(2)$.
\end{enumerate}
\end{rem}
\noindent
Let $\F^{\bullet}$ be
     the filtration on $\H^*(\Xgf,\qzp(2))$
      resulting from the Hochschild-Serre spectral sequence
         \eqref{TOR.2.3}.
The following fact due to Salberger (\cite{Sal}, Main Lemma 3.9)
      will play key roles in our proof of the above two propositions:
\begin{lem}[{{\bf Salberger}}]\label{lem2-1}
The following group has a finite exponent$:$
\[ \Naa\H^3(\Xgf,\qzp(2)) \cap \Fbb \H^3(\Xgf,\qzp(2)). \]
\end{lem}
\subsection{Proof of Proposition \ref{prop2-2}}\label{sect2.4}
For (1), see Remark \ref{rem2-3}\,(1).
We prove (2).
Put
\[ \H^3:=\H^3(\Xgf,\qzp(2)) \quad\hbox{ and }\quad
    \vG :=\phi(\CH^2(\Xgf,1)\otimes\qzp) \subset \H^3 \]
(cf.\ \eqref{TOR.2.1}).
Let $\F^{\bullet}$ be the filtration on $\H^3$
      resulting from the spectral sequence sequence \eqref{TOR.2.3},
        and put $\Naa\H^3:= \Naa\H^3(\Xgf,\qzp(2))$.
We have $\vG \subset (\Faa\H^3)_{\Div}=(\H^3)_{\Div}$ by {\bf H0}, and
there is a filtration on $\H^3$
\[ 0 \subset \vG+(\Fbb\H^3)_{\Div} \subset (\Faa\H^3)_{\Div} \subset \H^3. \]
By \eqref{TOR.2.1}, the inclusion $\Naa\H^3 \subset \H^3$
       induces an inclusion $\CH^2(\Xgf)_{\ptor} \subset \H^3/\vG$.
We show that the image of this inclusion is cofinitely generated,
      using the above filtration on $\H^3$.
It suffices to show the following lemma:
\begin{lem}\label{claim2-2}
\begin{itemize}
\item[(1)]
The kernel of $\CH^2(\Xgf)_{\ptor} \to \H^3/(\vG+(\Fbb\H^3)_{\Div})$
    is finite.
\item[(2)]
The image of $\CH^2(\Xgf)_{\ptor} \to \H^3/(\Faa\H^3)_{\Div}$ is finite.
\item[(3)]
Put $M:=(\Faa\H^3)_{\Div}/(\vG+(\Fbb\H^3)_{\Div})$.
Then the assumption of Proposition $\ref{prop2-2}\,(2)$ implies that
$M$ is cofinitely generated.
\end{itemize}
\end{lem}
\begin{pf}
(1) There is an exact sequence
\[ 0 \to (\Naa\H^3  \cap (\Fbb\H^3)_{\Div})/(\vG \cap (\Fbb\H^3)_{\Div})
     \to \CH^2(\Xgf)_{\ptor} \to \H^3/(\vG+(\Fbb\H^3)_{\Div}). \]
Hence (1) follows from Lemma \ref{lem2-1}
     and Corollary \ref{cor1-5}\,(1).
\par
(2)
Let $U_0$ and $\cX_{U_0} \ra U_0$ be as in \S\ref{sect2.2}.
For non-empty open $U \subset U_0$,
    there is a commutative diagram up to a sign
\[ \xymatrix{
\Naa\H^3(\cX_U,\qzp(2)) \ar[r] \ar@{_{(}->}[d] & \CH^2(\cX_U)\otimes \zp \ar[d]^{\varrho} \\
\H^3(\cX_U,\qzp(2)) \ar[r] & \H^4(\cX_U,\zp(2)) } \]
by the same argument as for [CTSS], \S1, Proposition 1.
Here the top arrow is the composite of
    $\Naa\H^3(\cX_U,\qzp(2)) \ra \CH^2(\cX_U)_{\ptor}$
      (cf.\ Lemma \ref{lem2-0}) with the natural inclusion.
The bottom arrow is a Bockstein map and the right vertical arrow is the cycle class map of $\cX_U$.
Taking the inductive limit with respect to all non-empty
       $U \subset U_0$,
we obtain a commutative diagram (up to a sign)
\[\xymatrix{ \Naa\H^3 \ar[r] \ar@{_{(}->}[d] & \CH^2(\Xgf)\otimes \zp \ar[d]^{\varrho_{\ind}} \\
\H^3 \ar[r] & \H^4_{\ind}(\Xgf,\zp(2)), }\]
where $\H^4_{\ind}(\Xgf,\zp(2))$ is defined as the inductive
     limit of $\H^4(\cX_U,\zp(2))$ with respect to $U \subset U_0$.
Now this diagram yields a commutative diagram (up to a sign)
\[ \xymatrix{
\CH^2(\Xgf)_{\ptor} \ar[r] \ar[d] & \CH^2(\Xgf)\otimes \zp \ar[d]^{\varrho_{\cont}} \\
\H^3/(\Faa\H^3)_{\Div} \; \ar@{^{(}->}[r] & \H^4_{\cont}(\Xgf,\zp(2)), }\]
where $\H^*_{\cont}(\Xgf,\zp(2))$ denotes the continuous \'etale cohomology \cite{Ja1}
  and the bottom arrow is injective by {\bf H0} and loc.\ cit., Theorem (5.14).
The image of $\varrho_{\cont}$ is finitely generated over $\zp$ by \cite{Sa}, Theorem (4-4).
This proves (2).
\par
(3) We put
\[ N:= (\Faa\H^3)_{\Div}/\{ \vG+ (\Fbb\H^3 \cap (\Faa\H^3)_{\Div} ) \}
   = \Coker\big(\regXkqzp\big){}_{\Div}, \]
which is cofinitely generated by assumption.
There is an exact sequence
\[ (\Fbb\H^3 \cap (\Faa\H^3)_{\Div})/(\Fbb\H^3)_{\Div} \lra M \lra N \lra 0, \]
where the first group has a finite exponent
     by Corollary \ref{cor1-5}\,(1),
        $N$ is divisible and cofinitely generated, and $M$ is divisible.
Hence $M$ is cofinitely generated by Lemma \ref{lem1-3}\,(2).
This completes the proof of Lemma \ref{claim2-2}
    and Proposition \ref{prop2-2}.
\end{pf}
\subsection{Proof of Proposition \ref{prop2-1}}\label{sect2.5}
We put
\[ \NFaa\H^3(\Xgf,\qzp(2))
    := \Naa\H^3(\Xgf,\qzp(2)) \cap \Faa\H^3(\Xgf,\qzp(2)). \]
Note that $\Naa\H^3(\Xgf,\qzp(2))_{\Div}
       = \NFaa\H^3(\Xgf,\qzp(2))_{\Div}$ by {\bf H0}.
There is an edge homomorphism of the spectral sequence \eqref{TOR.2.3}
\begin{equation}\label{TOR.4.1}
   \psi : \Faa\H^3(\Xgf,\qzp(2)) \lra \H^1(k,\H^2(\Xggf,\qzp(2))).
\end{equation}
The composite of $\phi$
    in \eqref{TOR.2.1} and $\psi$
       agrees with $\regXkqzp$.
Thus by Lemma \ref{lem2-1}, we are reduced to the following lemma,
  which generalizes \cite{LS}, Lemma (5.7) and extends \cite{La0}, Lemma (3.3):
\stepcounter{thm}
\begin{lem}\label{lem2-1a}
Assume either $p\geq 5$ or \eqref{intro.1'}.
Then $\psi(\NFaa\H^3(\Xgf,\qzp(2))_{\Div})$ is contained in $\H^1_g(k,\H^2(\Xggf,\qzp(2)))$.
\end{lem}
\noindent
We start the proof of this lemma.
The assertion is obvious
        under the assumption \eqref{intro.1'}.
Hence we are done if $k$ is $\ell$-adic local with $\ell \not = p$ (cf.\ Remark \ref{rem2-1}\,(1)).
It remains to deal with the following two cases:
\begin{quote}
(1) $k$ is $p$-adic local with $p \geq 5$.
\par\noindent
(2) $k$ is global and $p \geq 5$.
\end{quote}
Put $A := \H^2(\Xggf,\qzp(2))$ for simplicity.
We first reduce the case (2) to the case (1).
Suppose that $k$ is global.
Then there is a commutative diagram
\[ \xymatrix{
\NFaa\H^3(\Xgf,\qzp(2))_{\Div} \ar[r] \ar[d] & \H^1(k,A) \ar[d] \\
\prod{}_{v \in S^1} \ \NFaa\H^3(\Xkv,\qzp(2))_{\Div}
    \ar[r] & \prod{}_{v \in S^1} \ \H^1(k_v,A), }\]
where the vertical arrows are natural restriction maps.
By this diagram and the definition of $\H^1_{g}(k,A)$,
      the case (2) is reduced to the case (1).
We prove the case (1).
We first reduce the problem to the case where $\Xgf$ has semistable reduction.
By the alteration theorem of de Jong \cite{dJ},
     there exists a proper flat generically finite
        morphism $X' \ra \Xgf$ such that
          $X'$ is projective smooth over $k$
            and has a proper flat regular model
             over the integral closure $\OO'$ of $\OO_k$ in $\vG(X',\O_{X'})$
                with semistable reduction.
There is a commutative diagram
\[\xymatrix{ \NFaa\H^3(\Xgf,\qzp(2))_{\Div} \ar[r] \ar[d] & \H^1(k,A)  \ar[r] \ar[d] & \H^1(k,A)/\H^1_{g}(k,A) \ar[d] \\
\NFaa\H^3(X',\qzp(2))_{\Div} \ar[r] & \H^1(L,A') \ar[r] & \H^1(L,A')/\H^1_{g}(L,A'), }\]
where we put $L:=\Frac(\OO')$ and
   $A':=\H^2(X'\otimes_L \ol k,\qzp(2))$,
     and the vertical arrows are natural restriction maps.
Our task is to show that the composite of the upper row is zero.
Because $X'$ and $\Xgf$ are proper smooth varieties over $k$,
    the restriction map $r: A \ra A'$ has a quasi-section
       $s: A' \ra A$ with $s \circ r = d \cdot \id_A$,
         where $d$ denotes the extension degree of
           the function field of $X'\otimes_L \ol k$ over that of $\Xggf$.
Hence by the functoriality of $\H^1_{g}(k,A)$ in $A$,
     the right vertical arrow in the above diagram has finite kernel,
      and the problem is reduced to showing that
        the composite of the lower row is zero.
Thus we are reduced to the case where $\Xgf$ has a proper
    flat regular model $\cX$ over $S=\Spec(\OO_k)$ with semistable reduction.
We prove this case.
Let $j:\Xgf \hra \cX$ be the natural open immersion.
There is a natural injective map
\[ \xymatrix{ \alpha_r : \H^3(\cX\hspace{-2pt},\tau_{\leq 2} Rj_*\mupr 2)
     \; \ar@{^{(}->}[r] & \H^3(\Xgf,\mupr 2) } \]
induced by the natural morphism
$\tau_{\leq 2} Rj_*\mupr 2 \ra Rj_*\mupr 2$.
By Theorem \ref{thm1-1}, it suffices to show the following two lemmas
 (see also Remark \ref{rem2-2} below):
\begin{lem}\label{lem2-2}
$\Naa\H^3(\Xgf,\mupr 2) \subset \Image(\alpha_r)$
     for any $r \geq 1$.
\end{lem}
\begin{lem}\label{lem2-3}
Put
\begin{align*}
\H^3(\cX\hspace{-2pt},\tau_{\leq 2} Rj_*\qzp(2)) & :=
    \varinjlim_{r \geq 1}~
       \H^3(\cX\hspace{-2pt},\tau_{\leq 2} Rj_*\mupr 2), \quad \hbox{ and }\\
\H^3(\cX\hspace{-2pt},\tau_{\leq 2} Rj_*\qzp(2))^0 & :=
     \ker\big( \H^3(\cX\hspace{-2pt},\tau_{\leq 2} Rj_*\qzp(2)) \ra
              \H^3(\Xggf,\qzp(2)) \big).
\end{align*}
Then the canonical map
\[ \H^3(\cX\hspace{-2pt},\tau_{\leq 2} Rj_*\qp(2))^0 \lra
    \H^3(\cX\hspace{-2pt},\tau_{\leq 2} Rj_*\qzp(2))^0 \]
has finite cokernel,
    where $\H^3(\cX\hspace{-2pt},\tau_{\leq 2} Rj_*\qp(2))^0$ is as we defined in $\S\ref{sect1.2}$.
\end{lem}
\begin{pf*}{\it Proof of Lemma $\ref{lem2-2}$}
We use the following fact due to Hagihara (\cite{Sat2}, A.2.4, A.2.6),
       whose latter vanishing will be used later in \S\ref{sect5}:
\begin{lem}[{{\bf Hagihara}}]\label{lem2-h}
Let $n,r$ and $c$ be integers with $n \geq 0$ and $r,c \geq 1$.
Then for any $q \leq n + c$
   and any closed subscheme $Z \subset Y$ with  $\codim_\cX(Z) \geq c$,
   we have
$$
   \H^q_Z(\cX\hspace{-2pt},\tau_{\leq n}Rj_*\mupr n)=0
   =\H^{q+1}_Z(\cX\hspace{-2pt},\tau_{\geq n+1}Rj_*\mupr n).
$$
\end{lem}
\medskip
To show Lemma \ref{lem2-2},
we compute the local-global spectral sequence
\[ E_1^{u,v}=\bigoplus_{x\in \cX^a} \
   \H^{u+v}_x(\cX\hspace{-2pt},\tau_{\leq 2} Rj_*\mupr 2) \Lra
   \H^{u+v}(\cX\hspace{-2pt},\tau_{\leq 2} Rj_*\mupr 2). \]
By the first part of Lemma \ref{lem2-h}
   and the smooth purity for points on $\Xgf$,
     we have
\[ E_1^{u,v}=
   \begin{cases}
   \H^v(K,\mupr 2)  \qquad & \text{(if $u=0$)} \\
   \displaystyle \bigoplus{}_{x\in \Xgf^u}~ \H^{v-u}(x,\mupr {2-u}) \qquad &
                                \text{(if $v\leq 2$).}\\
\end{cases} \]
Repeating the same computation as in the proof of Lemma \ref{lem2-0}, we obtain
\[  \Naa\H^3(\Xgf,\mupr 2) \simeq E_2^{1,2}=E_{\infty}^{1,2} \hra \H^3(\cX\hspace{-2pt},\tau_{\leq 2} Rj_*\mupr 2), \]
which implies Lemma \ref{lem2-2}.
\end{pf*}
\begin{rem}\label{rem2-2}
Lemma $\ref{lem2-2}$ extends a result of Langer-Saito $($\cite{LS}, Lemma $(5.4))$ to regular semistable families and
removes the assumption in \cite{La0}, Lemma $(3.1)$
 concerning Gersten's conjecture for algebraic $K$-groups.
Therefore the same assumption in loc.\ cit., Theorem {\rm A} has been removed as well.
\end{rem}
\begin{pf*}{\it Proof of Lemma \ref{lem2-3}}
By the Bloch-Kato-Hyodo theorem
    on the structure of $p$-adic vanishing cycles (\cite{bk}, \cite{Hy}),
     there is a distinguished triangle of the following form in $D^b(\cX_{\et})$
         (cf.\ \cite{Sat2}, (4.3.3)):
\[  \tau_{\leq 2} Rj_*\mu_{p^r}^{\otimes 2}
     \lra \tau_{\leq 2} Rj_*\mu_{p^{r+s}}^{\otimes 2}
      \lra \tau_{\leq 2} Rj_*\mu_{p^{s}}^{\otimes 2}
       \lra (\tau_{\leq 2} Rj_*\mu_{p^r}^{\otimes 2})[1] \]
Taking \'etale cohomology groups, we obtain a long exact sequence
\addtocounter{equation}{5}
\begin{equation}\label{TOR.4.2}
\begin{CD}
   \dotsb @.~ \lra~ @. \H^q(\cX\hspace{-2pt},\tau_{\leq 2} Rj_*\mu_{p^r}^{\otimes 2})
      @.~ \lra ~@. \H^q(\cX\hspace{-2pt},\tau_{\leq 2} Rj_*\mu_{p^{r+s}}^{\otimes 2})
      @.~ \lra ~@. \H^q(\cX\hspace{-2pt},\tau_{\leq 2} Rj_*\mu_{p^{s}}^{\otimes 2})\\
      @. \lra @. \H^{q+1}(\cX\hspace{-2pt},\tau_{\leq 2} Rj_*\mu_{p^r}^{\otimes 2})
      @. \lra @. \hspace{-77pt}\dotsb.
\end{CD}
\end{equation}
We claim that
   $\H^q(\cX\hspace{-2pt},\tau_{\leq 2}Rj_*\mupr 2)$ is finite
     for any $q$ and $r$.
Indeed, the claim is reduced to the case $r=1$
   by the exactness of \eqref{TOR.4.2}
    and this case follows from the Bloch-Kato-Hyodo theorem
       mentioned above and
        the properness of $\cX$ over $S$.
Hence taking the projective limit of \eqref{TOR.4.2} with respect to $r$
     and then taking the inductive limit with respect to $s$
       we obtain a long exact sequence
\[ \begin{CD}
   \dotsb @.~ \lra~ @. \H^q(\cX\hspace{-2pt},\tau_{\leq 2} Rj_*\zp(2))
      @.~ \lra ~@. \H^q(\cX\hspace{-2pt},\tau_{\leq 2} Rj_*\qp(2))
      @.~ \lra ~@. \H^q(\cX\hspace{-2pt},\tau_{\leq 2} Rj_*\qzp(2))\\
      @. \lra @. \H^{q+1}(\cX\hspace{-2pt},\tau_{\leq 2} Rj_*\zp(2))
      @. \lra @. \hspace{-80pt}\dotsb,
\end{CD} \]
where $\H^q(\cX\hspace{-2pt},\tau_{\leq 2} Rj_*\zp(2))$
     is finitely generated over $\zp$ for any $q$.
The assertion in the lemma easily follows from this exact sequence
    and a similar long exact sequence of \'etale cohomology groups
       of $\Xggf$.
The details are straight-forward and left to the reader.
\end{pf*}
This completes the proof of Lemma \ref{lem2-1a},
   Proposition \ref{prop2-1} and Theorem \ref{thm2-1}.

\newpage

\section{Cycle class map and unramified cohomology}\label{sect3}
\medskip
Let $k, S, p, \cX$ and $K$ be as in the notation \ref{intro}.8.
In this section we give a brief review of $p$-adic \'etale Tate twists
   and provide some preliminary results on cycle class maps.
The main result of this section is Corollary \ref{prop0-0} below.
\subsection{$\bs{p}$-adic \'etale Tate twist}\label{sect3.1}
Let $n$ and $r$ be positive integers.
We recall here the fundamental properties (S1)--(S7) listed below
    of the object
$\pnurX = \pnurX_\cX \in D^b(\cX_{\et},\pnuz)$ introduced by the second author \cite{Sat2}.
The properties (S1), (S2), (S3) and (S4)
    characterizes $\pnurX$ uniquely
      up to a unique isomorphism
     in $D^b(\cX_{\et},\pnuz)$.
\medskip
\begin{enumerate}
\item[{\bf (S1)}]
 {\it There is an isomorphism $t : \pnurX \vert_V
   \simeq \mu_{p^r}^{\otimes n}$ on $V:=\cX[p^{-1}]$.}
\smallskip
\item[{\bf (S2)}]
 {\it $\pnurX$ is concentrated in $[0,n]$.}
\smallskip
\item[{\bf (S3)}]
 {\it Let $Z \subset \cX$ be a locally closed regular subscheme of pure codimension $c$ with $\ch(Z)=p$.
Let $i : Z \ra \cX$ be the natural immersion.
Then there is a canonical Gysin isomorphism}
\[ \gys^n_{i}: \logwitt Z r {n-c} [-n-c] \isom
     \tau_{\leq n+c} Ri^!\pnurX \quad \hbox{ in } \; D^b(Z_{\et},\pnuz), \]
where $\logwitt Z r q$ denotes the \'etale subsheaf of the logarithmic part of
       the Hodge-Witt sheaf $\witt Z r q$ $($\cite{Bl0}, \cite{Il}$)$.
\smallskip
\item[{\bf (S4)}]
 {\it For $x\in \cX$ and $q\in \bZ_{\ge 0}$, we define $\pnuz(q) \in D^b(x_\et, \Z/p^r\Z)$ as
\[ \pnuz(q) := \begin{cases}
              \mu_{p^r}^{\otimes q} & (\hbox{if }\, \ch(x)\not=p)\\
              \logwitt x r q [-q] & (\hbox{if }\, \ch(x) =p). \end{cases} \]
Then for $y,x\in \cX$ with $c:=\codim(x)=\codim(y)+1$,
there is a commutative diagram
\[ \xymatrix{
\H^{n-c+1}(y,\pnuz(n-c+1)) \ar[r]^{\quad -\dval} \ar[d]_{\gys^n_{i_y}}
 & \H^{n-c}(x,\pnuz(n-c)) \ar[d]^{\gys^n_{i_x}} \\
\H^{n+c-1}_y(\cX\hspace{-2pt},\pnurX) \ar[r]^{\dloc}
 & \H^{n+c}_x(\cX\hspace{-2pt},\pnurX). } \]
Here for $z \in \cX$, $\gys^n_{i_z}$ is induced by
    the Gysin map in {\rm (S3)}
    $($resp.\ the absolute purity \cite{RZ}, \cite{Th}, \cite{Fu}$)$
    if $\ch(z)=p$
       $($resp.\ $\ch(z)\not=p)$.
The arrow $\dloc$ denotes the boundary map in localization theory
   and $\dval$ denotes the boundary map of Galois cohomology groups
      due to Kato \cite{Ka3}, $\S1$.}
\smallskip
\item[{\bf (S5)}]
 {\it Let $Y$ be the union of the fibers of $\cX/S$
      of characteristic $p$.
We define the \'etale sheaf $\nu_{Y,r}^{n-1}$ on $Y$ as
\[ \nu_{Y,r}^{n-1}:= \Ker\Big(\dval:
   \bigoplus{}_{y\in Y^0}\ i_{y*}\logwitt y r {n-1} \lra
   \bigoplus{}_{x\in Y^1}\ i_{x*}\logwitt x r {n-2}\Big), \]
where for $y \in Y$, $i_y$ denotes the canonical map $y \hra Y$.
Let $i$ and $j$ be as follows$:$
\[ \begin{CD}
      V=\cX[p^{-1}] @>{j}>> \cX @<{i}<< Y.
\end{CD} \]
Then there is a distinguished triangle
   in $D^b(\cX_{\et},\pnuz)$
\[ \begin{CD}
i_*\nu_{Y,r}^{n-1}[-n-1] @>{g}>>
\pnurX @>{t'}>> \tau_{\leq n}Rj_*\mu_{p^r}^{\otimes n} @>{\sigma}>>
    i_*\nu_{Y,r}^{n-1}[-n],
\end{CD} \]
where $t'$ is induced by the isomorphism $t$ in {\rm (S1)}
  and the acyclicity property {\rm (S2)}.
The arrow $g$ arises from the Gysin morphisms in {\rm (S3)},
    $\sigma$ is induced by the boundary maps of Galois cohomology groups
         $($cf.\ {\rm (S4))}.}
\smallskip
\item[{\bf (S6)}]
 {\it There is a canonical distinguished triangle
   of the following form in $D^b(\cX_{\et})$}:
\[ \begin{CD}
   \pnmurX
     @>>> \psnX @>{\delta_{s,r}}>> \pnurX[1]
     @>{\ul{p}{}^s}>> \pnmurX [1].
\end{CD} \]
\smallskip
\item[{\bf (S7)}]
{\it $\H^i(\cX\hspace{-2pt},\pnurX)$ is finite for any $r$ and $i$ $($by the properness of $\cX)$}.
\end{enumerate}
\begin{rem}
These properties of $\pnurX$ deeply rely on the computation on the \'etale sheaf of $p$-adic vanishing cycles
      due to Bloch-Kato \cite{bk} and Hyodo \cite{Hy}.
\end{rem}
\begin{lem}\label{lem4-1}
Put
\begin{align*}
 \H^q(\cX\hspace{-2pt},\zpnX) & := \us{r \ge 1}\varprojlim \ \H^q(\cX\hspace{-2pt},\pnurX), \qquad
 \H^q(\cX\hspace{-2pt},\qzpnX) := \us{r \ge 1}\varinjlim \ \H^q(\cX\hspace{-2pt},\pnurX), \\
 \H^q(\cX\hspace{-2pt},\qpnX) & := \H^q(\cX\hspace{-2pt},\zpnX) \otimes_{\zp} \qp.
\end{align*}
Then there is a long exact sequence of $\zp$-modules
\[ \begin{CD}
 \dotsb @.~ \lra~ @. \H^q(\cX\hspace{-2pt},\zpnX) @.~ \lra ~@. \H^q(\cX\hspace{-2pt},\qpnX)
 @.~ \lra ~@. \H^q(\cX\hspace{-2pt},\qzpnX) \\
 @. \lra @. \H^{q+1}(\cX\hspace{-2pt},\zpnX) @. \lra @. \hspace{-55pt}\dotsb.
\end{CD} \]
$\H^q(\cX\hspace{-2pt},\zpnX)$ is finitely generated over $\zp$,
$\H^q(\cX\hspace{-2pt},\qzpnX)$ is cofinitely generated over $\zp$,
 and $\H^q(\cX\hspace{-2pt},\qpnX)$ is finite-dimensional over $\qp$.
\end{lem}
\begin{pf}
The assertions immediately follow from (S6) and (S7). The details are straight-forward and left to the reader.
\end{pf}
\subsection{Cycle class map}
Let us review the definition of the cycle map
\[ \varrho^n_r: \CH^n(\cX)/p^r \lra \H^{2n}(\cX\hspace{-2pt},\pnurX). \]
Consider the local-global spectral sequence
\begin{equation}\label{INJ.2.1}
   E_1^{u,v}=\bigoplus_{x\in \cX^u} \
     \H^{u+v}_{x}(\cX\hspace{-2pt},\pnurX)
      \Lra \H^{u+v}(\cX\hspace{-2pt},\pnurX).
\end{equation}
By (S3) and the absolute cohomological purity \cite{Fu}
   (cf.\ \cite{RZ}, \cite{Th}), we have
\begin{equation}\label{INJ.2.2}
   E_1^{u,v} \simeq
   \bigoplus_{x\in \cX^u} \, \H^{v-u}(x,\pnuz(n-u))
     \;\;\quad\text{ for } v \leq n.
\end{equation}
This implies that there is an edge homomorphism
$E_2^{n,n}\to \H^{2n}(\cX\hspace{-2pt},\pnurX)$ with
\begin{align*} E_2^{n,n} & \simeq
 \Coker \Big(\dval : \bigoplus{}_{y \in \cX^{n-1}} \ \H^1(y ,\pnuz(1)) \lra
   \bigoplus{}_{x\in \cX^n} \ \H^0(x,\pnuz)\Big) \\
      & = \CH^n(\cX)/p^r, \end{align*}
where $\dval$ is as in (S4). We define $\varrho^n_r$ as the composite map
\[ \varrho^n_r : \CH^n(\cX)/p^r \simeq E_2^{n,n} \lra \H^{2n}(\cX\hspace{-2pt},\pnurX). \]
In what follows, we restrict our attention to the case $n=2$.
\addtocounter{thm}{2}
\begin{lem}\label{lem3-1}
Let $Y$ be the fiber of $\cX \to S$ over a closed point of $S$, and let $K$ be the function field of $\cX$.
Put
\begin{align*}
 \Naa\H^i(\cX\hspace{-2pt},\pnuX) &:=  \Ker\big(
    \H^i(\cX\hspace{-2pt},\pnuX) \to \H^i(K,\mupnu 2) \big), \\
   \Nbb\H^i_Y(\cX\hspace{-2pt},\pnuX) &:=
   \Ker\Big(\H^i_Y(\cX\hspace{-2pt},\pnuX) \to
    \bigoplus{}_{y \in Y^0}~ \H^i_y(\cX\hspace{-2pt},\pnuX) \Big).
\end{align*}
\begin{itemize}
\item[(1)]
$\Naa\H^3(\cX\hspace{-2pt},\pnuX)$
   is isomorphic to the cohomology of the Gersten complex modulo $p^r$
\[ K^M_2(K)/p^r \lra
    \bigoplus_{y\in \cX^1} \ \kappa(y)^{\times}/p^r \lra
    \bigoplus_{x\in \cX^2} \ \pnuz, \]
and there is an exact sequence
\[   0 \lra \CH^2(\cX\hspace{-2pt},1)/p^r \lra \Naa\H^3(\cX\hspace{-2pt},\pnuX)
     \lra {}_{p^r}\CH^2(\cX) \lra 0 \]
See $\S\ref{intro}.4$ for the definition of $\CH^2(\cX\hspace{-2pt},1)$.
\item[(2)]
There are isomorphisms
\begin{align*}
 \H^3_Y(\cX\hspace{-2pt},\pnuX) & \simeq
    \Ker\Big( \dval :
   \bigoplus{}_{y\in \cX^1 \cap Y} \ \k(y)^{\times}/p^r \lra
   \bigoplus{}_{x\in \cX^2 \cap Y} \ \pnuz
    \Big),\\
   \Nbb\H^4_Y(\cX\hspace{-2pt},\pnuX) & \simeq
   \Coker\Big( \dval :
   \bigoplus{}_{y\in \cX^1 \cap Y} \ \k(y)^{\times}/p^r \lra
   \bigoplus{}_{x\in \cX^2 \cap Y} \ \pnuz \Big) \\
    & = \CH_{d-2}(Y)/p^r,
\end{align*}
where $d$ denotes the Krull dimension of $\cX$.
\end{itemize}
\end{lem}
\begin{pf}
(1) follows from the same argument as for Lemma \ref{lem2-0}.
One can prove (2) in the same way as for (1), using the spectral sequence
\[ E_1^{u,v}= \bigoplus_{x \in \cX^u \cap Y} \ \H^{u+v}_x(\cX\hspace{-2pt},\pnurX)
    \Lra \H^{u+v}_Y(\cX\hspace{-2pt},\pnurX) \]
and the purity isomorphism
\[ E_1^{u,v} \simeq
   \bigoplus_{x\in \cX^u \cap Y} \, \H^{v-u}(x,\pnuz(n-u))
     \;\;\quad\text{ for } v \leq n \]
instead of \eqref{INJ.2.1} and \eqref{INJ.2.2}.
The details are straight-forward and left to the reader.
\end{pf}
\begin{cor}\label{cor3-1}
${}_{p^r}\CH^2(\cX)$ is finite for any $r \geq 1$,
     and $\CH^2(\cX)_{\ptor}$ is cofinitely generated.
\end{cor}
\begin{pf}
The finiteness of ${}_{p^r}\CH^2(\cX)$ follows from the exact sequence
      in Lemma \ref{lem3-1}\,(1) and (S7) in \S\ref{sect3.1}.
The second assertion follows from Lemma \ref{lem4-1} and the facts that
     $\CH^2(\cX)_{\ptor}$ is a subquotient of $\H^3(\cX\hspace{-2pt},\qzpX)$.
\end{pf}
The following lemma will be used in the proof of Theorem \ref{cor0-1}.
\begin{lem}\label{lem3-2}
Assume that the natural inclusion
\[\xymatrix{ i_0 : \Naa\H^3(\cX\hspace{-2pt},\qzpX) \; \ar@{^{(}->}[r] & \H^3(\cX\hspace{-2pt},\qzpX) }\]
has finite cokernel.
Then there exists a positive integer $r_0$
    such that the kernel of the map
\[ \varrho^2_{\ptor,r}: \CH^2(\cX)_{\ptor} \lra \H^4(\cX\hspace{-2pt},\pnuX) \]
agrees with $(\CH^2(\cX)_{\ptor})_{\Div}$ for any $r \geq r_0$.
\end{lem}
\begin{pf}
The following argument is essentially the same as the proof of
   \cite{CTSS}, Corollaire 3.
We recall it for the convenience of the reader.
By the exact sequence in Lemma \ref{lem3-1}\,(1), we have
\[ \Cotor(\Naa\H^3(\cX\hspace{-2pt},\qzpX)) \simeq \Cotor(\CH^2(\cX)_{\ptor}). \]
By (S4) in \S\ref{sect3.1} and the same argument as [CTSS], \S1,
   one can show the commutativity of the following diagram up to a sign:
\[\xymatrix{  & \Naa\H^3(\cX\hspace{-2pt},\qzpX) \ar[r] \ar@{_{(}->}[d]_{i_0}
 & \CH^2(\cX)_{\ptor} \ar[d]^{\varrho^2_{\ptor,r}} \\
\H^3(\cX\hspace{-2pt},\qzpX) \ar[r]^{\times p^r}
& \H^3(\cX\hspace{-2pt},\qzpX) \ar[r]^{\delta_{\infty,r}} & \H^4(\cX\hspace{-2pt},\pnuX), }\]
where the lower row is an exact sequence
  induced by the distinguished triangle
\[ \begin{CD}
   \qzpX @>{\times p^r}>> \qzpX  @>{\delta_{\infty,r}}>>  \pnuX [1]
         @>>> \qzpX [1] \end{CD} \]
obtained by taking the inductive limit of the distinguished triangle of (S6) with respect to $s>0$.
The above diagram induces the following commutative diagram up to a sign:
\[ \xymatrix{
 & \Cotor(\Naa\H^3(\cX\hspace{-2pt},\qzpX)) \ar[r]^{\quad\sim} \ar@{_{(}->}[d]_{\ol {i_0}}
 & \Cotor(\CH^2(\cX)_{\ptor}) \ar[d]^{\ol{\varrho^2_{\ptor,r}}} \\
\Cotor(\H^3(\cX\hspace{-2pt},\qzpX)) \ar[r]^{\times p^r}
 & \Cotor(\H^3(\cX\hspace{-2pt},\qzpX)) \ar[r]^{\quad\ol{\delta_{\infty,r}}} & \H^4(\cX\hspace{-2pt},\pnuX),\\
} \]
where the lower row remains exact.
 $\H^4(\cX\hspace{-2pt},\pnuX)$ and $\Cotor(\H^3(\cX\hspace{-2pt},\qzpX))$ are finite by (S7) and Lemma \ref{lem4-1}.
Hence $\ol{\delta_{\infty,r}}$ is injective for any $r$
  with $p^r \cdot \Cotor(\H^3(\cX\hspace{-2pt},\qzpX)) = 0$.
The finiteness of $\Coker(i_0)$ implies
the injectivity of $\ol {i_0}$.
Thus we obtain Lemma \ref{lem3-2}.
\end{pf}
\begin{rem}\label{rem3-1}
If $k$ is $\ell$-adic local with $\ell \not = p$, then we have $\qzpX=\qzp(2)$ by definition and
\[ \H^3(\cX\hspace{-2pt},\qzpX) = \H^3(\cX\hspace{-2pt},\qzp(2)) \simeq \H^3(Y,\qzp(2)) \] 
by the proper base-change theorem, where $Y$ denotes the closed fiber of $\cX \ra S$.
The last group is finite by Deligne's proof of the Weil conjecture \cite{De}.
Hence $\varrho^2_{\ptor,r}$ for $\cX$ is injective for a sufficiently large $r \geq 1$ by Lemma $\ref{lem3-2}$.
On the other hand, if $k$ is global or $p$-adic local,
   then $\H^3(\cX\hspace{-2pt},\qzpX)$ is not in general finite.
Therefore we consider
   the finiteness of the group $\H^3_{\ur}(K,\Xgf;\qzp(2))$
   to investigate the injectivity of $\varrho^2_{\ptor,r}$.
\end{rem}
\begin{cor}\label{prop0-0}
If $\HBX$ is finite,
then there is a positive integer $r_0$ such that
    $\ker(\varrho^2_{\ptor,r})=(\CH^2(\cX)_{\ptor})_{\Div}$ for any $r \geq r_0$.
\end{cor}
\begin{pf}
Let $i_0$ 
  be as in Lemma \ref{lem3-2}.
Since $\Coker(i_0)$ is a subgroup of $\HBX$ (cf.\ \eqref{UC.0.1} below),
the assumption implies that $\Coker(i_0)$ is finite.
Hence the assertion follows from Lemma \ref{lem3-2}.
\end{pf}
\begin{rem}
By the spectral sequence \eqref{INJ.2.1} and the isomorphisms in \eqref{INJ.2.2} with $n=2$,
 there is an exact sequence
\addtocounter{equation}{6}
\begin{equation}\label{UC.0.1}
\begin{CD}
0 @.~\; \lra \;~@. \Naa\H^3(\cX\hspace{-2pt},\qzpX) @>{i_0}>>
     \H^3(\cX\hspace{-2pt},\qzpX) @.~\; \lra \;~@. \H^3_{\ur}(K,\qzp(2)) \\
    @.~ \lra ~@. \CH^2(\cX) \otimes \qzp  @>{\varrho^2_{\qzp}}>>
          \H^4(\cX\hspace{-2pt},\qzpX) @.~ \lra ~@. \dotsb. \qquad \qquad \quad \;\;
\end{CD}
\end{equation}
Because the groups $\H^*(\cX\hspace{-2pt},\qzpX)$ are cofinitely generated $($cf.\ Lemma $\ref{lem4-1})$,
this exact sequence implies that $\H^3_{\ur}(K,\qzp(2))$ is cofinitely generated
      if and only if $\CH^2(\cX) \otimes \qzp$ is cofinitely generated.
\end{rem}
We mention here some remarks on unramified cohomology groups.
\stepcounter{thm}
\begin{rem}\label{rem:ur}
\begin{enumerate}
\item[(1)]
For $n=0$, we have
\[ \HnrXz=\H^1(\cX\hspace{-2pt},\qzp). \]
If $k$ is global, then $\HnrXz$ is finite by a theorem of Katz-Lang \cite{KL}.
\item[(2)]
For $n=1$, we have
\[  \HnrXo = \Br(\cX)_{\ptor}. \]
If $k$ is global, the finiteness of $\HnrXo$ is equivalent to
 the finiteness of the Tate-Shafarevich group of the Picard variety of $\Xgf\,($cf.\ \cite{gr}, {\rm III}, \cite{Ta}$)$.
\item[(3)]
For $n=d:=\dim(\cX)$, $\HnrXd$ agrees with a group considered by Kato \cite{Ka3},
    who conjectured that
\[ \HnrXd=0 \;\; \hbox{ if $p\not=2$ or $k$ has no embedding into $\bR$}. \]
His conjecture is a generalization, to
    higher-dimensional proper arithmetic schemes, of the corresponding classical fact on
    the Brauer groups of local and global integer rings.
The $d=2$ case is proved in \cite{Ka3} and
   the $d=3$ case is proved in \cite{JS}.
\end{enumerate}
\end{rem}
We prove here that $\H^3_{\ur}(K,\qzp(2))$ is related with
  the torsion part of the cokernel of a cycle class map, assuming its finiteness.
This result will not be used in the rest of this paper, but shows an arithmetic meaning of $\H^3_{\ur}(K,\qzp(2))$.
See also Appendix B below for a zeta value formula for threefolds over finite fields using unramified cohomology.
\begin{prop}\label{prop5-0}
Assume that $\H^3_{\ur}(K,\qzp(2))$ is finite.
Then the order of
\[ \Coker \big(\varrho^2_{\zp} : \CH^2(\cX) \otimes \zp \lra \H^4(\cX\hspace{-2pt},\zptX) \big){}_{\ptor} \]
agrees with that of $\H^3_{\ur}(K,\qzp(2))$.
\end{prop}
\begin{pf}
By Lemma \ref{lem4-1}, $\H^4(\cX\hspace{-2pt},\zptX)$ is finitely generated over $\zp$,
 and $\Coker(\varrho^2_{\zp})_{\ptor}$ is finite.
Consider the following commutative diagram with exact rows (cf.\ Lemma \ref{lem4-1}):
\[ {\small
\begin{CD}
0 @.~ \ra ~@. \CH^2(\cX)_{\ptor}
   @.~ \ra ~@. \CH^2(\cX) \otimes \zp @.~ \ra ~@. \CH^2(\cX) \otimes \qp
   @.~ \ra ~@. \CH^2(\cX) \otimes \qzp @.~\ra ~@. 0\\
 @. @. @V{a}VV  @.  @V{\varrho^2_{\zp}}VV
   @. @V{\varrho^2_{\qp}}VV @. @V{\varrho^2_{\qzp}}VV \\
0 @.~ \ra ~@. \Cotor(\H^3(\cX\hspace{-2pt},\qzpX))
   @.~ \ra ~@. \H^4(\cX\hspace{-2pt},\zptX) @.~ \ra ~@. \H^4(\cX\hspace{-2pt},\qptX)
   @.~ \ra ~@. \H^4(\cX\hspace{-2pt},\qzpX),
\end{CD} } \]
where $a$ denotes the map obtained from the short exact sequence in Lemma \ref{lem3-1}\,(1).
See the proof of Lemma \ref{lem3-2} for the commutativity of the left square.
By the finiteness assumption on $\H^3_{\ur}(K,\qzp(2))$,
we see that
\[ \Coker(a) \simeq \gr_{\NF}^0\H^3(\cX\hspace{-2pt},\qzpX) :=
 \H^3(\cX\hspace{-2pt},\qzpX)/\Naa\H^3(\cX\hspace{-2pt},\qzpX) \]
(cf.\ Lemma \ref{lem3-1}\,(1))
and that the natural map
$\Ker(\varrho^2_{\qp}) \ra \Ker(\varrho^2_{\qzp})$ is zero (cf.\ \eqref{UC.0.1}).
Hence by a diagram chase on the above diagram, we obtain
   a short exact sequence
\[  0 \lra \gr_{\NF}^0\H^3(\cX\hspace{-2pt},\qzpX) \lra
 \Coker(\varrho^2_{\zp})_{\ptor} \lra \Ker(\varrho^2_{\qzp}) \lra 0. \]
Comparing this sequence with \eqref{UC.0.1}, we obtain the assertion.
\end{pf}
\newpage

\section{Finiteness of an unramified cohomology group}\label{sect4}
\medskip
\subsection{Finiteness of $\bs{\HBX}$}\label{sect5-1}
Let $k, S, p, \cX$ and $K$ be as in the notation \ref{intro}.8.
In this and the next section, we prove the following result,
  which implies the finiteness assertion on $\HBX$ in Theorem \ref{thm0-0}\,(1).
See \S\ref{sect2.1} for the condition {\bf H0}.
\begin{thm}\label{thm4-1}
Assume {\bf H0}, {\bf H1*} and either $p\geq 5$ or the equality
\begin{equation}\tag{$*_g$}\label{intro.1''}
\H^1_g(k,\H^2(\Xggf,\qzp(2)))_{\Div} =
    \H^1(k,\H^2(\Xggf,\qzp(2)))_{\Div}.
\end{equation}
Then $\HBX$ is finite.
\end{thm}
\noindent
In this section we reduce Theorem \ref{thm4-1}
   to Key Lemma \ref{keylem3-1} stated in \S\ref{sect5-4} below. We will prove the key lemma in \S\ref{sect5}.
We first prove Theorem \ref{cor0-1} admitting Theorem \ref{thm4-1}.
\par
\smallskip
\begin{pf*}
{\it Proof of Theorem \ref{cor0-1}}
The assumption $\H^2(\Xgf,\O_{\Xgf})=0$ implies {\bf H1*} and $(*_g)$
 (cf.\ Fact \ref{fact0-0}, Remark \ref{lem2-01}, Remark \ref{rem2-1}\,(1)).
Hence $\HBX$ is finite by Theorem \ref{thm4-1}.
By Corollary \ref{prop0-0},
   there is a positive integer $r_0$ such that
  \[ \ker(\varrho^2_{\ptor,r})=(\CH^2(\cX)_{\ptor})_{\Div} \quad \hbox{ for any $r \geq r_0$.} \]
Thus it remains to check that $\CH^2(\cX)_{\ptor}$ is finite,
   which follows from the finiteness of $\CH^2(\Xgf)_{\ptor}$
     (cf.\ Theorem \ref{thm2-1}) and \cite{CTR2}, Lemma 3.3.
This completes the proof.
\end{pf*}
\subsection{Proof of Theorem \ref{thm4-1}, Step 1}
We reduce Theorem \ref{thm4-1}
    to Proposition \ref{keyprop} below.
Let $\Naa\H^3(\Xgf,\qzp(2))$
   (resp.\ $\gr_{\NF}^0 \H^3(\Xgf,\qzp(2))$) be the kernel
     (resp.\ the image) of the natural map
$$
   \H^3(\Xgf,\qzp(2)) \lra \H^3(K,\qzp(2)).
$$
In view of Lemma \ref{lem3-1},
   there is a commutative diagram with exact rows
\begin{equation}\label{INJ.3.1}
\begin{CD}
   \Naa\H^3(\Xgf,\qzp(2)) @. ~ \hra~ @. \H^3(\Xgf,\qzp(2))
     @. ~ \ra ~ @. \gr_{\NF}^0 \H^3(\Xgf,\qzp(2)) \\
   @V{\delta_1}VV @.  @V{\delta_2}VV @. @V{\db}VV \\
   \us{v\in S^1}{\bigoplus} ~ \Nbb\H^4_{\Yv}(\cX\hspace{-2pt},\qzpX) @. \; \hra \; @.
   \us{v\in S^1}{\bigoplus} ~ \H^4_{\Yv}(\cX\hspace{-2pt},\qzpX) @. \; \ra \; @.
   \us{v\in S^1}{\bigoplus}~ \us{y\in \Yv^0}{\bigoplus}~
      \H^4_{y}(\cX\hspace{-2pt},\qzpX),
\end{CD}
\end{equation}
where the arrows $\delta_2$ and $\db$ arise from boundary maps in localization theory
   and $\delta_1$ is induced by the right square.
Note that we have
$$
\ker(\db) = \HBX.
$$
\stepcounter{thm}
\begin{prop}\label{keyprop}
Assume {\bf H0}, {\bf H1*} and either $p\geq 5$ or $(*_g)$.
Then $\Ker(\db)_{\Div}=0$.
\end{prop}
\noindent
The proof of this proposition
    will be started in \S\ref{sect4.3} below
           and finished in the next section.
We first finish the proof of Theorem \ref{thm4-1}, admitting Proposition 
\ref{keyprop}.
It suffices to show:
\begin{lem}\label{lem3-4}
Assume {\bf H1*} if $k$ is global.
Then $\Ker(\db)$ is cofinitely generated.
\end{lem}
\begin{pf}
The case that $k$ is local is obvious,
  because $\H^3(\Xgf,\qzp(2))$ is cofinitely generated.
Assume that $k$ is global.
We use the notation fixed in \ref{intro}.\GN.
By Lemma \ref{lem4-1}, $\H^3(\cX\hspace{-2pt},\qzpX)$ is cofinitely generated.
Hence it suffices to show $\Coker(\delta_1)$ is cofinitely generated,
    where $\delta_1$ is as in \eqref{INJ.3.1}.
There is a commutative diagram
\[ \xymatrix{ \CH^2(\Xgf,1)\otimes\qzp \ar[r]^{\partial\quad\;\;\;} \ar[d]
   & \bigoplus{}_{v\in S^1} \ \CH_{d-2}(\Yv)\otimes\qzp \ar[d]^{\hspace{-2pt}\wr} \\
\Naa\H^3(\Xgf,\qzp(2)) \ar[r]^{\delta_1\quad\;\;\;}
 & \bigoplus{}_{v \in S^1} \ \Nbb\H^4_{\Yv}(\cX\hspace{-2pt},\qzpX), } \]
where the right vertical isomorphism follows from Lemma \ref{lem3-1}\,(2)
    and $\partial$ is the boundary map of the localization sequence of higher Chow groups.
See \eqref{TOR.2.1} for the left vertical arrow.
Since $\Nbb\H^4_{\Yv}(\cX\hspace{-2pt},\qzpX)$ is cofinitely generated
     for any $v \in S^1$,
it suffices to show that
   for a sufficiently small non-empty open subset $U \subset S$,
       the cokernel of the boundary map
\[ \partial_U : \CH^2(\Xgf,1)\otimes\qzp
   \lra  \bigoplus_{v\in (U)^1}\ \CH_{d-2}(\Yv)\otimes\qzp \]
   is cofinitely generated.
Note that $\CH_{d-2}(\Yv)=\CH^1(\Yv)$
       if $\Yv$ is smooth.
Now let $U$ be a non-empty open subset of $S \sm \varSigma$
     for which $\cX \times _S U \ra U$ is smooth.
Put $A:=\H^2(\Xggf,\qzp(2))$, viewed as a smooth sheaf on $U_{\et}$.
There is a commutative diagram up to a sign
\[ \xymatrix{ \CH^2(\Xgf,1)\otimes\qzp \ar[r]^{\partial_U\quad} \ar[d]_{\reg_{\qzp}}
  & \bigoplus{}_{v\in U^1} \ \CH^1(\Yv) \otimes \qzp \ar[d]^{\tau_U} \\
   \H^1(k,A)  \ar[r]^{\delta_U\quad} & \bigoplus{}_{v\in U^1} \ A(-1)^{G_{\Fv}}.  } \]
See \S\ref{sect1.5} for $\delta_U$.
The right vertical arrow $\tau_U$ is defined as the composite map
\[ \CH^1(\Yv) \otimes\qzp \; \hra \; \H^2(\Yv,\qzp(1))
     \os{\epsilon}\lra
      \H^2(\Ynrv,\qzp(1))^{G_{\Fv}}=A(-1)^{G_{\Fv}}, \]
where the first injective map is the cycle class map
    for divisors on $\Yv$.
Note that $\Coker(\partial_U)$ is divisible and that
$\Ker(\tau_U)$ has a finite exponent by the isomorphism
\[ \Ker(\epsilon) \simeq \H^1(\Fv,\H^1(\Ynrv,\qzp(1)))
\simeq \H^1(\Fv,\Cotor(\H^1(\Xggf,\qzp(1)))) \]
for $v \in U^1$,
where the first isomorphism follows from the Hochschild-Serre spectral
sequence for $\Yv$.
Hence
   to prove that $\Coker(\partial_U)$ is cofinitely generated,
   it suffices to show that the map
\[ \partial':=\tau_U \circ \partial_U : \CH^2(\Xgf,1)\otimes\qzp
   \lra \bigoplus_{v\in U^1}\ \big(A(-1)^{G_{\Fv}}\big){}_{\Div} \]
has cofinitely generated cokernel
         (cf.\ Lemma \ref{lem1-3}\,(2)).
Finally since $\Coker(\reg_{\qzp})_{\Div}$
        is cofinitely generated by {\bf H1*}
          (cf.\ Lemma \ref{lem1-1}),
     $\partial'$ has cofinitely generated cokernel
         by Lemma \ref{lem1-2}\,(1).
Thus we obtain Lemma \ref{lem3-4}.
\end{pf}
\subsection{Proof of Theorem \ref{thm4-1}, Step 2}\label{sect4.3}
We construct a key commutative diagram \eqref{INJ.4.2} below
       and prove Lemma \ref{lem2-1b},
   which play key roles in our proof of Proposition \ref{keyprop}.
We need some preliminaries.
We suppose that $k$ is global
    until the end of Lemma \ref{lem2-2'}.
Let $\varSigma \subset S$ be the set of the closed points on $S$
       of characteristic $p$.
For non-empty open $U \subset S$, put
\[ \XU:=X \times_S U \qaq \XUp:=\XU \times_{S} {(S \sm \varSigma)}. \]
Let $j_U:\XUp \to \XU$ be the natural open immersion.
There is a natural injective map
\[ \xymatrix{
\alpha_{U,r}: \H^3(\XU,\tau_{\leq 2} Rj_{U*}\mupr 2)
    \; \ar@{^{(}->}[r] & \H^3(\XUp,\mupr 2) } \]
induced by the canonical morphism
$\tau_{\leq 2} Rj_{U*}\mupr 2\to Rj_{U*}\mupr 2$.
\begin{lem}\label{lem2-2'}
We have
$\Naa\H^3(\XUp,\mupr 2)\subset \Image(\alpha_{U,r}).$
\end{lem}
\begin{pf}
We compute the local-global spectral sequence
\[ E_1^{u,v}=\bigoplus_{x\in (\XU)^u} \
\H^{u+v}_x(\XU,\tau_{\leq 2} Rj_{U*}\mupr 2) \Lra
\H^{u+v}(\XU,\tau_{\leq 2} Rj_{U*}\mupr 2). \]
By the absolute cohomological purity \cite{Fu} and Lemma \ref{lem2-h}\,(1),
we have
\[ E_1^{u,v} \simeq \begin{cases}
       \H^v(K,\mupr 2)  & \quad \hbox{ (if }\; u=0)\\
        \displaystyle \bigoplus{}_{x\in (\XUp)^u}~
          \H^{v-u}(x,\mupr {2-u}) &\quad \hbox{ (if }\; v \leq 2).\\
                    \end{cases} \]
Repeating the same computation as in the proof of Lemma \ref{lem2-0},
     we obtain
\[   \Naa\H^3(\XUp,\mupr 2)
      \simeq E_2^{1,2}=E_{\infty}^{1,2} \hra
\H^3(\XU,\tau_{\leq 2} Rj_{U*}\mupr 2), \]
which completes the proof of Lemma \ref{lem2-2'}.
\end{pf}
\smallskip
Now we suppose that $k$ is either local or global, and
    define the group $\WW$ as follows:
\stepcounter{equation}
\begin{equation}\label{INJ.4.1}
    \WW := \begin{cases}
   \H^3(\Xgf,\qzp(2)) &
     \hbox{(if $k$ is $\ell$-adic local with }\ell \not = p)\\
   \H^3(\cX\hspace{-2pt},\tau_{\leq 2} Rj_*\qzp(2)) &
     \hbox{(if $k$ is $p$-adic local)} \\
   \us{\varSigma\subset U \subset S}\varinjlim \
   \H^3(\XU,\tau_{\leq 2} Rj_{U*}\qzp(2))
         & \hbox{(if $k$ is global),}
        \end{cases}
\end{equation}
where $j$ in the second case
 denotes the natural open immersion $\Xgf \hra \cX$, and
    the limit in the last case is taken over all non-empty open subsets $U \subset S$ which contain $\varSigma$.
By Lemma \ref{lem2-2} and Lemma \ref{lem2-2'},
    there are inclusions
\[ \Naa\H^3(\Xgf,\qzp(2)) \subset \WW \subset
        \H^3(\Xgf,\qzp(2)) \]
and a commutative diagram
\begin{equation}\label{INJ.4.2}
\xymatrix{
\NFaa\H^3(\Xgf,\qzp(2))_{\Div} \;\ar@{^{(}->}[r] \ar[rd]_{\cdiv} &\quad (\WW^0)_{\Div} \ar[d]^{\adiv \qquad }  \\
& \H^1(k,\H^2(\Xggf,\qzp(2))). }
\end{equation}
Here $\NFaa\H^3(\Xgf,\qzp(2))$ is as we defined in \S\ref{sect2.5},
    and we put
\begin{equation}\label{INJ.4.3}
   \WW^0 := \ker \big(\WW \lra \H^3(\Xggf,\qzp(2)) \big).
\end{equation}
The arrows $\adiv$ and $\cdiv$ are
     induced by the edge homomorphism \eqref{TOR.4.1}.
We show here the following lemma,
     which is stronger than Lemma \ref{lem2-1a}:
\addtocounter{thm}{3}
\begin{lem}\label{lem2-1b}
Assume either $p\geq 5$ or \eqref{intro.1''}.
Then $\Image(\adiv)\subset \H^1_g(k,\H^2(\Xggf,\qzp(2)))$.
\end{lem}
\begin{rem}
We will prove the equality
    $\Image(\adiv) = \H^1_g(k,\H^2(\Xggf,\qzp(2)))_{\Div}$
      under the same assumptions,
       later in Lemma $\ref{lem5-2}$.
\end{rem}
\noindent
The following corollary of Lemma \ref{lem2-1b}
      will be used later in \S\ref{sect5-4}:
\begin{cor}\label{cor2-1}
Assume {\bf H0}, {\bf H1*} and
      either $p\geq 5$ or \eqref{intro.1''}.
Then we have
\[ \Image(\cdiv)=\Image(\adiv)= \H^1_g(k,\H^2(\Xggf,\qzp(2)))_{\Div}. \]
\end{cor}
\begin{pf*}{\it Proof of Lemma $\ref{lem2-1b}$}
The assertion under the second condition is rather obvious.
In particular, we are done if $k$ is $\ell$-adic local with $\ell \not = p$
       (cf.\ Remark \ref{rem2-1}\,(1)).
If $k$ is $p$-adic local with $p \geq 5$,
    the assertion follows from Theorem \ref{thm1-1} and Lemma \ref{lem2-3}.
Before proving the global case, we show the following sublemma:
\begin{sublem}\label{lem2-4}
Let $k$ be an $\ell$-adic local field with $\ell \not =p$.
Let $\cX$ be a proper smooth scheme over $S:=\Spec(\OO_k)$.
Put $A :=\H^i(\Xggf,\qzp(n))$ and
\[ \HnrXkk {i+1} n := \Image \big(\H^{i+1}(\cX\hspace{-2pt},\qzp(n)) \to \H^{i+1}(\Xgf,\qzp(n))\big). \]
Then we have
\[ \H^1_f(k,A) \subset \Image\big(\Faa\H^{i+1}(\Xgf,\qzp(n)) \cap \HnrXkk {i+1} n \to \H^1(k,A)\big) \]
and the quotient is annihilated by $\#(A/A_{\Div})$,
    where $\F^{\bullet}$ denotes the filtration induced by
       the Hochschild-Serre spectral sequence \eqref{TOR.2.3}.
\end{sublem}
\begin{pf}
Put $\varLambda:=\qzp$, and let $\rf$ be the residue field of $k$.
By the proper smooth base change theorem,
   $G_k$ acts on $A$ through the quotient $G_\rf$.
It suffices to show the following two claims:
\begin{enumerate}
\item[(i)]
{\it We have
\[ \Image\big(\Faa\H^{i+1}(\Xgf,\varLambda(n)) \cap \H^{i+1}_{\ur}(\Xgf, \varLambda(n)) \ra \H^1(k,A) \big)
  =\Image\big(\H^1(\rf,A) \ra \H^1(k,A)\big). \] }
\item[(ii)]
{\it We have
\[ \H^1_f(k,A) \subset \Image(\H^1(\rf,A)\to \H^1(k,A)) \]
   and the quotient is annihilated by $\#(A/A_{\Div})$.}
\end{enumerate}
We show these claims.
Let $Y$ be the closed fiber of $\cX/S$, and
   consider a commutative diagram with exact rows
\[ \xymatrix{
0 \ar[r] & \H^1(\rf,\H^i(\Ynr,\varLambda(n))) \ar[r] \ar@{_{(}->}[d]_{\sigma_1}
 & \H^{i+1}(Y,\varLambda(n)) \ar[r] \ar[d]_{\sigma_2}
 & \H^{i+1}(\Ynr,\varLambda(n)))^{G_\rf} \ar[d]_{\wr\hspace{-1.5pt}}^{\sigma_3} \\
0 \ar[r] & \H^1(k,A) \ar[r]
 & \H^{i+1}(\Xgf,\varLambda(n))/\Fbb \ar[r]
 & \H^{i+1}(\Xggf,\varLambda(n)))^{G_k},
} \]
where the exactness of the upper (resp.\ lower) row
      follows from the fact that $\cd(G_\rf)=1$ (resp.\ $\cd(G_k)=2$).
The arrows $\sigma_1$ and $\sigma_3$ are induced by the isomorphism
    $\H^*(\Ynr,\varLambda(n)) \simeq \H^*(\Xggf,\varLambda(n))$
        (proper smooth base change theorem).
The arrow $\sigma_2$ is induced by 
\[ \sigma'_2 : \H^{i+1}(Y,\varLambda(n)) \lisom \H^{i+1}(\cX\hspace{-2pt},\varLambda(n))
     \lra \H^{i+1}(\Xgf,\varLambda(n)). \]
Since
   $\Image(\sigma'_2)=\H^{i+1}_{\ur}(\Xgf, \varLambda(n))$
    by definition, the claim (i) follows from the above diagram.
The second assertion immediately follows from the fact that
    $\H^1_f(k,A)=\Image(\H^1(\rf,A)_{\Div} \ra \H^1(k,A)).$
This completes the proof of Lemma \ref{lem2-4}.
\end{pf}
\medskip
We prove Lemma \ref{lem2-1b} in the case that $k$ is global with $p \geq 5$.
Let $\WW$ and $\WW^0$ be as in \eqref{INJ.4.1} and \eqref{INJ.4.3}, respectively, and put
\[ A:=\H^2(\Xggf,\qzp(2)).  \]
Note that $(\WW^0)_{\Div}=\WW_{\Div}$ by {\bf H0}.
By a similar argument as for Lemma \ref{lem1-1}, we have 
\[ \WW_{\Div} = \varinjlim_{\varSigma\subset U \subset S} \ \H^3(\XU,\tau_{\leq 2} Rj_{U*}\qzp(2))_{\Div}. \]
Here the limit is taken over all non-empty open subsets $U \subset S$ which contain $\varSigma$,
 and $j_U$ denotes the natural open immersion $\XUp \hra \XU$.
By this equality and the definition of $H^1_g(k,A)$ (cf.\ Definition \ref{def1-1}),
it suffices to show the following sublemma:
\begin{sublem}
Let $U$ be an open subset of $S$ containing $\varSigma$,
    and fix an open subset $U'$ of $U \sm \varSigma$ for which $\cX_{U'} \ra U'$ is smooth $($and proper$)$.
Put $\WW_U := \H^3(\XU,\tau_{\leq 2} Rj_{U*}\qzp(2))$.
Then for any $x \in (\WW_U)_{\Div}$, its diagonal image
\[ \xt = (\xt_v)_{v \in S^1}
  \in \prod_{v\in (U')^1} \ \H^1(k_v,A)/\H^1_f(k_v,A) \times \prod_{v\in S \sm U'} \ \H^1(k_v,A)/\H^1_g(k_v,A) \]
is zero.
\end{sublem}
\medskip \noindent
{\it Proof.}
Since $(\WW_U)_{\Div}$ is divisible,
   it suffices to show that $\xt$ is killed by a positive integer
   independent of $x$.
By Lemma \ref{lem2-4}, $\xt_v$ with $v\in (U')^1$
     is killed by $\#(A/A_{\Div})$.
Next we compute $\xt_v$ with $v \in \varSigma$.
Let $\Xv$ and $j_v : \Xkv \hra \Xv$ be as in \ref{intro}.\GN, and put
\[ \WW_v := \H^3(\Xv, \tau_{\leq 2}Rj_{v*}\qzp(2)). \]
By {\bf H0} over $k$, we have
\[ \Image \big( (\WW_U)_{\Div} \ra \WW_v \big)
 \subset \ker\big(
   \WW_v \ra \H^3(\Xggf,\qzp(2) \big){}_{\Div}. \]
Hence Theorem \ref{thm1-1} and Lemma \ref{lem2-3}
   imply that $\xt_v=0$ for $v\in \varSigma$.
Finally, because the product of the other components
\[ \prod_{v\in S \sm (U'\cup\varSigma)}\  \H^1(k_v,A)/\H^1_g(k_v,A) \]
is a finite group, we see that all local components
 of $\xt$ is annihilated by
 a positive integer independent of $x$.
This completes the proof of the sublemma and Lemma \ref{lem2-1b}.
\end{pf*}
\subsection{Proof of Theorem \ref{thm4-1}, Step 3}\label{sect5-4}
We reduce Proposition \ref{keyprop} to Key Lemma \ref{keylem3-1} below.
We replace the conditions
     in Proposition \ref{keyprop} with another condition
\par \vspace{3pt}
\begin{quote}
{\bf N1:}
{\it We have $\Image(\adiv)=\Image(\cdiv)$ in \eqref{INJ.4.2},
  and $\Coker\big(\reg_{\qzp}\big){}_{\Div}$
     is cofinitely generated over $\zp$.
Here $\reg_{\qzp}$ denotes the regulator map \eqref{TOR.1.1}.}
\end{quote}
\par \vspace{3pt} \noindent
Indeed, assuming {\bf H0}, {\bf H1*} and
      either $p\geq 5$ or \eqref{intro.1''}, we obtain {\bf N1},
 by Corollary \ref{cor2-1} and the fact that the quotient $\H(k,A)_{\Div}/\H^1_g(k,A)_{\Div}$,
 with $A=\H^2(\Xggf,\qzp(2))$, is cofinitely generated over $\zp$ for
       (cf.\ Lemma \ref{lem1-1}).
Thus Proposition \ref{keyprop} is reduced to the following:
\begin{keylem}\label{keylem3-1}
Assume {\bf H0} and {\bf N1}.
Then we have $\Ker(\db)_{\Div}=0$.
\end{keylem}
\noindent
This lemma will be proved in the next section.
\newpage

\section{Proof of the key lemma}\label{sect5}
\medskip
\subsection{Proof of Key Lemma \ref{keylem3-1}}
Let
\[ \deltad: \H^3(\Xgf,\qzp(2)) \lra
    \bigoplus_{v \in S^1} \,
    \bigoplus_{y\in (\Yv)^0} \ \H^4_{y}(\cX\hspace{-2pt},\qzpX) \]
be the map induced by $\db$ in \eqref{INJ.3.1}.
Put
\[ \varTheta :=
    \H^3(\Xgf,\qzp(2))\left/
      \big(\Naa\H^3(\Xgf,\qzp(2))_{\Div}\big)\right. \]
and let $\varThetaur \subset\varTheta$ be the image of $\Ker(\deltad)$.
Note that we have
\[ \varThetaur = \ker \Big(
     \varTheta \to \gr_{\NF}^0 \H^3(\Xgf,\qzp(2))
            \os{\db}{\to} \bigoplus{}_{v \in S^1}\,
                  \bigoplus{}_{y\in (\Yv)^0} \ \H^4_y(\cX\hspace{-2pt},\qzpX)
                   \Big) \]
and a short exact sequence
\[ 0 \lra \Cotor(\Naa\H^3(\Xgf,\qzp(2)))
     \lra \varThetaur \lra \Ker(\db) \lra 0. \]
If $k$ is global,
   the assumption of Proposition \ref{prop2-2}\,(2) is satisfied
      by the condition {\bf N1}.
Hence $\Cotor(\Naa\H^3(\Xgf,\qzp(2)))$ is finite
      in both cases $k$ is local and global
    (cf.\ Proposition \ref{prop2-2}, \eqref{TOR.2.1}).
By the above short exact sequence and Lemma \ref{lem1-3}\,(3),
    our task is to show
\[ \phantom{ \;{\bf H0}\; {\bf N1} \; } \varThetaur_{\Div}=0,
     \quad \hbox{assuming \; {\bf H0} \; and \; {\bf N1}}. \]
Let $\F^{\bullet}$ be the filtration on $\H^3(\Xgf,\qzp(2))$
   resulting from the Hochschild-Serre spectral sequence
      \eqref{TOR.2.3}.
We define the filtration $\F^{\bullet}$ on $\varTheta$ as that induced by $\F^{\bullet}\hspace{-1.5pt}\H^3(\Xgf,\qzp(2))$,
 and define the filtration $\F^{\bullet} \varThetaur \subset \varThetaur$ as the pull-back of $\F^{\bullet}\varTheta$.
Since {\bf H0} implies the finiteness of $\gr_\F^0\varThetaur$,
it suffices to show
\begin{equation}\label{KEY.0.1}
   (\F^1\varThetaur)_{\Div}=0,
        \quad \hbox{assuming \;{\bf N1}.}
\end{equation}
The following lemma will play key roles:
\stepcounter{thm}
\begin{lem}\label{lem4-2l}
Suppose that $k$ is local.
Then the following composite map has finite kernel$:$
\[ \tdelta : \H^2(k,\H^1(\Xggf,\qzp(2))
        \lra \H^3(\Xgf,\qzp(2)) \os{\deltad}{\lra}
  \bigoplus_{y\in Y^0} \ \H^4_y(\cX\hspace{-2pt},\qzpX). \]
Here the first map is obtained by the Hochschild-Serre spectral sequence
      \eqref{TOR.2.3} and the fact that $\cd(k)=2$
        $($cf.\ $\S\ref{sect1.6})$.
Consequently,
    the group $\Fbb\H^3(\Xgf,\qzp(2)) \cap \ker(\deltad)$ is finite.
\end{lem}
\noindent
Admitting this lemma,
   we will prove \eqref{KEY.0.1} in \S\S\ref{sect5.1}--\ref{sect5.2}.
We will prove Lemma \ref{lem4-2l} in \S\ref{sect5.3}.
\subsection{Case $\bs{k}$ is local}\label{sect5.1}
We prove \eqref{KEY.0.1}
   assuming that $k$ is local and that Lemma \ref{lem4-2l} holds.
Let $\rf$ be the residue field of $k$.
By Lemma \ref{lem4-2l},
    $\F^2\varThetaur$ is finite.
We prove that
    $\Image(\F^1 \varThetaur \ra \gr_F^1\varTheta)$ is finite,
       which is exactly the finiteness of $\gr_F^1 \varThetaur$
       and implies \eqref{KEY.0.1}.
Let $\WW$ and $\WW^0$ be as in \eqref{INJ.4.1} and \eqref{INJ.4.3},
     respectively.
{\bf N1} implies
\begin{equation}\label{KEY.1.1}
\grF 1 \varTheta  \simeq \Faa\H^3(\Xgf,\qzp(2))\big/ \big( (\WW^0)_{\Div}+ \Fbb\H^3(\Xgf,\qzp(2)) \big).
\end{equation}
If $p\not=\ch(\rf)$,
   then the group on the right hand side is clearly finite.
If $p=\ch(\rf)$, then Lemma \ref{lem4-1l}
below implies that the image of
$\F^1\varThetaur \ra \grF 1 \varTheta$ is a subquotient of $\Cotor(\WW^0)$,
which is finite by the proof of Lemma \ref{lem2-3}.
Thus we are reduced to
\stepcounter{thm}
\begin{lem}\label{lem4-1l}
If $p=\ch(\rf)$, then
$\Ker(\deltad) \subset \WW$.
\end{lem}
\begin{pf}
Let the notation be as in the notation \ref{intro}.\LN.
Recall that $\WW=\H^3(\cX\hspace{-2pt},\tau_{\leq 2}Rj_*\qzp(2))$ by definition.
There is a commutative diagram
   with distinguished rows in $D^b(\cX\hspace{-2pt},\pnuz)$
\[ \xymatrix{
\pnuX \ar[r] \ar[d]_t & Rj_*\mupnu 2 \ar[r] \ar@{=}[d]
 & Ri_*Ri^!\pnuX[1] \ar[r] \ar[d]_{Ri_*Ri^!(t)[1]} & \pnuX[1] \ar[d]_{t[1]} \\
\tau_{\leq 2} Rj_*\mupnu 2 \ar[r] & Rj_*\mupnu 2 \ar[r]
  & Ri_*Ri^!(\tau_{\leq 2} Rj_*\mupnu 2)[1] \ar[r] & (\tau_{\leq 2} Rj_*\mupnu 2)[1], }\]
where $t$ is as in (S5) in \S\ref{sect3.1}.
The central square of this diagram
   gives rise to the left square of the following commutative diagram
     (whose rows are not exact):
\stepcounter{equation}
\begin{equation}\label{KEY.1.2}
{\small
\begin{CD}
\H^3(\Xgf,\qzp(2)) @. \; \lra \; @. \H^4_Y(\cX\hspace{-2pt},\qzpX) @. \; \lra \; @.
    \us{y\in Y^0}{\bigoplus} \ \H^4_y(\cX\hspace{-2pt},\qzpX)\\
@|  @. @VVV @. @V{\lambda}VV \\
\H^3(\Xgf,\qzp(2)) @. \; \os{\epsilon_1}\lra \; @.
    \H^0(Y,i^*R^3j_*\qzp(2))
      @. \; \os{\epsilon_2}\lra \; @.
\us{y\in Y^0}{\bigoplus} \ \H^0(y,i^*R^3j_*\qzp(2)).
\end{CD}
}
\end{equation}
Here we have used the isomorphism
$\tau_{\geq 3} Rj_*\mupnu 2 \simeq Ri_*Ri^!\tau_{\leq 2} Rj_*\mupnu 2[1]$.
The composite of the upper row is $\deltad$.
We have
$\Ker(\epsilon_1)=\WW$ obviously,
      and $\epsilon_2$ is injective by the second part of Lemma \ref{lem2-h}.
Hence we have $\ker(\deltad) \subset
     \ker(\epsilon_2 \circ \epsilon_1) = \ker(\epsilon_1)=\WW$.
\end{pf}
\subsection{Case $\bs{k}$ is global}\label{sect5.2}
We prove \eqref{KEY.0.1}
   assuming that $k$ is global and that Lemma \ref{lem4-2l} holds.
Let $\WW$ and $\WW^0$ be as in \eqref{INJ.4.1} and \eqref{INJ.4.3}, respectively.
{\bf N1} implies
\begin{equation}\label{KEY.2.1}
\grF 1 \varTheta  \simeq \Faa\H^3(\Xgf,\qzp(2))\big/
    \big( (\WW^0)_{\Div}+ \Fbb\H^3(\Xgf,\qzp(2)) \big).
\end{equation}
We first prove the following lemma:
\stepcounter{thm}
\begin{lem}\label{lem4-1g}
$\Ker(\deltad) \subset \WW$.
\end{lem}
\begin{pf}
We use the notation in \ref{intro}.\GN.
By the same argument as for the proof of Lemma \ref{lem4-1l},
   we obtain a commutative diagram analogous to \eqref{KEY.1.2}
\[ {\small
\begin{CD}
\H^3(\Xgf,\qzp(2)) @.~\lra~@.
   \underset{v\in S^1}{\bigoplus} \ \H^4_{\Yv}(\cX\hspace{-2pt},\qzpX)
   @.~\lra~@.
   \underset{v\in S^1}{\bigoplus} \,
   \underset{y\in (\Yv)^0}{\bigoplus} \ \H^4_y(\cX\hspace{-2pt},\qzpX)\\
@| @.  @VVV @. @VVV \\
\H^3(\Xgf,\qzp(2)) @.~\os{\epsilon_1}\lra~@.
   \underset{v\in \varSigma}{\bigoplus} \ \H^0(Y_v,i_v^*R^3j_{v*}\qzp(2))
     @.~\os{\epsilon_2}{\lra}~@.
   \underset{v\in \varSigma}{\bigoplus} \, \underset{y\in (\Yv)^0}{\bigoplus} \
     \H^0(y,i_v^*R^3j_{v*}\qzp(2))\\
\end{CD} } \]
The composite of the upper row is $\deltad$.
The assertion follows from the facts that $\Ker(\epsilon_1)=\WW$
and that $\epsilon_2$ is injective (cf.\ Lemma \ref{lem2-h}).
\end{pf}
\par
We prove \eqref{KEY.0.1}.
By Lemma \ref{lem1-3}\,(4), it suffcies to show that
\[ (\F^2\varThetaur)_{\Div}=0 =(\gr_F^1 \varThetaur)_{\Div}. \]
Since $\Fbb\H^3(\Xgf,\qzp(2)) \cap \ker(\deltad)$ has a finite exponent
 (Corollary \ref{cor1-5}\,(2), Lemma \ref{lem4-2l}), we have $(\F^2\varThetaur)_{\Div}=0$.
We show $(\grF 1\varThetaur)_{\Div}=0$.
By \eqref{KEY.2.1} and Lemma \ref{lem4-1g}, we have
\[  \grF 1\varThetaur \subset \varXi:= \WW^0/((\WW^0)_{\Div}+Z) \quad
 \text{ with } \; Z := \WW^0 \cap \Fbb\H^3(\Xgf,\qzp(2)). \]
By Corollary \ref{cor1-5}\,(1), $\Cotor(Z)$ has a finite exponent, which implies
\[ (\grF 1 \varThetaur)_{\Div} \subset \varXi_{\Div} = \Cotor(\WW^0)_{\Div} =0 \]
(cf.\ Lemma \ref{lem1-3}\,(3)). Thus we obtain \eqref{KEY.0.1}.
\subsection{Proof of Lemma \ref{lem4-2l}}\label{sect5.3}
The case that $k$ is $p$-adic local follows from \cite{Sat1}, Theorem 3.1, Lemma 3.2\,(1) (cf.\ \cite{Ts3}).
More precisely, $\cX$ is assumed in \cite{Sat1}, \S3 to have strict semistable reduction,
     but one can remove the strictness assumption easily.
The details are left to the reader.
\par
We prove Lemma \ref{lem4-2l}
   assuming that $k$ is $\ell$-adic local with $\ell \not= p$.
Note that in this case $\cX/S$ may not have semistable reduction.
If $\cX/S$ has strict semistable reduction,
   then the assertion is proved in \cite{Sat1}, Theorem 2.1.
We prove the general case. Put
\[ \Lam := \qzp \]
for simplicity. By the alteration theorem of de Jong \cite{dJ},
    we take a proper generically finite morphism
    $f:\tX \ra X$ such that $\tX$ has strict semistable reduction over
      the normalization $\tS=\Spec(\OO_{\tk})$ of $S$ in $\tX$.
Note that $\tdelta$ is the composite of a composite map
\[ \xdelta : \H^2(k,\H^1(\Xggf,\Lam(2)) \lra
     \H^3(\Xgf,\Lam(2)) \os{\dloc}{\lra}
       \H^4_Y(\cX\hspace{-2pt},\Lam(2)) \]
with a pull-back map
\begin{equation}\label{KEY.4.1}
   \H^4_Y(\cX\hspace{-2pt},\Lam(2)) \lra \bigoplus_{y\in Y^0} \ \H^4_y(\cX\hspace{-2pt},\Lam(2)).
\end{equation}
Here the arrow $\dloc$ is the boundary map in localization theory. There is a commutative diagram
\[ \xymatrix{
\H^2(k,\H^1(\Xggf,\Lam(2))) \ar[r]^{\quad \xdelta} \ar[d]_{f^*} & \H^4_Y(\cX\hspace{-2pt},\Lam(2)) \ar[d]^{f^*} \\
\H^2(\tk,\H^1(\tXggf,\Lam(2))) \ar[r]^{\quad \xdelta '} & \H^4_{\tY}(\tX,\Lam(2)), } \]
where $\tXggf:=\tX \otimes_{\OO_{\tk}} \ol k$ and $\tY$ denotes the closed fiber of $\tX/\tS$.
We have already shown that $\Ker(\xdelta')$ is finite,
   and a standard norm argument shows that
      the left vertical arrow has finite kernel.
Thus $\Ker(\xdelta)$ is finite as well.
It remains to show
\stepcounter{thm}
\begin{lem}\label{lem4-3}
$\Image(\xdelta) \cap \Nbb\H^4_Y(\cX\hspace{-2pt},\Lam(2))$ is finite,
  where $\Nbb\H^4_Y(\cX\hspace{-2pt},\Lam(2))$
    denotes the kernel of the map \eqref{KEY.4.1}.
\end{lem}
\begin{pf}
First we note that
\[ \Image(\xdelta) \subset \Image \big( \H^1(\rf,\H^3_{\Ynr}(\Xnr,\Lam(2)))
      \ra \H^4_Y(\cX\hspace{-2pt},\Lam(2)) \big). \]
Indeed, this follows from the fact that $\xdelta$ factors as follows:
\begin{align*}
  & \H^2(k,\H^1(\Xggf,\Lam(2))) \simeq
     \H^1(\rf,\H^1(\knr,\H^1(\Xggf,\Lam(2)))) \\
  & \lra \H^1(\rf,\H^2(\Xknr,\Lam(2))) 
    \lra \H^1(\rf,\H^3_{\Ynr}(\Xnr,\Lam(2)))
   \lra \H^4_Y(\cX\hspace{-2pt},\Lam(2)).
\end{align*}
Hence it suffices to show the finiteness of the kernel of the composite map
\[ \upsilon : \H^1(\rf,\H^3_{\Ynr}(\Xnr,\Lam(2)))
    \lra \H^4_Y(\cX\hspace{-2pt},\Lam(2)) \lra \bigoplus_{y \in Y^0} \ \H^4_y(\cX\hspace{-2pt},\Lam(2)). \]
There is a commutative diagram with exact rows and columns
\[ {\small
\xymatrix{
    & \H^1(\rf,\H^3_{\Ynr}(\Xnr,\Lam(2))) \ar[d] \\
  \CH_{d-2}(Y)\otimes\Lam \; \ar@{^{(}->}[r] \ar[d]_{\iota}
   & \H^4_Y(\cX\hspace{-2pt},\Lam(2)) \ar[r] \ar[d]
   & \bigoplus{}_{y \in Y^0} \ \H^4_y(\cX\hspace{-2pt},\Lam(2)) \ar[d] \\
  \CH_{d-2}(\Ynr)\otimes\Lam \; \ar@{^{(}->}[r]
   & \H^4_{\Ynr}(\Xnr,\Lam(2)) \ar[r]
   & \bigoplus{}_{\eta \in (\ol Y)^0} \ \H^4_\eta(\Xnr,\Lam(2)),
}} \]
where the horizontal rows arise from the isomorphisms
\[ \Nbb\H^4_Y(\cX\hspace{-2pt},\Lam(2))\simeq \CH_{d-2}(Y)\otimes\Lam, \quad
   \Nbb\H^4_{\Ynr}(\Xnr,\Lam(2))\simeq \CH_{d-2}(\Ynr)\otimes\Lam \]
 with $d:=\dim(\cX)$ (cf.\ Lemma \ref{lem3-1}\,(2)).
The middle vertical exact sequence arises
from the Hochschild-Serre spectral sequence for $\Xnr/X$.
A diagram chase shows that $\Ker(\upsilon)\simeq \Ker(\iota)$,
    and we are reduced to showing the finiteness of $\Ker(\iota)$.
Because the natural restriction map
    $\CH_{d-2}(Y)/\CH_{d-2}(Y)_{\tor} \to \CH_{d-2}(\Ynr)/\CH_{d-2}(\Ynr)_{\tor}$
   is injective by the standard norm argument,
     the finiteness of $\Ker(\iota)$ follows from
the following general lemma:
\begin{lem}\label{lem4-4}
Let $e$ be a positive integer and
let $Z$ be a scheme which is
   separated of finite type over $F:=\ol {\rf}$ with
    $\dim(Z)\leq e$.
Then the group $\CH_{e-1}(Z)/\CH_{e-1}(Z)_{\tor}$
    is a finitely generated abelian group.
\end{lem}
\noindent
{\it Proof of Lemma $\ref{lem4-4}$.}
Obviously we may suppose that $Z$ is reduced.
We first reduced the problem to the case where $Z$ is proper.
Take a dense open immersion $Z \hookrightarrow Z'$ with $\ol Z'$ is proper.
Writing $Z''$ for $Z' \sm Z$,
there is an exact sequence
\[ \CH_{e-1}(Z'')\lra \CH_{e-1}(Z')\lra \CH_{e-1}(Z)\lra 0, \]
where $\CH_{e-1}(Z'')$ is finitely generated free abelian group
     because $\dim(Z'') \leq e-1$.
Let $f:\tZ \to Z$ be the normalization.
Since $f$ is birational and finite,
   one easily sees that the cokernel of
$f_*: \CH_{e-1}(\tZ) \to \CH_{e-1}(Z)$ is finite.
Thus we may assume $Z$ is a proper normal variety of dimension $e$ over $F$.
Since $F$ is algebraically closed,
    $Z$ has an $F$-rational point.
Now the theory of Picard functor (cf.\ \cite{Mu}, \S5)
   implies the functorial isomorphisms
$\CH_{e-1}(Z)\simeq \Pic_{Z/F}(F)$,
where $\Pic_{Z/F}$ denotes the Picard functor for $Z/F$.
This functor is representable
     by a group scheme and fits into the exact sequence of group schemes
\[ 0 \lra \Pic^{\tau}_{Z/F} \lra \Pic_{Z/F} \lra \NS_{Z/F}\lra 0, \]
where $\Pic^{\tau}_{Z/F}$ is quasi-projective over $F$ and
the reduce part of $\NS_{Z/F}$ is associated with a finitely generated
    abelian group.
Since $F$ is the algebraic closure of
     a finite field, the group $\Pic^{\tau}_{Z/F}(F)$ is torsion.
Lemma \ref{lem4-4} follows immediately from these facts.
\end{pf}
This completes the proof of Lemma \ref{lem4-3}, Lemma \ref{lem4-2l}
      and the key lemma \ref{keylem3-1}. \qed
\newpage

\section{Converse result}\label{sect6}
\medskip
\subsection{Statement of the result}
Let $\rf$ be a finite field, and let $Z$ be a proper smooth geometrically integral variety over $\rf$.
We say that {\it the Tate conjecture holds in codimension $1$ for $Z$}, if the \'etale cycle class map
\[ \CH^1(Z) \otimes \bQ_{\ell} \lra \H^2(Z \otimes_{\rf}\ol {\rf},\bQ_\ell(1))^{G_\rf} \]
is surjective for a prime number $\ell \ne \ch(\rf)$ (\cite{Ta}, \cite{Ta2}).
By \cite{Mi0}, Theorem 4.1,
this condition is equivalent to that the $\ell$-primary torsion part of the Grothendieck-Brauer group
 $\Br(Z)=\H^2(Z,\Gm)$ is finite for any prime number $\ell$ including $\ch(\rf)$. 

Let $k, S, p, \cX$ and $K$ be as in the notation \ref{intro}.8.
In this section, we prove the following result,
  which implies Theorem \ref{thm0-0}\,(2)
   (see \S\ref{sect2.1} for {\bf H0}):
\begin{thm}\label{thm6-1}
Assume {\bf H0} and either $p\geq 5$ or the equality
\begin{equation}\tag{$*_g$}\label{intro.1'''}
\H^1_g(k,\H^2(\Xggf,\qzp(2)))_{\Div} =
    \H^1(k,\H^2(\Xggf,\qzp(2)))_{\Div}.
\end{equation}
Assume further the following three conditions$:$
\begin{quote}
{\bf F1 :} $\CH^2(\Xgf)_{\ptor}$ is finite.
\par\noindent
{\bf F2 :} $\HBX$ is finite.
\par\noindent
{\bf T :}
The reduced part of every closed fiber of $\cX/S$ has simple normal crossings on $\cX$,
   and the Tate conjecture holds in codimension $1$
        for the irreducible components of those fibers.
\end{quote}
Then {\bf H1*} holds.
\end{thm}
\subsection{Proof of Theorem \ref{thm6-1}}
Let
\[ \deltad: \H^3(\Xgf,\qzp(2)) \lra
\bigoplus_{v\in S}\, \bigoplus_{y\in (\Yv)^0}
   \ \H^4_y(\cX\hspace{-2pt},\qzpX). \]
be the map induced by $\db$ in \eqref{INJ.3.1}.
Let $\WW$ be as in \eqref{INJ.4.1}.
We need the following lemma:
\begin{lem}\label{lem5-1}
Assume that {\bf T} holds.
Then we have
\[ \WW_{\Div} \subset \Ker(\deltad)+\Fbb\H^3(\Xgf,\qzp(2)). \]
\end{lem}
\noindent
This lemma will be proved in \S\S\ref{sect6.2}--\ref{sect6.3} below.
We first finish the proof of Theorem \ref{thm6-1}, admitting Lemma \ref{lem5-1}.
The assumption {\bf F1} implies
   \[ \CH^2(\Xgf,1)\otimes\qzp = \Naa\H^3(\Xgf,\qzp(2))_{\Div}. \]
The assumption {\bf F2} implies the equality
   \[ \Naa\H^3(\Xgf,\qzp(2))=\ker(\deltad) \] up to a finite group.
Hence by Lemma \ref{lem5-1}, {\bf F1} and {\bf F2}
   imply the equality
\[ \Image(\reg_{\qzp})  = \Image(\adiv), \]
   where $\adiv$ is as in \eqref{INJ.4.2}.
Thus we are reduced to the following lemma stronger than
    Lemma \ref{lem2-1b}:
\begin{lem}\label{lem5-2}
Assume either $p \geq 5$ or $(*_g)$.
Then we have
\[ \Image(\adiv) = \H^1_g(k,\H^2(\Xggf,\qzp(2)))_{\Div}. \]
\end{lem}
\begin{pf}
If $k$ is local, then the assertion follows from
    Theorem \ref{thm1-1} and Lemma \ref{lem2-3}.
We show the inclusion $\Image(\adiv) \supset \H^1_g(k,\H^2(\Xggf,\qzp(2)))_{\Div}$,
 assuming that $k$ is global (the inclusion in the other direction has been proved in Lemma \ref{lem2-1b}).
Let $\WW^0$ be as in \eqref{INJ.4.3} and
    put $A:=\H^2(\Xggf,\qzp(2))$.
By Lemma \ref{lem1-3}\,(3),
   it is enough to show the following:
\begin{enumerate}
\item[(i)]
{\it The image of the composite map
\[ \xymatrix{
    \WW^0 \; \ar@{^{(}->}[r] & \Faa\H^3(\Xgf,\qzp(2)) \ar[r]^{\qquad \psi} & \H^1(k,A) } \]
contains $\H^1_g(k,A)_{\Div}$, where the arrow $\psi$ is as in \eqref{TOR.4.1}. }
\item[(ii)]
{\it The kernel of this composite map is cofinitely generated up to a group of finite exponent.}
\end{enumerate}
(ii) follows from Corollary \ref{cor1-5}\,(1).
We prove (i) in what follows. We use the notation fixed in \ref{intro}.\GN.
Let $U\subset S$ be a non-empty open subset which contains $\varSigma$
     and for which $\cX_U \ra U$ is smooth outside of $\varSigma$.
Let $j_U : \XUp \ra \XU$ be the natural open immersion.
Put $U':= U \sm \varSigma$ and $\Lam:=\qzp$.
For $v \in \varSigma$, put
\[ M_v  :=  \Faa\H^3(\Xkv,\Lam(2))
    / (\H^3(\Xv,\tau_{\leq 2}R(\jv)_*\Lam(2))^0)_{\Div}, \]
where the superscript $0$ means the subgroup of elements
    which vanishes in $\H^3(\Xggf,\Lam(2)) \simeq \H^3(\Xkv \otimes_{\kv} \ol {\kv},\Lam(2))$.
We construct a commutative diagram with exact rows
\[\xymatrix{
0 \ar[r] & \ker(r_U) \ar[r] \ar[d]_{c_U}  
 & \Faa\H^3(\XUp,\varLambda(2))  \ar[r]^{\qquad r_U} \ar[d]_{\psi_U}
  & \bigoplus{}_{v \in \varSigma}\; M_v \ar[d]_{b_{\varSigma}} \\
0 \ar[r] & \ker(a_U) \ar[r] & \H^1(U,A) \ar[r]^{a_U\qquad}
  & \bigoplus{}_{v \in \varSigma}\; \H^1_{/g}(k_v,A), } \]
where $\Faa$ on $\H^3(\XUp,\varLambda(2))$ means the filtration resulting from the
  Hochschild-Serre spectral sequence \eqref{TOR.1.2} for $\XUp$, and
 $\psi_U$ is an edge homomorphism of that spectral sequence.
The arrows $r_U$ and $a_U$ are natural pull-back maps, and we put
\begin{align*}
 \H^1_{/g}(k_v,A) & :=\H^1(k_v,A)/\H^1_g(k_v,A).
\end{align*}
The existence of $b_{\varSigma}$ follows from the local case of Lemma \ref{lem5-2},
 and $c_U$ denotes the map induced by the right square.
Note that $\ker(a_U)$ contains $\H^1_{f,U}(k,A)$.
Now let
\[ c : \WW^\dagger :=\us{\varSigma \subset U \subset S}{\varinjlim}\ \ker(r_U)
 \lra \us{\varSigma \subset U \subset S}{\varinjlim}\ \ker(a_U) \]
be the inductive limit of $c_U$, where $U$ runs through all non-empty open subsets of $S$
      which contains $\varSigma$ and for which $\cX_U \ra U$ is smooth
         outside of $\varSigma$.
Because the group on the right hand side contains $\H^1_{g}(k,A)$,
 it remains to show that
\begin{enumerate}
\item[(iii)]
{\it  $\Coker(c)$ has a finite exponent. }
\item[(iv)]
{\it $\WW^\dagger$ is contained in $\WW^0$.}
\end{enumerate}
(iv) is rather straight-forward and left to the reader.
We prove (iii). For $U \subset S$ as above,
applying the snake lemma to the above diagram, we see that the kernel of the natural map
\[  \Coker(c_U) \lra \Coker(\psi_U) \]
is a subquotient of $\ker(b_\varSigma)$. By the local case of Lemma \ref{lem5-2}, we have
\[ \ker(b_{\varSigma}) \simeq
    \bigoplus_{v \in \varSigma} \ \Image(\Fbb\H^3(\Xkv,\Lam(2)) \ra M_v) \]
    and the group on the right hand side is finite
      by Lemma \ref{lem5-3} below.
On the other hand, $\Coker(\psi_U)$ is zero if $p \geq 3$, and killed by $2$ if $p=2$.
Hence passing to the limit, we see that $\Coker(c)$ has a finite exponent.
This completes the proof of Lemma \ref{lem5-2}.
\end{pf}
\subsection{Proof of Lemma \ref{lem5-1}, Step 1}\label{sect6.2}
We start the proof of Lemma \ref{lem5-1}.
Our task is to show the inclusion
\begin{equation}\label{UC.3.1}
\deltad ( \WW_{\Div})
     \subset \deltad (\Fbb\H^3(\Xgf,\qzp(2)) ).
\end{equation}
If $k$ is global, then
   the assertion is reduced to the local case, because the natural map
\[ \Fbb\H^3(\Xgf,\qzp(2))) \lra \bigoplus_{v\in S^1} \
    \Fbb\H^3(\Xkv,\qzp(2))) \]
has finite cokernel by Corollary \ref{cor1-5}\,(2).
\par
Assume now that $k$ is local.
In this subsection,
    we treat the case that
      $k$ is $\ell$-adic local with $\ell \not = p$.
We use the notation fixed in \ref{intro}.\LN.
Recall that
     $Y$ has simple normal crossings on $\cX$ by the assumption {\bf T}.
Note that $\deltad$ factors as
\[ \H^3(\Xgf,\qzp(2)) \lra \H^4_Y(\cX\hspace{-2pt},\qzp(2)) \os{\iota}\lra
      \bigoplus_{y\in Y^0}\ \H^4_y(\cX\hspace{-2pt},\qzp(2)), \]
and that $\Image(\deltad) \subset \Image(\iota)$.
There is an exact sequence
\[ 0 \to \H^1(\rf,\H^3_{\Ynr}(\Xnr,\qzp(2)))
     \to \H^4_{Y}(\cX\hspace{-2pt},\qzp(2)) \to
          \H^4_{\Ynr}(\Xnr,\qzp(2))^{G_{\rf}} \to  0 \]
arising from a Hochschild-Serre spectral sequence.
We have
    $\Ker(\iota) \simeq \CH_{d-2}(Y)\otimes\qzp$ with $d:=\dim(\cX)$
      by Lemma \ref{lem3-1}\,(2).
Hence to show the inclusion \eqref{UC.3.1},
    it suffices to prove
\stepcounter{thm}
\begin{prop}\label{claim5-1}
\begin{enumerate}
\item[(1)]
Assume that {\bf T} holds.
Then the composite map
\stepcounter{equation}
\begin{equation}\label{UC.3.0}
\CH_{d-2}(Y)\otimes\qzp \lra \H^4_{Y}(\cX\hspace{-2pt},\qzp(2)) \lra
          \H^4_{\Ynr}(\Xnr,\qzp(2))^{G_{\rf}}
\end{equation}
is an isomorphism up to finite groups.
Consequently, we have
\[ \Image(\iota) \simeq \H^1(\rf,\H^3_{\Ynr}(\Xnr,\qzp(2))) \]
up to finite groups.
\item[(2)]
The image of the composite map
\[ \H^2(k,\H^1(\Xgf,\qzp(2))) \lra \H^3(\Xgf,\qzp(2)) \lra \H^4_{Y}(\cX\hspace{-2pt},\qzp(2)) \]
contains $\H^1(\rf,\H^3_{\Ynr}(\Xnr,\qzp(2)))_{\Div}$.
\end{enumerate}
\end{prop}
\noindent
We first show the following lemma:
\stepcounter{thm}
\begin{lem}\label{lem-inv}
\begin{enumerate}
\item[(1)]
Consider the Mayer-Vietoris spectral sequence
\stepcounter{equation}
\begin{equation}\label{UC.3.3}
   E_1^{u,v}=\H^{2u+v-2}(\Ynr{}^{(-u+1)},\qzp(u+1)) \Lra
\H^{u+v}_{\Ynr}(\Xnr,\qzp(2))
\end{equation}
obtained from the absolute purity
    $($cf.\ \cite{RZ}, \cite{Th}, \cite{Fu}$)$,
where $\Ynr{}^{(q)}$ denotes
      the disjoint union of $q$-fold intersections of
        distinct irreducible components of
          the reduced part of $\Ynr$.
Then there are isomorphisms up to finite groups
\begin{align*}
   \H^1(\rf,\H^3_{\Ynr}(\Xnr,\qzp(2))) & \simeq  \H^1(\rf,E_2^{-1,4}),\\
   \H^4_{\Ynr}(\Xnr,\qzp(2))^{G_{\rf}} & \simeq (E_2^{0,4})^{G_{\rf}}.
\end{align*}
\item[(2)]
As $G_{\rf}$-module,
   $\H^0(\knr,\H^2(\Xggf,\qp))$ has weight $\leq 2$.
\end{enumerate}
\end{lem}
\begin{pf*}
{\it Proof of Lemma \ref{lem-inv}}
(1)
Since $E_1^{u,v}=0$ for any $(u,v)$ with $u>0$ or $2u+v < 2$,
   there is a short exact sequence
\begin{equation}\label{UC.3.2}
   0 \lra E_2^{0,3} \lra \H^3_{\Ynr}(\Xnr,\qzp(2))) \lra E_2^{-1,4} \lra 0,
\end{equation}
and the edge homomorphism
\begin{equation}\label{UC.3.2'}
\xymatrix{  E_2^{0,4} \; \ar@{^{(}->}[r] & \H^4_{\Ynr}(\Xnr,\qzp(2))), }
\end{equation}
where we have $E_2^{-1,4}=\Ker(d_1^{-1,4})$
   and $E_2^{0,4}=\Coker(d_1^{-1,4})$
    and $d_1^{-1,4}$ is the Gysin map
       $\H^0(\Ynr{}^{(2)},\qzp) \ra \H^2(\Ynr{}^{(1)},\qzp(1)))$.
Note that $E_1^{u,v}$ is pure of weight $v-4$
       by Deligne's proof of the Weil conjecture \cite{De},
    so that $\H^i(\rf,E^{u,v}_{\infty})$ ($i=0,1$) is finite unless $v=4$.
The assertions immediately follow from these facts.
\par
(2)
By the alteration theorem of de Jong \cite{dJ},
    we may assume that $\cX$ is projective and has semistable reduction over $S$.
If $\Xgf$ is a surface, then the assertion is proved in \cite{RZ}.
Otherwise, take a closed immersion $\cX \hra \P^N_S=: \P$.
By \cite{JS}, Proposition 4.3\,(b), there exists a hyperplane $H \subset \P$
 which is flat over $S$ and for which ${\mathscr Z}:=\cX \times_{\P} H$ is regular with semistable reduction over $S$.
The restriction map
    $\H^2(\Xggf,\qp) \to \H^2(\Zggf,\qp)$\,($\Zggf:={\mathscr Z} \times_{\OO_k} \ol k$)
     is injective by the weak and hard Lefschetz theorems.
Hence the claim is reduced to the case of surfaces.
This completes the proof of the lemma.
\end{pf*}
\begin{pf*}{\it Proof of Proposition \ref{claim5-1}}
(1) Note that the composite map \eqref{UC.3.0} in question has finite kernel
     by Lemma \ref{lem4-4} and the arguments in the proof of Lemma \ref{lem4-3}.
We prove that \eqref{UC.3.0} has finite cokernel, assuming {\bf T}.
By the Kummer theory, there is a short exact sequence
\begin{equation}\label{UC.3.2''}
   0 \lra \Pic(\Ynr{}^{(1)})\otimes\qzp  \lra
    \H^2(\Ynr{}^{(1)},\qzp(1)) \lra \Br(\Ynr{}^{(1)})_{\ptor} \lra 0
\end{equation}
and the differential map $d_1^{-1,4}$
     of the spectral sequence \eqref{UC.3.3} factors through the Gysin map
\[ \H^0(\Ynr{}^{(2)},\qzp) \lra \Pic(\Ynr{}^{(1)})\otimes\qzp, \]
whose cokernel is $\CH_{d-2}(\Ynr)\otimes\qzp$.
Hence in view of the computations in the proof of Lemma \ref{lem-inv}\,(1),
   the Gysin map
     $\CH_{d-2}(\Ynr)\otimes\qzp \ra \H^4_{\Ynr}(\Xnr,\qzp(2))$
        (cf.\ Lemma \ref{lem3-1}\,(2))
         factors through the map \eqref{UC.3.2'} and
we obtain a commutative diagram
\[\xymatrix{
\CH_{d-2}(Y)\otimes\qzp \ar[r]^{\eqref{UC.3.0}\quad} \ar[d]
 & \H^4_{\Ynr}(\Xnr,\qzp(2))^{G_{\rf}} \\
(\CH_{d-2}(\Ynr)\otimes\qzp)^{G_{\rf}} \ar[r] & \;\; (E_2^{0,4})^{G_{\rf}}, \ar[u]_{\eqref{UC.3.2'}} } \]
where the left vertical arrow
   has finite cokernel (and kernel) by Lemma \ref{lem4-4}
         and a standard norm argument,
      and the right vertical arrow has finite cokernel
       (and is injective) by Lemma \ref{lem-inv}\,(1).
Thus it suffices to show that
     the bottom horizontal arrow has finite cokernel.
By the exact sequence \eqref{UC.3.2''},
   we obtain a short exact sequence
\[ 0 \lra \CH_{d-2}(\Ynr)\otimes\qzp \lra E_2^{0,4} \lra \Br(\Ynr{}^{(1)})_{\ptor} \lra 0. \]
Our task is to show that
   $\big(\Br(\Ynr{}^{(1)})_{\ptor}\big){}^{G_{\rf}}$ is finite,
      which follows from the assumption {\bf T}
        and the finiteness of the kernel of the natural map
\[ \H^1(G_{\rf},\Pic(\Ynr{}^{(1)})\otimes\qzp)  \lra \H^1(G_{\rf},\H^2(\Ynr{}^{(1)},\qzp(1))) \]
(cf.\ Lemma \ref{lem:app} in \S\ref{sect6.4} below).
Thus we obtain the assertion.
\par
(2) Since $\cd(\knr)=1$, there is a short exact sequence
\[ 0 \ra \H^1(\knr,\H^1(\Xggf,\qzp(2))) \ra \H^2(\Xknr,\qzp(2)) \ra
          \H^2(\Xggf,\qzp(2))^{G_{\knr}} \ra 0 \]
arising from a Hochschild-Serre spectral sequence.
By Lemma \ref{lem-inv} the last group has weight $\leq -2$,
   and we have isomorphisms up to finite groups
\begin{equation}\label{UC.3.4}
\begin{CD}
\H^2(k,\H^1(\Xggf,\qzp(2)))
   @. \; \simeq \; @. \H^1(\rf,\H^1(\knr,\H^1(\Xggf,\qzp(2))))\\
   @. \simeq @. \H^1(\rf,\H^2(\Xknr,\qzp(2))).
       \phantom{AAAi}
\end{CD}
\end{equation}
Now we plug the short exact sequence \eqref{UC.3.2}
     into the localization exact sequence
\[ \H^2(\Xknr,\qzp(2)) \lra \H^3_{\Ynr}(\Xnr,\qzp(2)))
      \os{\alpha}\lra \H^3(\Xnr,\qzp(2)). \]
Note that $\H^3(\Xnr,\qzp(2))\simeq \H^3(\Ynr,\qzp(2))$, so that
   it has weight $\leq -1$ (cf.\ \cite{De}).
Let $E_2^{u,v}$ be as in \eqref{UC.3.3}.
Since $E_2^{-1,4}$ is pure of weight $0$,
the induced map
\[ E_2^{-1,4} \lra \H^3(\Xnr,\qzp(2)))/\alpha(E_2^{0,3}) \]
has finite image.
Hence the composite map
\[ \H^2(\Xknr,\qzp(2)) \lra \H^3_{\Ynr}(\Xnr,\qzp(2))) \lra E_2^{-1,4} \]
    has finite cokernel, and the following map has finite cokernel as well:
\[ \H^1(\rf,\H^2(\Xknr,\qzp(2))) \lra  \H^1(\rf,E_2^{-1,4}). \]
Now Proposition \ref{claim5-1}\,(2) follows from this fact together with
      \eqref{UC.3.4} and the first isomorphism in Lemma \ref{lem-inv}\,(1).
\end{pf*}
\addtocounter{thm}{5}
\begin{rem}\label{rem5-2}
Let $J$ be the set of
   the irreducible components of $Y{}^{(2)}$ and put
\[ \Delta:=\Ker\big(g': \Z^J \to \NS(Y{}^{(1)})\big) \quad
     \hbox{ with } \quad \NS(Y{}^{(1)}):=
       \bigoplus_{y\in Y^0} \ \NS(Y_y), \]
where for $y \in Y^0$,
     $Y_y$ denotes the closure $\ol{\{y\}}\subset Y$
      and $\NS(Y_y)$ denotes its N\'eron-Severi group.
The arrow $g'$ arises from the Gysin map $\Z^J \to \Pic(Y{}^{(1)})$.
One can easily show, assuming {\bf T}
     and using Lemma $\ref{lem:app}$ in $\S\ref{sect6.4}$ below,
    that the corank of $\H^1(\rf,E_2^{-1,4})$
     over $\zp$ is equal to
    the rank of $\Delta$ over $\Z$.
Hence Proposition $\ref{claim5-1}\,(2)$ implies the inequality
\stepcounter{equation}
\begin{equation}\label{UC.3.5}
\dim_{\qp}(\H^2(k,\H^1(\Xggf,\qp(2)))\geq \dim_{\Q}(\Delta\otimes\Q),
\end{equation}
which will be used in the next subsection.
\end{rem}
\subsection{Proof of Lemma \ref{lem5-1}, Step 2}\label{sect6.3}
We prove Lemma \ref{lem5-1},
    assuming that $k$ is $p$-adic local
      (see \ref{intro}.\LN~ for notation).
We first show the following lemma:
\begin{lem}\label{lem5-3}
We have
\[ \Fbb\H^3(\Xgf,\qzp(2)) \subset \H^3(\cX\hspace{-2pt},\tau_{\leq 2}Rj_*\qzp(2)) \;(=\WW). \]
\end{lem}
\begin{pf*}
{\it Proof of Lemma \ref{lem5-3}}
By (S5) in \S\ref{sect3.1}, there is a distinguished triangle
  in $D^b(\cX_{\et},\Z/p^r\Z)$
\stepcounter{equation}
\begin{equation}\label{UC.4.1}
\begin{CD}
i_*\nu_{Y,r}^1[-3] @>{g}>>
    \pnuX @>{t'}>> \tau_{\leq 2}Rj_*\mupr 2 @>>> i_*\nu_{Y,r}^1[-2].
\end{CD}
\end{equation}
Applying $Ri^!$ to this triangle,
    we obtain a distinguished triangle in $D^b(Y_{\et},\Z/p^r\Z)$
\begin{equation}\label{UC.4.2}
\begin{CD}
\nu_{Y,r}^1[-3] @>{\gys_i^2}>>
    Ri^!\pnuX @>{Ri^!(t)}>>
      i^*(\tau_{\geq 3}Rj_*\mupr 2)[-1] @>>> \nu_{Y,r}^1[-2],
\end{CD}
\end{equation}
where $\gys_i^2:=Ri^!(g)$
    and we have used the natural isomorphism
\[ i^*(\tau_{\geq 3}Rj_*\mupr 2)[-1] \simeq
         Ri^!(\tau_{\leq 2}Rj_*\mupr 2). \]
Now let us recall the commutative diagram \eqref{KEY.1.2}:
\[ \xymatrix{
\H^3(\Xgf,\qzp(2)) \ar[r] \ar@{=}[d] & \H^4_Y(\cX\hspace{-2pt},\qzpX) \ar[r] \ar[d]
   & \bigoplus{}_{y\in Y^0} \ \H^4_y(\cX\hspace{-2pt},\qzpX) \ar[d]_{\lambda} \\
\H^3(\Xgf,\qzp(2)) \ar[r]^{\epsilon_1\quad} & \H^0(Y,i^*R^3j_*\qzp(2)) \ar[r]^{\epsilon_2\qquad}
 & \bigoplus{}_{y\in Y^0}\ \H^0(y,i^*R^3j_*\qzp(2)).
} \]
where the middle and the right vertical arrows are
      induced by $Ri^!(t)$ in \eqref{UC.4.2}.
By the proof of Lemma \ref{lem4-1l},
   we have
    $\H^3(\cX\hspace{-2pt},\tau_{\leq 2}Rj_*\qzp(2))=\Ker(\epsilon_2\epsilon_1)$.
Hence it suffices to show the image of the composite map
\[ \tdelta : \H^2(k,\H^1(\Xggf,\qzp(2)) \lra
           \H^3(\Xgf,\qzp(2)) \os{\deltad}\lra
           \bigoplus_{y\in Y^0} \ \H^4_y(\cX\hspace{-2pt},\qzpX) \]
is contained in $\Ker(\lambda)$.
By the distinguished triangle \eqref{UC.4.2},
     $\Ker(\lambda)$ agrees with the image of the Gysin map
\[ \bigoplus_{y\in Y^0} \ \gys_{i_y}^2 :
    \bigoplus_{y\in Y^0} \ \H^1(y,\logwitt y \infty 1) \; \hra \;
    \bigoplus_{y\in Y^0} \ \H^4_y(\cX\hspace{-2pt},\qzpX).  \]
On the other hand,
   as is seen in \S\ref{sect5.3},
$\tdelta$ factors through the maps
\begin{align*}
\H^1(\rf,\H^0(\Ynr,\mlogwitt {\Ynr} {\infty} 1 )) & \lra
\H^1\Big(
   \rf,{\bigoplus}{}_{\eta \in \Ynr{}^0}~ \H^0(\eta,\logwitt \eta \infty 1)
      \Big) \\
  & \lra \bigoplus{}_{y\in Y^0}~\H^1(y,\logwitt y \infty 1).
\end{align*}
Thus we obtain the assertion.
\end{pf*}
We start the proof of Lemma \ref{lem5-1}, i.e.,
    the inclusion \eqref{UC.3.1}, assuming that $k$ is $p$-adic local.
The triangle \eqref{UC.4.1}
   gives rise to the upper exact row
   of the following diagram whose
      left square is commutative and
       whose right square is anti-commutative:
\[ \xymatrix{
\H^3(\cX\hspace{-2pt},\pnuX) \ar[r] \ar@{=}[d]
 & \H^3(\cX\hspace{-2pt},\tau_{\leq 2}Rj_*\mupr 2) \ar[r] \ar[d]
 & \H^1(Y,\nu_{Y,r}^1) \ar[d]^{\gys_i^2} \\
\H^3(\cX\hspace{-2pt},\pnuX) \ar[r] & \H^3(\Xgf,\mupr 2) \ar[r] & \H^4_Y(\cX\hspace{-2pt},\pnuX),
} \]
where $\gys_i^2$ is as in \eqref{UC.4.2}
    and the anti-commutativity of the right square follows from
       (S4) in \S\ref{sect3.1}.
Hence the map $\deltad$ restricted to
    $\WW=\H^3(\cX\hspace{-2pt},\tau_{\leq 2}Rj_*\qzp(2))$ factors as
\[ \WW \lra \H^1(Y,\nu^1_{Y,\infty}) \lra
    \H^4_Y(\cX\hspace{-2pt},\qzpX) \lra \bigoplus_{y\in Y^0}\ \H^4_y(\cX\hspace{-2pt},\qzpX), \]
where $\nu_{Y,\infty}^1:=\varinjlim_r~ \nu_{Y,r}^1$.
By Lemmas \ref{lem5-3} and \ref{lem4-2l},
   it suffices to show that the corank of
\begin{equation}\label{UC.4.3}
   \Image\Big(\H^1(Y,\nu^1_{Y,\infty})  \to \bigoplus{}_{y\in Y^0} \
    \H^4_y(\cX\hspace{-2pt},\qzpX)\Big)
\end{equation}
is not greater than $\dim_{\qp}\H^2(k,\H^1(\Xggf,\qp(2)))$.
We pursue an analogy to the case $p\not=\ch(\rf)$
    by replacing $\H^4_Y(\cX\hspace{-2pt},\qzp(2))$ with $\H^1(Y,\nu_{Y,\infty}^1)$.
There is an exact sequence
\[ 0 \lra \H^1(\rf,\H^0(\Ynr,\nu^1_{\ol Y,\infty}))
     \lra \H^1(Y,\nu^1_{Y,\infty}) \lra
         \H^1(\Ynr,\nu^1_{\ol Y,\infty})^{G_{\rf}} \lra 0 \]
arising from a Hochschild-Serre spectral sequence.
By \cite{sato1}, Corollary 1.5,
    there is a Mayer-Vietoris spectral sequence
\[ E_1^{a,b} =\H^{a+b}(\Ynr{}^{(1-a)},\logwitt {\Ynr{}^{(1-a)}} {\infty} {1+a})
      \Lra \H^{a+b}(\Ynr,\nu^1_{\Ynr,\infty}). \]
Note that $E_1^{a,b}$ is of weight $b-1$
   so that $\H^i(\rf,E_1^{a,b})$ is finite
unless $b=1$.
Thus we obtain isomorphisms up to finite groups
\begin{equation}\label{UC.4.4}
\begin{CD}
   \H^1(\rf,\H^0(\Ynr,\nu^1_{\Ynr,\infty}))
     @.\; \simeq \; @.\H^1(\rf,E_2^{-1,1}),\\
  \qquad \H^1(\Ynr,\nu^1_{\Ynr,\infty})^{G_{\rf}} @. \simeq @. (E_2^{0,1})^{G_{\rf}} \qquad
\end{CD}
\end{equation}
with
$E_2^{-1,1}=\Ker(d_1^{-1,1})$ and $E_2^{0,1}=\Coker(d_1^{-1,1})$,
   where $d_1^{-1,1}$ is the Gysin map
\[ \H^0(\Ynr{}^{(2)},\qzp) \lra \H^1(\Ynr{}^{(1)},\logwitt {\Ynr{}^{(1)}} {\infty} 1 ). \]
There is an exact sequence of $G_{\rf}$-modules
   (cf.\ \eqref{UC.A.1} below)
\[ 0\lra \Pic(\Ynr{}^{(1)})\otimes\qzp  \lra
     \H^1(\Ynr{}^{(1)},\logwitt {\Ynr{}^{(1)}} {\infty} 1)
      \lra \Br(\Ynr{}^{(1)})_{\ptor} \lra 0. \]
Hence we see that
     the group \eqref{UC.4.3} coincides with the image of
     $\H^1(\rf,\H^0(\Ynr,\nu^1_{Y,\infty}))$ up to finite groups
   by the same computation
     as for Proposition \ref{claim5-1}\,(1)
     and the weight arguments in \cite{CTSS}, \S2.2.
Now we are reduced to showing
\[ \dim_{\qp}\H^2(k,\H^1(\Xggf,\qp(2))) \geq
    \corank\big(\H^1(\rf,\H^0(\Ynr,\nu^1_{\ol Y,\infty}))\big)
     =\corank\big(\H^1(\rf,E_2^{-1,1})\big), \]
where the last equality follows from \eqref{UC.4.4}.
As is seen in Remark \ref{rem5-2},
   the right hand side is equal to
    $\dim_{\Q}(\Delta \otimes \Q)$ under the condition {\bf T}.
On the other hand, by \cite{Ja2}, Corollary 7,
the left hand side does not change when one replaces $p$
    with another prime $p'$.
Thus the desired inequality follows from \eqref{UC.3.5}.
This completes the proof of Lemma \ref{lem5-1}
    and Theorem \ref{thm6-1}.
\qed
\subsection{Appendix to Section \ref{sect6}}\label{sect6.4}
Let $Z$ be a proper smooth variety over a finite field $\rf$.
For a positive integer $m$,
   we define the object $\Z/m\Z(1) \in D^b(Z_{\et},\Z/m\Z)$ as
$$
   \Z/m\Z(1) := \mu_{m'} \oplus
      (\logwitt Z r 1 [-1])
$$
where we factorized $m$ as $m' \cdot p^r$ with $(p,m')=1$.
There is a distinguished triangle
   of Kummer theory for $\Gm:=\Gm,_Z$ in $D^b(Z_{\et})$
$$
\begin{CD}
    \Z/m\Z(1) @>>> \Gm  @>{\times m}>> \Gm @>>> \Z/m\Z(1)[1].
\end{CD}
$$
So there is a short exact sequence
      of $G_{\rf}$-modules
$$
\begin{CD}
0  \lra  \Pic(\ol Z)/m  \lra \H^2(\ol Z,\Z/m\Z(1)) \lra
           {}_m\Br(\ol Z) \lra 0,
\end{CD}
$$
  where $\ol Z := Z \otimes_{\rf} \ol {\rf}$.
Taking the inductive limit with respect to $m \geq 1$,
      we obtain a short exact sequence
      of $G_{\rf}$-modules
\begin{equation}\label{UC.A.1}
\begin{CD}
0  \lra  \Pic(\ol Z)\otimes\qz  \os{\alpha}{\lra} \H^2(\ol Z,\qz(1))
    \lra \Br(\ol Z) \lra 0.
\end{CD}
\end{equation}
Concerning the arrow $\alpha$, we prove the following lemma,
   which has been used in this section.
\stepcounter{thm}
\begin{lem}\label{lem:app}
The map
$\H^1(\rf,\Pic(\ol Z)\otimes\qz) \to \H^1(\rf,\H^2(\ol Z,\qz(1))$ induced by $\alpha$ has finite kernel.
\end{lem}
\begin{pf}
Note that $\Pic(\ol Z)\otimes\qz \simeq
     (\NS(\ol Z)/\NS(\ol Z)_{\tor}) \otimes \qz$.
By a theorem of Matsusaka \cite{matsusaka}, Theorem 4,
         the group $\Div(\ol Z)/\Div(\ol Z)_{\num}$
            is isomorphic to $\NS(\ol Z)/\NS(\ol Z)_{\tor}$,
              where $\Div(\ol Z)$ denotes the group
                of Weil divisors on $\ol Z$,
                $\Div(Z)_{\num}$ denotes the subgroup
                    of Weil divisors numerically equivalent to zero.
By this fact and
     the fact that $\NS(\ol Z)$
       is finitely generated,
   there exists a finite family $\{C_i \}_{i \in I}$ of
     proper smooth curves over $\rf$
       which are finite over $Z$
      and for which the kernel of the natural map
        $\NS(\ol Z)
          \ra \bigoplus_{i\in I} \ \NS(\ol C_i)$ with
            $\ol C_i:=C_i \otimes_{\rf} \ol {\rf}$
            is torsion.
Now consider a commutative diagram
\[ \xymatrix{
 \H^1(\rf, \NS(\ol Z) \otimes \qz) \ar[r] \ar[d] & \H^1(\rf,\H^2(\ol Z, \qz(1))), \ar[d] \\
   \bigoplus{}_{i \in I} \ \H^1(\rf, \NS(\ol C_i) \otimes \qz) \ar[r]^{\sim \;\;} &
        \bigoplus{}_{i \in I}\ \H^1(\rf,\H^2(\ol C_i, \qz(1))). } \]
By a standard norm argument,
    one can easily show that the left vertical map has finite kernel.
The bottom horizontal arrow is bijective,
       because $\Br(\ol C_i)=0$ for any $i \in I$
         by Tsen's theorem (cf.\ \cite{Se}, II.3.3).
Hence the top horizontal arrow has finite kernel and we obtain
      the assertion.
\end{pf}
\newpage
\appendix
\section{Relation with conjectures of Beilinson and Lichtenbaum}\label{sectA}
\medskip
In this appendix, the Zariski site $Z_{\zar}$ on a scheme $Z$ always means
    $(\et/Z)_{\zar}$, and $Z_{\et}$ means the usual small \'etale site.
Let $k, p, S, \cX$ and $K$ be as in the notation \ref{intro}.8.
\subsection{Motivic complex and conjectures}
Let $\Z(2)_{\zar}=\Z(2)^{\cX}_{\zar}$ be the motivic complex on $\cX_{\zar}$
   defined by using Bloch's cycle complex, and let $\Z(2)_{\et}$ be its \'etale sheafification,
which are, by works of Levine (\cite{Le1}, \cite{Le2}), considered as strong candidates for motivic complexes
of Beilinson-Lichtenbaum (\cite{Be}, \cite{Li}) in Zariski and \'etale topology, respectively (see also \cite{Li2}, \cite{Li3}).
We put
\[ \H^*_{\zar}(\cX\hspace{-2pt},\bZ(2)):=\H^*_{\zar}(\cX\hspace{-2pt},\bZ(2)_{\zar}), \qquad
\H^*_{\et}(\cX\hspace{-2pt},\bZ(2)):=\H^*_{\et}(\cX\hspace{-2pt},\bZ(2)_{\et}). \]
In this appendix, we observe that the finiteness of $\H_{\ur}^3(K,\qzp(2))$ is deduced from the following
    conjectures on motivic complexes:
\begin{conj}\label{conjA-1}
Let $\epsilon:\cX_{\et}\to \cX_{\zar}$
   be the natural continuous map of sites.
Then$:$
\begin{enumerate}
\item[(1)]
$($Beilinson-Lichtenbaum conjecture$)$.
We have
\[ \begin{CD} \bZ(2)_{\zar} \isom \tau_{\leq 2} R\epsilon_* \bZ(2)_{\et} \quad \hbox{ in } D(\cX_{\zar}).\end{CD} \]
\item[(2)]
$($Hilbert's theorem $90)$.
We have $R^3\epsilon_* \bZ(2)_{\et}=0$.
\item[(3)]
$($Kummer theory on $\cX[p^{-1}]_{\et})$.
We have $(\bZ(2)_{\et})|_{\cX[p^{-1}]} \otimes^{\L} \Z/p^r\Z \simeq \mu_{p^r}^{\otimes 2}$.
\end{enumerate}
\end{conj}
\noindent
This conjecture holds if $\cX$ is smooth over $S$ by a result of Geisser \cite{Ge}, Theorem 1.2
    and the Merkur'ev-Suslin theorem \cite{MS} (see also \cite{GL2}, Remark 5.9).
\begin{conj}\label{conjA-2}
Let $\gamma^2$ be the canonical map
\[ \gamma^2 : \CH^2(\cX)=\H^4_{\zar}(\cX\hspace{-2pt},\bZ(2)) \lra
\H^4_{\et}(\cX\hspace{-2pt},\bZ(2)). \]
Then the $p$-primary torsion part of $\Coker(\gamma^2)$ is finite.
\end{conj}
\noindent
This conjecture is based on Lichtenbaum's conjecture \cite{Li}
   that $\H^4_{\et}(\cX\hspace{-2pt},\bZ(2))$ is a finitely generated abelian group (by the properness of $\cX/S$).
The aim of this appendix is to prove the following:
\begin{prop}\label{propA-1}
If Conjectures $\mathrm{\ref{conjA-1}}$ and $\mathrm{\ref{conjA-2}}$ hold, then $\H^3_{\ur}(K,\qzp(2))$ is finite.
\end{prop}
\noindent
This proposition is reduced to the following lemma:
\begin{lem}\label{lemA-1}
\begin{enumerate}
\item[(1)]
If Conjecture $\mathrm{\ref{conjA-1}}$ holds, then for $r \geq 1$ there is an exact sequence
\[  0 \to \Coker\left({}_{p^r}\CH^2(\cX) \os{\alpha_r}{\to} {}_{p^r}\H^4_{\et}(\cX\hspace{-2pt},\bZ(2))\right) \to
\H^3_{\ur}(K,\prz(2)) \to \Ker(\varrho^2_{r}) \to 0, \]
where $\alpha_r$ denotes the map induced by $\gamma^2$ and $\varrho^2_{r}$ denotes the cycle class map
\[ \varrho^2_{r}: \CH^2(\cX)/p^r \lra \H^{4}_{\et}(\cX\hspace{-2pt},\pnuX). \]
\item[(2)]
If Conjectures $\mathrm{\ref{conjA-1}}$ and $\mathrm{\ref{conjA-2}}$ hold,
then $\Coker(\alpha_{\qzp})$ and $\Ker(\varrho^2_{\qzp})$ are finite,
   where $\alpha_{\qzp}:= \varinjlim{}_{r \ge 1} \ \alpha_r$ and $\varrho^2_{\qzp}:=\varinjlim{}_{r \ge 1} \ \varrho^2_r$.
\end{enumerate}
\end{lem}
\noindent
To prove this lemma, we need the following sublemma,
    which is a variant of Geisser's arguments in \cite{Ge}, \S6:
\begin{sublem}\label{lemA-2}
Put $\prz(2)_{\et} :=\Z(2)_{\et} \otimes^{\L}\Z/p^r\Z$.
If Conjecture $\mathrm{\ref{conjA-1}}$ holds, then there is a unique isomorphism
\[ \prz(2)_{\et} \isom \pnuX \quad \hbox{ in } D(\cX_{\et},\Z/p^r\Z) \]
that extends the isomorphism in Conjecture $\mathrm{\ref{conjA-1}}\,(3)$.
\end{sublem}
\noindent
We prove Sublemma \ref{lemA-2} in \S\ref{appA-2} below and Lemma \ref{lemA-1} in \S\ref{appA-3} below.
\subsection{Proof of Sublemma \ref{lemA-2}}\label{appA-2}
By Conjecture \ref{conjA-1}\,(3), we have only to consider the case where $p$ is not invertible on $S$.
Let us note that
\par
\vspace{3pt}
\begin{quote}
$(*)$ \quad $\Z/p^r\Z(2)_{\et}$ is concentrated in degrees $\leq 2$
\end{quote}
\par
\vspace{3pt}
\noindent
by Conjecture \ref{conjA-1}\,(1) and (2).
Let $V$, $Y$, $i$ and $j$ be as follows:
\[ \begin{CD} V:=\cX[p^{-1}] @>{j}>> \cX @<{i}<< Y, \end{CD} \]
where $Y$ denotes the union of the fibers of $\cX/S$ of characteristic $p$.
In \'etale topology, we define $Ri^!$ and $Rj_*$ for unbounded
    complexes by the method of Spaltenstein \cite{Spa}.
We will prove
\begin{equation}\label{eq:Y}
      \tau_{\leq 3}Ri^! \Z/p^r\Z(2)_{\et} \simeq \nu_{Y,r}^1 [-3]
        \quad \hbox{ in } D(Y_{\et},\Z/p^r\Z),
\end{equation}
using $(*)$
    (see (S5) in \S\ref{sect3.1} for $\nu_{Y,r}^1$).
We first prove Sublemma \ref{lemA-2} admitting this isomorphism.
Since $(\Z(2)_{\et})|_V \otimes^{\L}\Z/p^r\Z \simeq \mu_{p^r}^{\otimes 2}$
     by Conjecture \ref{conjA-1}\,(3),
   we obtain a distinguished triangle from \eqref{eq:Y} and $(*)$
\[ \begin{CD}  i_* \nu_{Y,r}^1 [-3] @>>> \Z/p^r\Z(2)_{\et} @>>>
         \tau_{\leq 2} Rj_* \mu_{p^r}^{\otimes 2} @>>> i_* \nu_{Y,r}^1 [-2]. \end{CD} \]
Hence comparing this distinguished triangle with that of (S5) in \S\ref{sect3.1},
 we obtain the desired isomorphism in the sublemma, whose uniqueness follows from \cite{Sat2},\,Lemmas\,1.1 and 1.2\,(1).

In what follows, we prove \eqref{eq:Y}.
Put $\K:=\Z(2)_{\zar} \otimes ^{\L} \Z/p^r\Z$
     and $\LL:= \Z/p^r\Z(2)_{\et}$ for simplicity.
Let $\epsilon:\cX_{\et}\to \cX_{\zar}$ be as in Conjecture \ref{conjA-1}.
In Zariski topology, we define $Ri_{\zar}^!$ and $Rj_{\zar *}$ for unbounded
    complexes in the usual way by the finiteness of cohomological dimension.
Because $\LL=\epsilon^*\K$
   is concentrated in degrees $\leq 2$ by $(*)$,
     there is a commutative diagram
      with distinguished rows in $D^b(\cX_{\et},\Z/p^r\Z)$
\[\xymatrix{
 \e^* \K \ar[r] \ar@{=}[d] & \tau_{\leq 2} \e^*Rj_{\zar*}j_{\zar}^*\K \ar[r] \ar[d]_{\alpha}
    & (\tau_{\leq 3}\e^* i_{\zar*}Ri_{\zar}^!\K)[1] \ar[r] \ar[d]_{\beta} & \e^* \K[1] \ar@{=}[d] \\
   \LL \ar[r] & \tau_{\leq 2} Rj_{\et*}j_{\et}^*\LL \ar[r]
    & (\tau_{\leq 3}i_{\et*}Ri_{\et}^!\LL)[1] \ar[r] & \LL[1], } \]
where the upper (resp.\ lower) row is obtained from the localization triangle
    in the Zariski (resp.\ \'etale) topology and the arrows
       $\alpha$ and $\beta$ are canonical base-change morphisms.
Since $\alpha$ is an isomorphism (\cite{MS}, \cite{SV}, \cite{GL2}), $\beta$ is an isomorphism as well.
Hence \eqref{eq:Y} is reduced to showing
\begin{equation}\label{eq:Y2}
\tau_{\leq 3} Ri_{\zar}^!\K \simeq \e_{Y*} \nu_{Y,r}^1[-3]
     \quad \hbox{ in } D(Y_{\zar},\Z/p^r\Z),
\end{equation}
where $\e_Y : Y_{\et} \ra Y_{\zar}$ denotes the natural continuous map of sites and we have used the base-change
   isomorphism $\e^*i_{\zar*} = i_{\et*}\e_Y^*$ (\cite{Ge}, Proposition 2.2\,(a)).
Finally we show \eqref{eq:Y2}. Consider the local-global spectral sequence in the Zariski topology
\[ E_1^{u,v}= \bigoplus_{x \in \cX^u \cap Y} \ R^{u+v}i_{x*}(Ri_x^!Ri_{\zar}^!\K) \Lra R^{u+v}i_{\zar}^! \K, \]
where for $x \in Y$, $i_x$ denotes the natural map $x \to Y$.
We have
\[ E_1^{u,v} \simeq \begin{cases} \displaystyle \bigoplus{}_{x\in \cX^u \cap Y}\ i_{x*}\e_{x*}\logwitt x r {2-u}
    \quad&\text{ (if } v = 2)\\ 0 \quad&\text{ (otherwise)} \end{cases} \]
by the localization sequence of higher Chow groups \cite{Le1} and results of Geisser-Levine (\cite{GL1}, Proposition 3.1,
 Theorem 7.1), where for $x \in Y$, $\e_x$ denotes the natural continuous map $x_{\et} \to x_{\zar}$ of sites.
By this description of $E_1$-terms and the compatibility of boundary maps (\cite{GL2}, Lemma 3.2,
     see also \cite{Sz}, Appendix), we obtain \eqref{eq:Y2}.
This completes the proof of Sublemma \ref{lemA-2}.
\subsection{Proof of Lemma \ref{lemA-1}}\label{appA-3}
(1) By Sublemma \ref{lemA-2}, there is an exact sequence
\[ 0 \lra \H^3_{\et}(\cX\hspace{-2pt},\bZ(2))/p^r \lra \H^3_{\et}(\cX\hspace{-2pt},\pnuX) \lra
    {}_{p^r}\H^4_{\et}(\cX\hspace{-2pt},\bZ(2)) \lra 0. \]
By Conjecture \ref{conjA-1}\,(1) and (2), we have
\[ \H^3_{\et}(\cX\hspace{-2pt},\bZ(2))
\simeq \H^3_{\zar}(\cX\hspace{-2pt},\bZ(2)) \simeq \CH^2(\cX\hspace{-2pt},1). \]
Thus we get an exact sequence
\[ 0 \lra \CH^2(\cX\hspace{-2pt},1)/p^r \lra \H^3_{\et}(\cX\hspace{-2pt},\pnuX) \lra
     {}_{p^r}\H^4_{\et}(\cX\hspace{-2pt},\bZ(2)) \lra 0. \]
On the other hand, there is an exact sequence
\[ 0 \lra \Naa\H^3_{\et}(\cX\hspace{-2pt},\pnuX) \lra \H^3_{\et}(\cX\hspace{-2pt},\pnuX) \lra
     \H^3_{\ur}(K,\prz(2)) \lra  \Ker(\varrho^2_{r})\lra 0 \]
which is a variant of \eqref{UC.0.1} (see Lemma \ref{lem3-1} for $\Naa$).
In view of the short exact sequence in Lemma \ref{lem3-1}\,(1),
we get the desired exact sequence.
\par
(2) By Conjecture \ref{conjA-1}\,(1) and (2),
   the map $\gamma^2$ in Conjecture \ref{conjA-2} is injective.
Hence we get an exact sequence
\[ 0 \lra \Coker(\alpha_r) \lra {}_{p^r}\Coker(\gamma^2) \lra \CH^2(\cX)/p^r \os{\gamma^2/p^r}{\lra}
    \H^4_{\et}(\cX\hspace{-2pt},\bZ(2))/p^r. \]
Noting that the composite of $\gamma^2/p^r$ and the injective map
\[ \xymatrix{
 \H^4_{\et}(\cX\hspace{-2pt},\bZ(2))/p^r \; \ar@{^{(}->}[r] & \H^{4}_{\et}(\cX\hspace{-2pt},\pnuX)
}\]
obtained from Sublemma \ref{lemA-2} coincides with $\varrho^2_{r}$, we get a short exact sequence
\[ 0 \lra \Coker(\alpha_r) \lra
{}_{p^r}\Coker(\gamma^2) \lra \Ker(\varrho^2_r) \lra 0, \]
    which implies the finiteness of $\Coker(\alpha_{\qzp})$
      and $\Ker(\varrho^2_{\qzp})$
      under Conjecture \ref{conjA-2}.
This completes the proof of Lemma \ref{lemA-1}
    and Proposition \ref{propA-1}.
\qed
\newpage
\section{Zeta value of threefolds over finite fields}\label{appB}
In this appendix B, all cohomology groups of schemes are taken over the \'etale topology.
Let $X$ be a projective smooth geometrically integral threefold over a finite field $\bF_q$,
 and let $K$ be the function field of $X$.
We define the unramified cohomology $\H^{n+1}_{\ur}(K,\qz(n))$ in the same way as in $\ref{intro}.8$.
We show that the groups
\[ \HurXb =\Br(X) \quad \hbox{ and } \quad \HurXc \]
are related with the value of the Hasse-Weil zeta function $\zeta(X,s)$ at $s=2$:
\[ \zeta^* (X,2):=\lim_{s \to 2} \zeta(X,s)(1-q^{2-s})^{-\varrho_2}, \qquad
  \hbox{where } \;\varrho_2:=\ord_{s=2}\;\zeta(X,s). \]
Let
\[ \theta \;:\; \CH^2(X) \lra \Hom(\CH^1(X),\bZ) \]
be the map induced by the intersection pairing and the degree map
\[ \CH^2(X) \times \CH^1(X) \lra \CH^3(X)=\CH_0(X) \os{\deg}\lra \bZ. \]
The map $\theta$ has finite cokernel by a theorem of Matsusaka \cite{matsusaka}, Theorem 4.
We define \[ \cR:=|\Coker(\theta)|. \]
We prove the following formula (compare with the formula in \cite{Ge2}):
\begin{mthm}\label{thmB-1}
Assume that $\Br(X)$ and $\HurXc$ are finite.
Then $\zeta^* (X,2)$ equals the following rational number
    up to a sign$:$
\[ q^{\chi(X,\cO_X,2)} \cdot \frac{|\HurXc|}{|\Br(X)| \cdot \cR} \cdot
\prod_{i=0}^3 |\CH^2(X,i)_{\tor}|^{(-1)^i} \cdot
\prod_{i=0}^1 |\CH^1(X,i)_{\tor}|^{(-1)^i}, \]
where $\CH^2(X,i)$ and $\CH^1(X,i)$ denote Bloch's higher Chow groups \cite{Bl2} and
 $\chi(X,\cO_X,2)$ denotes the following integer$:$
\[ \chi(X,\cO_X,2) := \sum_{i,j}(-1)^{i+j}\,(2-i)\,\dim_{\bF_q}\H^j(X,\Omega_X^i) \quad (0 \le i \le 2,\;\;0 \le j \le 3). \]
\end{mthm}
This theorem follows from a theorem of Milne (\cite{Mi1}, Theorem 0.1) and Proposition \ref{propB-1} below.
For integers $i , n \ge 0$, we define
\[ \H^i(X,\bhZ(n)) := \prod_{\text{all \,}\ell} \ \H^i(X,\bZ_{\ell}(n)), \]
where $\ell$ runs through all prime numbers, and
 $\H^i(X,\bZ_p(n))$\,$(p:=\ch(\bF_q))$ is defined as
\[ \H^i(X,\bZ_p(n)) := \varprojlim_{r \ge 1} \ \H^{i-n}(X,\logwitt X r n). \]
\begin{mprop}\label{propB-1}
\begin{itemize}
\item[(1)]
We have
\begin{equation}\label{eq:svz1}
\tag{B.3} \H^i(X,\bhZ(2)) \simeq
   \begin{cases}
    \CH^2(X,4-i)_{\tor} \quad & (i=0,1,2,3) \\
    (\CH^1(X,i-6)_{\tor})^*  \quad & (i = 6,7),
   \end{cases}
\end{equation}
where for an abelian group $M$, we put \[ M^* := \Hom(M,\qz). \]
Furthermore, $\CH^1(X,j)_{\tor}$ and $\CH^2(X,j)_{\tor}$ are finite for any $j \ge 0$, and we have
\[ \CH^1(X,j)_{\tor}=0 \; \hbox{ for } \; j \ge 2 \quad \hbox{ and }\quad \CH^2(X,j)_{\tor}=0 \; \hbox{ for } \; j \ge 4. \]
\item[(2)]
Assume that $\Br(X)$ is finite. Then we have
\[ \H^5(X,\bhZ(2))_{\tor} \simeq \Br(X)^*, \]
and the cycle class map
\[ \CH^2(X)\otimes\bZ_{\ell} \lra \H^4(X,\bZ_{\ell}(2)) \]
has finite cokernel for any prime number $\ell$.
\item[(3)]
Assume that $\Br(X)$ and $\HurXc$ are finite.
Then the following map given by the cup product with the canonical element $1 \in \bhZ \simeq \H^1(\bF_q,\bhZ)$
 has finite kernel and cokernel{\rm :}
\[ \epsilon^4 : \H^4(X,\bhZ(2)) \lra \H^5(X,\bhZ(2)), \]
and we have the following equality of rational numbers$:$
\[ \frac{|\ker(\epsilon^4)|}{|\Coker(\epsilon^4)|} = \frac{|\HurXc| \cdot |\CH^2(X)_{\tor}|}{|\Br(X)| \cdot \cR}. \]
\end{itemize}
\end{mprop}
\begin{pf*}{\it Proof of Proposition \ref{propB-1}}
(1) By standard arguments on limits, there is a long exact sequence
\[ \begin{CD}
   \dotsb @.~ \lra~ @. \H^i(X,\bhZ(2))
      @.~ \lra ~@. \H^i(X,\bhZ(2)) \otimes_{\bZ} \bQ
      @.~ \lra ~@. \H^i(X,\qz(2))\\
      @. \lra @. \H^{i+1}(X,\bhZ(2))
      @. \lra @. \hspace{-67pt}\dotsb.
\end{CD} \]
By \cite{CTSS}, p.\ 780, Th\'eor\`eme 2, p.\ 782, Th\'eor\`eme 3, we see that
\begin{center}
 $\H^i(X,\bhZ(2))$ and $\H^i(X,\qz(2))$ are finite for $i \ne 4,5$.
\end{center}
Hence we have
\[ \H^i(X,\bhZ(2))\simeq \H^{i-1}(X,\qz(2)) \quad\text{for } \; i \ne 4,5,6. \]
On the other hand, there is an exact sequence
\[ 0 \lra \CH^2(X,5-i) \otimes \qz \lra \H^{i-1}(X,\qz(2)) \lra \CH^2(X,4-i)_{\tor} \lra 0 \]
for $i \le 3$ (\cite{MS}, \cite{SV}, \cite{GL1}, \cite{GL2}),
where $\CH^2(X,5-i) \otimes \qz$ must be zero because it is divisible and finite.
Thus we get the isomorphism \eqref{eq:svz1} for $i \le 3$,
    the finiteness of $\CH^2(X,j)_{\tor}$ for $j \ge 1$ and the vanishing of $\CH^2(X,j)_{\tor}$ for $j \ge 4$.
The finiteness of $\CH^2(X,0)_{\tor}= \CH^2(X)_{\tor}$ (cf.\ \cite{CTSS}, p.\ 780, Th\'eor\`eme 1)
follows from the exact sequence
\[ 0 \lra \CH^2(X,1) \otimes \qz \lra \Naa\H^3(X,\qz(2)) \lra \CH^2(X)_{\tor} \lra 0 \]
(cf.\ Lemma \ref{lem2-0}), where we put
\[ \Naa\H^3(X,\qz(2)):=\ker(\H^3(X,\qz(2)) \to \H^3(K,\qz(2))). \]
As for the case $i=6,7$ of \eqref{eq:svz1}, we have
\[ \H^i(X,\bhZ(2))^* \simeq \H^{7-i}(X,\qz(1)) \]
 by a theorem of Milne \cite{Mi1}, Theorem 1.14\,(a).
It remains to show
\[ \CH^1(X,j)_{\tor} \simeq \H^{1-j}(X,\qz(1)) \quad \hbox { for } \; j \ge 0, \]
which can be checked by similar arguments as before.
\par
(2) We have $\H^5(X,\bhZ(2))^* \simeq \H^2(X,\qz(1))$ and an exact sequence
\begin{equation}\label{eqB-2}\tag{B.4}
 0 \lra \CH^1(X)\otimes\qz \lra \H^2(X,\qz(1)) \lra \Br(X) \lra 0. \end{equation}
Hence we have $(\H^5(X,\bhZ(2))_{\tor})^* \simeq \Br(X)$, assuming $\Br(X)$ is finite.
To show the second assertion for $\ell \ne \ch(\bF_q)$, it is enough to show that the cycle class map
\def\crys{{\mathrm{crys}}}
\[ \CH^2(X)\otimes\ql \lra \H^4(\ol{X},\ql(2))^{\vG}\]
is surjective, where $\vG:=\Gal(\ol {\bF_q}/\bF_q)$.
The assumption on $\Br(X)$ implies the bijectivity of the cycle class map
\[ \CH^1(X)\otimes\ql \isom \H^2(\ol{X},\ql(1))^{\vG} \]
by \cite{Ta2}, (4.3) Proposition (see also \cite{Mi0}, Theorem 4.1),
 and the assertion follow from \cite{Ta2}, (5.1) Proposition.
As for the case $\ell=\ch(\bF_q)$,
 one can easily pursue an analogy using crystalline cohomology, whose details are left to the reader.
\par
(3) 
The finiteness assumption on $\Br(X)$ implies the condition {\bf SS}$(X,1,\ell)$ in \cite{Mi1} for all prime numbers $\ell$
 by loc.\ cit., Proposition 0.3.
Hence {\bf SS}$(X,2,\ell)$ holds by the Poincar\'e duality,
 and $\epsilon^4$ has finite kernel and cokernel by loc.\ cit., Theorem 0.1.
\par
To show the equality assertion, we put
\[ \wh{\CH}{}^2(X):=\varprojlim{}_{n \ge 1} \ \CH^2(X)/n, \]
and consider the following commutative square (cf.\ \cite{Mi2}, Lemma 5.4):
\[ \xymatrix{
\wh{\CH}{}^2(X) \ar[rr]^{\varTheta \qquad} \ar[d]_{\alpha} && \Hom(\CH^1(X),\bhZ) \\
\H^4(X,\bhZ(2)) \ar[rr]^{\epsilon^4} && \H^5(X,\bhZ(2)) \ar[u]_{\beta},  } \]
where the top arrow $\varTheta$ denotes the map induced by $\theta$.
The arrow $\alpha$ denotes the cycle class map of codimension $2$, and $\beta$ denotes
 the Pontryagin dual of the cycle class map with $\qz$-coefficients in \eqref{eqB-2}.
The arrow $\alpha$ is injective  (cf.\ \eqref{UC.0.1}) and we have
\[ |\Coker(\alpha)|=|\HurXc| \]
by the finiteness assumption on $\HurXc$ and (2) (cf.\ Proposition \ref{prop5-0}).
The arrow $\beta$ is surjective and we have
\[ \ker(\beta)=\H^5(X,\bhZ(2))_{\tor} \os{(2)}{\simeq} \Br(X)^*, \]
by Milne's lemma (\cite{Mi2}, Lemma 5.3) and the isomorphism $\CH^1(X)\otimes\bhZ \simeq \H^2(X,\bhZ(1))$
 (cf.\ \cite{Ta2}, (4.3) Proposition), where we have used again the finiteness assumption on $\Br(X)$.
Therefore in view of the finiteness of $\ker(\epsilon^4)$, the map $\varTheta$ has finite kernel and we obtain
\[ \ker(\varTheta) = \wh{\CH}{}^2(X)_{\tor}=\CH^2(X)_{\tor}, \]
where we have used the finiteness of $\CH^2(X)_{\tor}$ in (1).
Finally the assertion follows from
 the following equality concerning the above diagram:
\[ \frac{|\ker(\varTheta)|}{|\Coker(\varTheta)|} =\frac{|\ker(\alpha)|}{|\Coker(\alpha)|}\cdot
 \frac{|\ker(\epsilon^4)|}{|\Coker(\epsilon^4)|} \cdot \frac{|\ker(\beta)|}{|\Coker(\beta)|} \]
This completes the proof of Proposition \ref{propB-1} and Theorem \ref{thmB-1}.
\end{pf*}

\newpage 

\end{document}